\documentclass[10pt]{article}
\usepackage{graphicx}
\usepackage{amsmath,amssymb,amsthm,amsfonts}
\usepackage{color}
\usepackage{amssymb}
\usepackage{cite}
\newtheorem{thm}{Theorem}[section]

\newtheorem{lem}[thm]{Lemma}
\newtheorem{prop}[thm]{Proposition}

\theoremstyle{definition}
\newtheorem{defn}{Definition}[section]
\theoremstyle{remark}
\usepackage{appendix}
\newtheorem{rem}{Remark}[section]

\numberwithin{equation}{section}

\DeclareMathSymbol{\C}{\mathalpha}{AMSb}{"43}

\newcommand{\eps}{\varepsilon}

\newcommand{\lam}{\lambda}

\newcommand{\R}{{\mathbb{R}}}
\newcommand{\h}{{\mathcal{H}}}
\newcommand{\inte}{\int_{\mathbb{R}^2}}

\def\a{\alpha_a}

\def\e{\varepsilon_a}

\def\R{{\mathbb R}}
\def\C{{\mathbb C}}

\newcommand{\bsub}{\begin{subequations}}
\newcommand{\esub}{\end{subequations}$\!$}

\begin{document}
\title{The Nonexistence of Vortices for Rotating Bose-Einstein Condensates with Attractive Interactions}
	
\author{
Yujin Guo\thanks{School of Mathematics and Statistics, and Hubei Key Laboratory of
Mathematical Sciences, Central China Normal University, Wuhan, 430079, P. R. China.  Emails: \texttt{yjguo@wipm.ac.cn}; \texttt{yguo@ccnu.edu.cn}. Y. J.  Guo is partially supported by NSFC under Grants No. 11671394 and 11931012.
		},
\, Yong Luo\thanks{University of Chinese Academy of Sciences, Beijing 100190, P. R. China;  Wuhan Institute of Physics and Mathematics, Chinese Academy of Sciences, P.O. Box 71010, Wuhan 430071, P. R. China.  Email: \texttt{luoyong.wipm@outlook.com}. Y. Luo is partially supported by the Project funded by China Postdoctoral Science Foundation No. 2019M662680.
 }
\, and\, Wen Yang\thanks{Wuhan Institute of Physics and Mathematics,
Chinese Academy of Sciences, P.O. Box 71010, Wuhan 430071, P. R. China; Innovation Academy for Precision Measurement Science and Technology, Chinese Academy of Sciences, Wuhan 430071, P. R. China.  Email: \texttt{wyang@wipm.ac.cn}. W. Yang is partially supported by NSFC under Grant No. 11801550.}
	}
	
\date{\today}
	
	\smallbreak \maketitle

\begin{abstract} This article is devoted to {studying} the model of two-dimensional attractive Bose-Einstein condensates in a trap $V(x)$  rotating at the velocity $\Omega $. {This model} can be described by the complex-valued  Gross-Pitaevskii energy functional. It is shown that there exists a critical rotational velocity $0<\Omega ^*:=\Omega^*(V)\leq \infty$, depending on the general trap $V(x)$, {such that} for any rotational velocity $0\le \Omega <\Omega ^*$, minimizers (i.e., ground states) exist if and only if $ a<a^*=\|w\|^2_2$, where $a>0$ denotes the absolute product for the number of particles times the scattering length, and $w>0$ is the unique positive solution  of {$\Delta w-w+w^3=0$} in $\R^2$. If $V(x)=|x|^2$ and $ 0<\Omega <\Omega ^*(=2)$ is fixed, we prove that, up to a constant phase,  all minimizers must be real-valued, unique and {free of} vortices as $a \nearrow a^*$, by analyzing the refined limit behavior of minimizers and employing the non-degenerancy of $w$.
	
\end{abstract}
	
\vskip 0.2truein
	
Key words: Bose-Einstein condensates; rotational velocity; Gross-Pitaevskii functional; mass concentration
	
\vskip 0.2truein

\tableofcontents

\section{Introduction}

\indent In physical experiments of rotating Bose-Einstein condensates (BECs), a large number of dilute (bosonic) atoms are confined in rotating traps and cooled to the sufficiently low temperatures (cf. \cite{Abo,A,B,CR,D,F,Lewin}). The condensation of a large fraction of particles into the same one-particle state {can be} observed, {when the temperature of the system decreases to a critical value.} These BECs in rotating traps display various interesting quantum phenomena, such as the appearance of quantized vortices \cite{Abo,B,CC,F,WG,ZW}, the center-of-mass  rotation \cite{Abo,LC,F,WG}, the effective lower dimensional behavior in strongly elongated traps \cite{A,Anderson,BC,CR,D,F}, and so on. Therefore, the properties and their applications of rotating BECs have been a focus of international interest in physics and mathematics since the first physical achievement of rotating BECs in the late 1990s. Especially, the complex structures, including the center-of-mass motions and the quantized vortices of all kinds, of rotating BECs with repulsive interactions have been observed and analyzed extensively over the past two decades, see \cite{Abo,A,CD,CP,CR,CRY,D,F} and the references therein.

In contrast to the repulsive case, when the forces between the cold atoms in the condensates are attractive,  {the system  then collapses as the number of the cold atoms increases beyond a critical value, see \cite{Hulet1,Hulet2,HM,KM, Hulet3} or \cite[Sec.~III.B]{D}}. Because of the above distinct mechanism, the system of rotating BECs in the attractive case presents more {complicated} phenomena and structures, some of which were recently explored by theoretical analysis and numerical simulations, see \cite{CC,D,F,LC} and the references therein. For example, the collapse cycles and implosion were investigated in three dimensional attractive BECs under rotation.  On the other hand, the bright soliton propagation and interactions were however observed in one dimensional attractive BECs under rotation. The main goal of this paper is to investigate mathematically the ground states of two dimensional attractive BECs in rotating traps, which were analyzed recently in \cite{BC,GS,Lewin,ANS} and the references therein.

As derived rigorously in \cite{Lewin}, the ground states of two-dimensional attractive BECs in  rotating traps can be equivalently described (see also \cite{BC,ANS,P}) by the minimizers of the following complex-valued  variational problem:
\begin{equation}\label{def:ea}
	e_F(a):=\inf _{\{u\in \h, \, \|u\|^2_2=1 \} } F_a(u) \,,
\end{equation}
where the Gross-Pitaevskii (GP) functional $F_a(u)$ is given by
\begin{equation}
	F_a(u):=\int _{\R ^2} \big(|\nabla u|^2+V(x)|u|^2\big)dx-\frac{a}{2}\int _{\R ^2}|u|^4dx-
	\Omega \int_{\R ^2}x^{\perp}\cdot (iu,\, \nabla u)dx,  \ \ u\in \h , \label{f}
\end{equation}
and the space $\h$ is defined as
\begin{equation}
	\h :=  \Big \{u\in  H^1(\R ^2, \mathbb{C}):\ \int _{\R ^2}
	V(x)|u|^2 dx<\infty\Big \}. \label{1.H}
\end{equation}
Here $x^{\perp} =(-x_2,x_1)$ with $x=(x_1,x_2)\in \R^2$, and
	$(iu,\, \nabla u)=i(u\nabla \bar u-\bar u\nabla u)/2$.
The parameter $ a>0$ in (\ref{f}) characterizes the absolute product of the scattering length $\nu$ of the two-body interaction times the number $N$ of particles in the condensates, and $\Omega \ge 0$ describes the rotational velocity of the trap $V(x)$. Alternatively, one may impose the constraint $\int_{\mathbb{R}^2} |u(x)|^2dx=N>0$, but the latter case can be easily  reduced to the previous one with $a$ being replaced by $a/N$. Therefore, in this paper we {shall} work with $e_F(a)$ instead. Our main interest of investigating $e_F(a)$ is {on two folds: on one hand, we provide an accurate limit description of minimizers for $e_F(a)$ with more general potentials $V(x)$; on the other hand, more importantly we shall prove analytically that up to a constant phase, all minimizers of $e_F(a)$ are real valued, unique and free of vortices for some special case.


When there is no rotation for the trap, $i.e.$ $\Omega = 0$, the existence, uniqueness, symmetric breaking and the refined limit behavior of real-valued constraint minimizers for $e_F(a)$ were studied recently in \cite{GLW,GS,GWZZ,GZZ} and the references therein, see also Theorem \ref{thmA} below for some details. Moreover, the stability and instability of real-valued constraint minimizers for $e_F(a)$ with $\Omega = 0$ were also studied in \cite{BC,Zhang,Z} and the references therein.
For the above non-rotational case, it turns out that the analysis  of $e_F(a)$ is {connected well} with the following nonlinear
scalar field equation
\begin{equation}
	-\Delta u+ u-u^3=0\  \mbox{  in } \  \R^2,\  \mbox{ where }\ u\in H^1(\R ^2,\R).  \label{Kwong}
\end{equation}
Recall from \cite{K,W} that, up to translations,  (\ref{Kwong}) admits a unique positive solution, which is radially symmetric. For convenience, {we denote the unique positive solution of (\ref{Kwong}) by $w=w(|x|)>0$.} Note also from \cite[Lemma 8.1.2]{C} that $w=w(|x|)>0$ satisfies
\begin{equation}\label{1:id}
\inte |\nabla w |^2dx  =\inte w ^2dx=\frac{1}{2}\inte w ^4dx.
\end{equation}
Moreover, it was proved in \cite{W} that the equal sign of the following
Gagliardo-Nirenberg inequality
\begin{equation}\label{GNineq}
	\inte |u(x)|^4 dx\le \frac 2 {\|w\|_2^{2}} \inte |\nabla u(x) |^2dx \inte |u(x)|^2dx ,\
\  u \in H^1(\R ^2, \R)\,
\end{equation}
is achieved for $u(x) = w(x)$.

More recently, the variational problem $e_F(a)$ under rotation $\Omega >0$ was discussed in \cite{BC,ANS,Lewin}, where the authors addressed the existence, non-existence and the limit behavior of complex-valued minimizers for $e_F(a)$ mainly in the case where $V(x)=|x|^2$. Recall from \cite{Lieb} the following diamagnetic
inequality:
for $\mathcal{A }=\frac{\Omega}{2}x^{\perp}$,
\begin{equation}  |\nabla u|^2 -\Omega \, x^{\perp}\cdot (iu,\, \nabla u) =
	|(\nabla -i\mathcal{A } )u|^2 -\frac{\Omega ^2}{4}  |x|^2|u|^2\ge \big| \nabla |u|\big|^2
	-\frac{\Omega ^2}{4}  |x|^2|u|^2, \quad u\in H^1(\R^2, \mathbb{C}) .
	\label{I:2:1A}
\end{equation}
By making full use of the inequality (\ref{I:2:1A}), one can note that the existence and non-existence results of \cite{Lewin} can be extended to the variational problem $e_F(a)$ with more general trapping potentials $V(x)$. {Concerning on this point}, we consider the general trapping potential  $0\le V(x)\in L^\infty_{loc}(\R^2)$  satisfying
\begin{equation}
	\underline{\lim} _{|x|\to\infty }\frac{V(x)}{|x|^2}>0.
\label{A:V}
\end{equation}
We also define the critical rotational velocity $\Omega ^*$ by
\begin{equation}
	\Omega ^*:=\sup \Big\{\Omega >0:\ \  V(x)-\frac{\Omega ^2}{4}|x|^2  \to\infty \,\ \mbox{as}\,\
	|x|\to\infty \Big\}.  \label{Omega}
\end{equation}
As an example, suppose $V(x)=|x|^s$ with $s\ge 2$, then we have
\begin{equation}
	\Omega ^* :=\arraycolsep=1.5pt \left\{ \begin{array}{ll}
2,\quad &  {\rm if}\quad s=2;\\ [2mm]
	\infty, \ &  {\rm if}\quad s>2,\end{array}\right.
\label{omega-1}
\end{equation}
which shows that depending on $V(x)$, both $0<\Omega ^*<\infty$ and $\Omega ^*=\infty$ can happen. Applying the inequalities
(\ref{GNineq}) and (\ref{I:2:1A}), the arguments of \cite{GS,Lewin} give essentially the following existence and nonexistence under the general assumption $(\ref{A:V})$:

\begin{thm}\label{thm1} Let $w$ be the unique positive solution of (\ref{Kwong}). Suppose $V(x)\in L^\infty_{\rm loc}(\R^2)$ satisfies
(\ref{A:V}). Then we have
\begin{enumerate}
\item If $0\le\Omega <\Omega ^*$ and $0\leq a< a^*:=\|w\|^2_2$, then there exists at least
one minimizer for $e_F(a)$.
\item If $0\le\Omega <\Omega ^*$ and $a \ge a^*:=\|w\|^2_2$, then there
is no minimizer for $e_F(a)$.
\item If $\Omega >\Omega ^*$, then for any $a\ge 0$, there is no minimizer for $e_F(a)$.
\end{enumerate}
Moreover, $e_F(a) > \inf _{\R^2}V_\Omega(x)$ for $0\le\Omega <\Omega ^*$ and $a<a^*$,
$\lim_{a\to a^*}
e_F(a) = e_F(a^*) = \inf _{\R^2}V_\Omega(x)$ for $0\le\Omega <\Omega ^*$, and
{$e_F(a) = -\infty$} for $\Omega \ge 0$ and $a>a^*$, where
\begin{equation}
V_\Omega (x):=
V(x)-\frac{\Omega ^2}{4}|x|^2.
\label{intro:v}
\end{equation}
\end{thm}

{For the reader's convenience, a sketch of the proof for Theorem \ref{thm1} is given in Section 2. From the above conclusion, we have seen that the existence and nonexistence of minimizers for $e_F(a)$
are {the same as the case without rotation, provided that the rotating velocity $\Omega$ of the trap is smaller than the critical value $\Omega ^*$.} In addition to} Theorem \ref{thm1}, we note from Remark
\ref{2:rem} below that for the particular case $\Omega =\Omega ^*<\infty$, the existence and nonexistence of minimizers for $e_F(a)$ seem complicated, which might depend  on the exact trapping profile of $V(x)$.


By the variational argument, if $u$ is a minimizer of $e_F(a)$, then there exists a Lagrange multiplier
$\mu = \mu (a, \Omega )\in \R$ such that $u$ is a complex-valued solution of
the following elliptic equation
\begin{equation}
	-\Delta u+V(x)u+i\, \Omega \, (x^{\perp}\cdot \nabla u)=\mu u+a|u|^2u\quad \mbox{in}
\, \  \R^2.  \label{eqn}
\end{equation}
Furthermore, the argument of \cite[Theorem 1]{GS} implies that if   $u_a$ is a   minimizer of
$e_F(a)$, then we  have
\begin{equation}
\|u_a\|_\infty\to \infty\ \ \mbox{and }\ \   \inte V_\Omega (x)|u_a|^2dx\to V_\Omega (x_0):
	=\inf _{x\in \R^2}
	V_\Omega (x)\quad \mbox{as}\ \ a\nearrow a^*,
\label{intro:lower}
\end{equation}
where $V_\Omega (x)$ is as in (\ref{intro:v}). This motivates us to investigate the limit behavior of  minimizers for $e_F(a)$ as $a\nearrow a^*$, for which we define

\begin{defn}
The function $h(x)\ge 0$ in $\R^2$ is homogeneous of degree $p\in\R^+$ (with respect to the origin), if there exists some $p>0$ such that
\begin{equation}\label{1:V}
		h(tx)=t^ph(x)\ \, \mbox{in}\ \, \R^2 \ \mbox{for any}\ t>0.
\end{equation}
\end{defn}

Following \cite[Remark 3.2]{Grossi}, the above definition implies that the homogeneous function $h(x)\in C(\R^2)$ of degree $p>0$ satisfies
\begin{equation}\label{1:Vh}
	0\le h(x)\le C|x|^p\,\, \mbox{in $\R^2$, where $C=\max_{x\in \partial B_1(0)}h(x)$.}
\end{equation}
{In addition to the assumptions that $V\in L^\infty_{\rm loc}(\R^2)$ satisfies
(\ref{A:V}) and $\Omega ^*$ is strictly positive}, we further assume that $V_{\Omega}(x)=V(x)-\frac{\Omega^2}{4}|x|^2$ satisfies additionally the following conditions:
\begin{enumerate}
\item [\rm($V$).]
 $V_{\Omega}(x)\geq0$, $V_{\Omega}(x)\le C|x|^q$ as $|x|\to\infty$ for some $q\ge 2$, $\{x\in\R^2:\,V_{\Omega}(x)=0\}=\{0\}$ and  $V_{\Omega}(x)=[1+o(1)] h(x)$ as $x\to 0$, where $h(x)\in C^\alpha (\R^2)$  is homogeneous of degree $p>0$ with $0<\alpha<\min\{p, 1\}$ and $H(y)=\inte h(x+y)w^2(x)dx$ admits a global minimum point $y_0\in \R^2$.
\end{enumerate}
Following the above assumption $(V)$, we define $\lambda \in (0,\infty]$ satisfying
\begin{equation}\label{def:li}
\lambda   =\arraycolsep=1.5pt\left\{\begin{array}{lll}
	\Big[\displaystyle\frac{p}{2}  \inte h(x+y_0)w^2(x)dx\Big]
	^{\frac{1}{2+p}}, \quad   &\mbox{if} & \quad 0<p<2;\\[4mm]
	\displaystyle\Big[   \inte \Big(h(x+y_0)+\frac{\Omega ^2 }{4}|x|^2\Big)w^2(x)dx  \Big]
	^{\frac{1}{4}}, \quad   &\mbox{if}& \quad p=2;\\[4mm]
\displaystyle\Big[\frac{\Omega ^2}{4} \inte |x|^2 w^2(x) dx\Big]^{\frac{1}{4}}, \quad   &\mbox{if}& \quad p>2;
\end{array}\right.
\end{equation}
where the point $y_0\in\R^2$ is as in the assumption $(V)$. Under the assumptions (\ref{A:V}) and $(V)$ for some $p>0$, {we have the following result concerning the $L^\infty-$uniform convergence as $a\nearrow a^*$}.

\begin{thm}\label{thm1.2}
Suppose $V(x)\in L^\infty_{\rm loc}(\R^2)$ satisfies (\ref{A:V}) and $(V)$ for some $p>0$, and assume $0<\Omega<\Omega^*$, where $\Omega ^*>0$ is defined as in (\ref{Omega}). If $u_a$ is a minimizer of $e_F(a )$, then we have
\begin{equation}
 w_a(x):=\frac{(a^*-a)^{\frac{1}{2+\gamma}}}{\lambda} u_{a}\Big( \frac{(a^*-a)^{\frac{1}{2+\gamma}}}{\lambda}x +x_a \Big)e^{-i\, \big(\frac{\Omega}{2\lambda}(a^*-a)^{\frac{1}{2+\gamma}}x\cdot x_a^\perp-\theta_a\big)} \to \frac { w( x)}{\sqrt{a^*}}\ \ \hbox{as}\,\ a\nearrow a^*
\label{I:con:a}
\end{equation}
uniformly in $L^\infty (\R^2, \mathbb{C})$, where $\gamma =\min \{p,2\}>0$ and $\theta_a \in [0,2\pi)$ is a properly chosen constant. Here $x_a\in\R^2$ is the unique global maximal point of $|u_{a}|$ as $a\nearrow a^*$.
\end{thm}

The proof of Theorem \ref{thm1.2} {relies heavily on} the non-degenerancy of $w$ and the delicate energy estimates of $e_F(a)$ as well.
As mentioned in \cite[Remark 2.2]{Lewin}, the {method} of \cite{Lewin} can be used to derive the convergence of (\ref{I:con:a}) in an $L^2$ sense for the general cases where the trapping potentials attain  their minimum at the origin and behave at least quadratically {both at zero and infinity}, for instance the quartic-quadratic potential $V (x) = |x|^2 + k|x|^4$ for $k \ge 0$. {Compared with those obtained in \cite{Lewin}, we remark that} our Theorem \ref{thm1.2} holds in the $L^\infty$ sense and is applicable to more general cases where trapping potentials attain their minimum at the origin and behave at least quadratically only at infinity, for instance the general
potential $V (x) =|x|^2 + k|x|^q$ for $k \ge 0$ and $q>0$.

{Finally, we consider a more refined limit behavior of minimizers for $e_F(a)$ as $a \nearrow a^*$, in the case where $V(x)=|x|^2$ and $0<\Omega<\Omega^*$ is fixed. {In this case, recall from (\ref{omega-1}) that $\Omega^*=2$.} Based on the $L^\infty-$uniform convergence of Theorem \ref{thm1.2}, we shall show the following {\em nonexistence of vortices} for $e_F(a)$ confined in the harmonic trap. }


\begin{thm}\label{thm1.2*}
Suppose $V(x)=|x|^2$ and $0<\Omega<\Omega^*=2$. Then, up to a constant phase, all minimizers of $e_F(a)$ are real-valued, unique and free of vortices as $a$ is sufficiently close to $a^*$ from below.	
\end{thm}

As the main result of the current paper, Theorem \ref{thm1.2*} can be generalized to the case {$V(x)=|x-A|^2$
for any point $A=(a_1,a_2)\in \R^2$. This follows directly by the transformation $v(x)=u(x+\frac{4A}{4-\Omega^2})e^{-i\frac{2\Omega}{4-\Omega^2}x\cdot A^{\perp}}$, where $A^{\perp}=(-a_2,a_1)$.} Theorem \ref{thm1.2*} and \cite[Theorem 1.4]{GLW} yield immediately that up to a constant phase $\theta_a\in [0, 2\pi)$, the unique minimizer $u_a$ of $e_F(a)$ satisfies
\begin{equation} \label{1:CC}
\begin{aligned}
u_a(x)=e^{-i\theta_a}u_a^0(x)\approx&\, \displaystyle\frac { \lambda}{\|w\|_2} e^{-i\theta_a} \Big\{  \displaystyle\frac{1}{(a^*-a)^\frac{1}{4}} w\Big(\frac{\lam  |x|}{(a^*-a)^\frac{1}{4}}\Big)+\displaystyle(a^*-a)^\frac{3}{4}\psi\Big(\frac{\lam |x|}{(a^*-a)^\frac{1}{4}}\Big)\\[4mm]
	&+\displaystyle(a^*-a)^\frac{7}{4} \phi_0\Big(\frac{\lam  |x|}{(a^*-a)^\frac{1}{4}}\Big)  +o\big((a^*-a)^\frac{7}{4}\big)\Big\}  \ \ \mbox{as}\, \ a\nearrow a^*,
\end{aligned}
\end{equation}
where $u_a^0(x)>0$ is the unique real minimizer of $e_F(a)$ at $\Omega =0$ as $a\nearrow a^*$, and $\lam>0$ is defined by (\ref{def:li}) with $p=2$. Here the real-valued radially symmetric functions $\psi,\, \phi_0\in C^2(\R^2)\cap L^\infty(\R^2)$ are unique. Therefore, the expansion (\ref{1:CC}) gives a positive answer for the question raised soon after \cite[Theorem 2.1]{Lewin}.

\subsection{Proof strategy of Theorem \ref{thm1.2*}}



The purpose of this subsection is to explain the general strategy of proving Theorem \ref{thm1.2*}, which can be summarized as the method of inductive symmetry. Roughly speaking, our basic idea is to {derive} the radial symmetry of minimizers as $a \nearrow a^*$ by analyzing the refined limit behavior of minimizers and employing the non-degenerancy of $w$.

To prove Theorem \ref{thm1.2*},  the first step is the $L^\infty-$uniform convergence already established in Theorem \ref{thm1.2}. Under the assumptions that $V=|x|^2$ and $0<\Omega<\Omega^*=2$, let $u_a$ be a minimizer of $e_F(a )$.
Stimulated by \eqref{I:con:a}, we consider
\begin{equation}
\arraycolsep=1.5pt
 \begin{array}{lll}
 w_a(x)&: =& \displaystyle\eps_a u_a\Big(\eps_ax  +x_a\Big)e^{-i\Omega \big(
 \frac{\eps_a}{2 } x\cdot x_a^\perp-\theta_a\big)} \\[3mm]
 &:=&\displaystyle {R}_a(x)+iI_a(x),\ \ \eps_a:=\Big(\int_{\R^2}|\nabla u_a|^{2}dx\Big)^{-\frac{1}{2}}>0,
 \end{array}
\label{AA:4:1}
\end{equation}
where $\theta_a\in [0,2\pi)$ is chosen properly as in (\ref{I:con:a}), and $x_a\in\R^2$ is the global maximal point of $|u_a|$.
Here ${R}_a$ and $  I_a $ are real-valued functions in $H^1(\R^2)$. The proof of Theorem \ref{thm1.2} then yields the following $L^\infty-$uniform convergence \begin{equation}
R_a\to \frac{w}{\sqrt{a^*}} \ \ \mbox{and}\ \ I_a\to 0  \ \mbox{uniformly in}\ \,  \R^2  \ \mbox{as}\ \ a\nearrow a^*,
\label{AA:4:2}
\end{equation}
where the non-degenerancy of $w$ is used.

{In the second step, we shall establish Proposition \ref{prop4.3} on the first estimates of the difference $R_a- \frac{w}{\sqrt{a^*}}$ and $I_a$ as $a \nearrow a^*$, by following (\ref{AA:4:2}) and analyzing the following system:
\begin{equation}\label{rev-4:3}
\left\{
\begin{aligned}
\mathcal{L}_a R_a&=\eps^2_a\Omega (x^\bot\cdot\nabla I_a)
\ \ &\mbox{in}\,\  \R^2,\\
\mathcal{L}_a I_a&=-\eps^2_a\Omega (x^\bot\cdot\nabla R_a)
\ \ &\mbox{in}\,\  \R^2,\\
\end{aligned}
\right.
\end{equation}
where
\begin{equation}\label{rev-4:3A}
\mathcal{L}_a:=-\Delta+\Big(\frac{\eps_a^4\Omega^2}{4}|x|^2+\eps_a^2V_\Omega (\eps_a x+x_a)-\e^2\mu _a-a|w_a|^2\Big).
\end{equation}
Here, the system \eqref{rev-4:3} is obtained from (\ref{eqn}) and (\ref{AA:4:1}).}

The third step is the first estimates of $R_a- v_a$ and $I_a$ as $a \nearrow a^*$, where $v_a>0$ is the rescaled positive minimizer of $e_F(a)$ at zero rotation $\Omega=0$. Note from \cite{GLW} that $v_a>0$ is radially symmetric and unique as $a\nearrow a^*$. We should emphasize that we are unable to derive the further refined estimates $R_a $ and $I_a$ as $a \nearrow a^*$ by further analyzing the difference $R_a- \frac{w}{\sqrt{a^*}}$, which is mainly due to the non-symmetry of $\mathcal{L}_a$. The main results of this step are given in Proposition \ref{prop4.4} by applying Proposition \ref{prop4.3} and some refined estimates of \cite{GLW}.

{Associated to each $ v_a$, there is a linearized operator $\mathcal{N}_a:\, D(\mathcal{N}_a)\subset L^2(\R^2) \longmapsto L^2(\R^2)$ defined as
\begin{equation}\label{A:G10}
\mathcal{N}_a:=-\Delta+\big(\e^4|x|^2-\e^2\beta_a-3av_a^2\big),
\end{equation}
where $\e >0$ is defined in (\ref{AA:4:1}), and the domain $D(\mathcal{N}_a)$ is defined by
$$D(\mathcal{N}_a)=H^2(\R^2)\cap \Big\{u\in L^2(\R^2):\inte |x|^4u^2dx< \infty \Big\}.$$
The fourth step is to study the operator $\mathcal{N}_a$ and establish the inductive estimates of Lemma \ref{lem4.6}.}

Following the estimates in previous steps, the last step is the complete proof of Theorem \ref{thm1.2*} by an induction process.  Heuristically speaking, we shall show that once we have a good control on the non-radial part of {$R_a-v_a$}, then a better control on the non-radial part of $R_a-v_a$ could be obtained. As a consequence,  we shall finally prove in Subsection 5.2 that the non-radial part of $R_a-v_a$ is arbitrarily small and
\[ \mbox{both}\ \, I_a(x)\equiv 0 \ \, \mbox{and}\ \, x_a\equiv 0\ \, \mbox{as}\ \,   a\nearrow a^*.\]
By the uniqueness of $v_a>0$ as $a\nearrow a^*$ (cf. \cite[Theorem 1.3]{GLW}), we are able to conclude from above and (\ref{AA:4:1}) that up to a constant phase, all minimizers of $e_F(a)$ must be real-valued, unique and free of vortices as $a\nearrow a^*$, and thereby Theorem \ref{thm1.2*} is proved. For the details of the proof of Theorem 1.3, we refer the reader to Section 5.

It is worthy of mentioning that the above ``inductive symmetry" argument can be applied for any radially symmetric potential like $|x|^2+t|x|^p$ with positive $t$ and $p\in(1,+\infty).$ Indeed, it can be generalized to the general radially symmetric function which admits polynomial growth at most, provided that the following two points are ensured: at first, Theorem 5.3 still holds for the given potential; secondly, we have to change the assumptions of Lemma 5.4 and Lemma \ref{rev-7} a little bit by requiring that $|T_{m,a}|$ and $|f_{i}|$ decay  slightly faster than $|\nabla T_{m,a}|$ and $|\nabla f_i|$, respectively, where $i=1, 2$.

{Before we end the introduction}, it deserves to remark that the nonexistence of vortices for rotating BECs with repulsive interactions were proved earlier in \cite{AJ,IM-1,IM-2, AA, CD, CP, CR, CRY, Ro} by different arguments, including jacobian estimates, vortex ball constructions, and so on. We also refer the celebrated monograph \cite{SSbook} to various kinds of techniques involved in the investigation of vortex structures for defocusing nonlinear problems. As far as we know, it seems that the above mentioned methods of studying vortex structures concentrate more on the analysis of the energy, which is not enough for our focusing nonlinear case. Indeed, besides the refined analysis of the energy,  the proof of our Theorem \ref{thm1.2*} relies more on the refined limit analysis of the governing system (\ref{rev-4:3}), for which one needs to make full use of the non-degenerancy of $w$.

This paper is organized as follows: in Section 2 we shall sketch briefly the proof of Theorem \ref{thm1}
under more general potentials satisfying (\ref{A:V}). The proof of Theorem \ref{thm1.2} is given in Section 3 by employing energy methods and blow-up analysis. Section 4 is devoted to the proof of Proposition \ref{prop4.3} on the refined limit behavior of minimizers for $e_F(a)$ as $a\nearrow a^*$. In Section 5, we first establish Proposition \ref{prop4.4}, and  the  linearized operator $\mathcal{N}_a$ is then analyzed in Subsection 5.1. We shall complete the proof of Theorem \ref{thm1.2*}  in Subsection 5.2 by the so-called method of inductive symmetry described above. In Appendix A, we shall prove two lemmas which are used in the proof of Theorem \ref{thm1.2*}.

\section{Existence of Minimizers}
In this section,  we shall address the proof of Theorem \ref{thm1} on the existence
and nonexistence of  minimizers for  $e_F(a)$. Since the proof of Theorem \ref{thm1} is overall similar to  \cite[Theorem 1]{GS}, we only sketch out the main differences.

Letting $w>0$ be the unique positive radial solution of (\ref{Kwong}), note that $w$  strictly decreases in $|x|$ (cf. \cite[Proposition 4.1]{GNN}), and admits the following exponential decay
\begin{equation}
	w(x) \, , \ |\nabla w(x)| = O(|x|^{-\frac{1}{2}}e^{-|x|}) \quad \text{as}\,\   |x|\to \infty ,   \label{4:exp}
\end{equation}
which are often used throughout the paper.
We next recall some results and notations from \cite{GS} on the minimizers of $e_F(a)$
in the special case where  $\Omega =0$. Equivalently,
we consider the following Gross-Pitaevskii energy functional
\begin{equation}
	E_a(u):=\int _{\R ^2} \big(|\nabla u|^2+V(x)|u|^2\big)dx-\frac{a}{2}\int _{\R ^2}|u|^4dx,
	\quad \ u\in \h    \label{non:f}
\end{equation}
under the unit mass constraint $\int _{\mathbb{R} ^2} |u|^2dx=1$, i.e.,
\begin{equation}
	e(a):=\inf _{\{u\in \,\h \ {\rm and }\ \inte |u|^2dx=1 \} } E_a(u),\label{non:ff}
\end{equation}
where $\h$ is defined as in (\ref{1.H}). One can note that $E_a(|u|)\le E_a(u)$ holds for any $u\in \h $. This property and
\cite[Theorem 1]{GS} yield immediately the following existence and nonexistence:

\begin{thm}\label{thmA} Let $w$ be the unique positive
radial solution of (\ref{Kwong}). Suppose $V\in L^\infty_{\rm loc}(\R^2)$ satisfies
$\lim_{|x|\to\infty} V(x) = \infty$. Then we have
\begin{enumerate}
\item If $0\leq a< a^*:=\|w\|^2_2$, then there exists at least one minimizer
			for  (\ref{non:ff}).
\item If  $a \ge a^*:=\|w\|^2_2$, then there is no  minimizer for (\ref{non:ff}).
\end{enumerate}
\end{thm}


To prove Theorem \ref{thm1} for the rotational case where $\Omega >0$, one needs to use the following diamagnetic
inequality (cf. \cite{Lieb}):
for $\mathcal{A }=\frac{\Omega}{2}x^{\perp}$,
\begin{equation}  |\nabla u|^2 -\Omega \, x^{\perp}\cdot (iu,\, \nabla u) =
	|(\nabla -i\mathcal{A } )u|^2 -\frac{\Omega ^2}{4}  |x|^2|u|^2\ge \big| \nabla |u|\big|^2
	-\frac{\Omega ^2}{4}  |x|^2|u|^2, \quad u\in H^1(\R^2, \mathbb{C}) .
	\label{2:2:1A}
\end{equation}
Define the critical rotational velocity $\Omega ^*$ as in (\ref{Omega}), such that
if $V(x)$ satisfies (\ref{A:V}), then $\Omega ^*$ exists and satisfies
$0<\Omega ^*\le\infty$.
Similar to those in \cite{Zhang}
and references therein,  one can establish that
for $V \in L_{\rm loc}^\infty(\R^2)$ satisfying $\lim_{|x|\to \infty} V(x) = \infty$, if  $2\le q<\infty$, then the embedding $\h  \hookrightarrow L^{q}(\R^2, \mathbb{C})$ is compact.
Applying the Gagliardo-Nirenberg inequality (\ref{GNineq}) and the diamagnetic
inequality (\ref{2:2:1A}), similar to \cite[Theorem 1]{GS} it is standard to use the above compactness property to derive the following existence.

\begin{thm}\label{2a:thm1} Let $w=w(|x|)$ be the unique positive radial solution
of (\ref{Kwong}). Suppose $V(x)\in L^\infty_{loc}(\R^2)$ satisfies (\ref{A:V}) such that $   \Omega ^*>0 $ exists, where $\Omega ^*$ is   defined as
in (\ref{Omega}).
Then for any $0<\Omega <\Omega ^* $ and  $0\leq a< a^*:=\|w\|^2_2$,   there exists at
least one minimizer for $e_F(a)$.
\end{thm}

We now make the following observation: by applying (\ref{GNineq}), the second inequality of (\ref{2:2:1A}) yields that
\begin{equation}
F_a(u) \ge \Big(1-\displaystyle\frac{a}{a^*}\Big)\inte \big| \nabla |u|
	\big|^2dx+\inte V_\Omega (x)|u|^2dx,\quad u\in \h ,
	\label{2:ineq}
\end{equation}
where $V_\Omega (x)=  V(x)-\frac{\Omega ^2}{4}|x|^2$. Choose a non-negative function $\varphi \in C_0^\infty(\R^2)$ such that $\varphi(x) = 1$ for $|x|\leq 1$. For any large $\tau>0$, we then consider the test function
\begin{equation}\label{def:trial}
	w_{\tau}(x) =  A_{\tau} \frac{\tau}{\|w\|_2} \varphi(x-x_0) w\big(\tau (x-x_0)\big)e^{i\Omega S(x)} \,,
\end{equation}
where $x_0\in \R^2$ is a suitable point, $S(x)=\frac{1}{2}x\cdot x_0^\perp$, and $A_{\tau}>0$ is chosen so that $\inte |w_\tau (x)|^2dx=1$.
Employing above estimates and the test function (\ref{def:trial}), the following nonexistence result can be established in the similar way of proving \cite[Theorem 1]{GS}.

\begin{thm}\label{2b:thm1} Let $w=w(|x|)$ be the unique positive radial solution
		of (\ref{Kwong}). Suppose $V(x)\in L^\infty_{loc}(\R^2)$ satisfies (\ref{A:V}) such that $   \Omega ^*>0 $ exists, where $\Omega ^*$ is  defined as in (\ref{Omega}).
Then there is no minimizer for $e_F(a)$, provided that either
\begin{enumerate}
			\item   $0<\Omega <\Omega ^* $ and $a \ge a^*:=\|w\|^2_2$, or
			\item   $\Omega >\Omega ^* $ and   $a \ge 0$.
\end{enumerate}
Moreover, $e_F(a) > \inf _{\R^2}V_\Omega(x)$ for $0\le\Omega <\Omega ^*$ and $a<a^*$, $\lim_{a \nearrow a^*}e_F(a) = e_F(a^*) = \inf _{\R^2}V_\Omega(x)$ for $0\le\Omega <\Omega ^*$, and $e(a) = -\infty$ for $\Omega \ge 0$ and $a>a^*$.
\end{thm}
	
Combining with Theorems \ref{thmA} and \ref{2a:thm1}, whose proofs are omitted for simplicity, Theorem \ref{2b:thm1} completes the proof of Theorem \ref{thm1}.
	
\begin{rem}\label{2:rem} As the complement of Theorem 1.1, we guess that the existence and nonexistence of minimizers in the case where $\Omega =\Omega ^*<\infty$ depend on $a$ and $V(x)$. Specifically, consider the trapping profile $V(x)=|x|^2$ and $\Omega =\Omega ^*(=2)$. When $a=0$, by the following inequality (see \cite[Remark 2.5]{EL})
$$\int_{\mathbb{R}^2} |\nabla u-ix^{\perp}u|^2dx\geq 2\int_{\mathbb{R}^2} |u|^2dx,$$
we have $e_F(0)\geq2$,  and one can further check that $e_F(0)\leq 2$ by using $u=\frac{1}{\sqrt{\pi}}e^{-\frac{|x|^2}{2}}$ as a test function. Thus, $e_F(0)$ admits a minimizer $u=\frac{1}{\sqrt{\pi}}e^{-\frac{|x|^2}{2}}$ in this case. While if $a=a^*$, we have $F_{a^*}(u)\geq0$,  and we can further verify that $e_F(a^*)=0$ by taking $u_\tau=\frac{1}{\sqrt{a^*}}\tau w(\tau x)$  as a test function and letting $\tau\to\infty$. Then, it is not difficult to check that there is no minimizer in the  case where $a=a^*$.

\end{rem}

\section{Mass Concentration as $a\nearrow a^*$}

By employing energy methods and blow-up analysis, the main purpose of this section is to establish Theorem \ref{thm1.2} on the $L^\infty-$uniform convergence of minimizers for $e_F(a)$ as $a\nearrow a^*$.
Recall from Theorem~\ref{thm1} that $e_F(a)$ admits minimizers and $e_F(a)\to \inf _{\R^2}V_\Omega(x)$ as $a\nearrow a^*$, where $V_\Omega(x)$ satisfies the assumption $(V)$ for any $p>0$, and $e_F(a)$ is the GP energy defined in $(\ref{def:ea})$.

To prove Theorem \ref{thm1.2}, we first address the following energy estimates of $e_F(a)$ as $a\nearrow a^*$.

\begin{lem}\label{lem3.2}
Suppose that $V(x)\in L^\infty_{\rm loc}(\R^2)$ satisfies (\ref{A:V}) and $(V)$ for some $p>0$, and assume $0<\Omega<\Omega^*$, where $\Omega ^*>0$ is defined in (\ref{Omega}).  Then there exist two positive constants $m$ and $M$, independent of $0<a<a^*$, such that
\begin{equation}
		0<m(a^*-a)^\frac{p}{p +2}\le e_F(a)\le M(a^*-a)^\frac{\gamma}{\gamma +2}\ \ \mbox{as}\ \
		a\nearrow a^*,
		\label{4:con:1}
\end{equation}
where $\gamma =\min \{p, 2\}>0$.
\end{lem}
	
\noindent{\bf Proof.} By (\ref{2:ineq}), we first note that a similar argument of \cite[Lemma 3]{GS} gives the lower bound of (\ref{4:con:1}).

To prove the upper bound of (\ref{4:con:1}),  similar  to   Theorem~\ref{2b:thm1}, we use a test function $w_\tau$ of the form (\ref{def:trial}) with $S(x)\equiv 0$, where $A_{\tau} >0$ satisfies
\begin{equation}\label{up1}
 	\frac 1{A_{\tau}^{2}} = \frac{1}{\|w\|_2^2} \int_{\R^2} w(x)^2 \varphi(x/ \tau  )^2 dx
 	= 1 + O( \tau^{-\infty}) \quad\text{as}\,\ \tau \to \infty,
\end{equation}
due to the exponential decay of $w$ in Eq.~(\ref{4:exp}), where and below we always use $f(t) = O(t^{-\infty})$ to denote a
function $f$ satisfying $\lim_{t\to \infty} |f(t)| t^s = 0$ for all $s>0$. Under the assumption $(V)$ for some $p>0$, choose $\mathcal{R}>0$ small enough so that
$V_\Omega(x) \leq C_0 |x|^p $ for $ |x|\leq \mathcal{R}$,
where $V_\Omega(x)$ is defined as in \eqref{intro:v}.
In view of this fact, we have
$$
	\int_{\R^2} V_\Omega (x) |w_\tau(x)|^2 dx \leq \frac{C_0  A_{\tau}^2}{\tau^{p}}
	\int_{\R^2} |x|^p w^2(x) dx \quad \mbox{as}\,\ \tau\to\infty .
$$
We then obtain from (\ref{def:trial}) that
\[\arraycolsep=1.5pt \begin {array}{lll}
	\displaystyle
	C_{\tau }(w_\tau):&=&\displaystyle\inte \Big(|\nabla w_\tau|^2 -\Omega\,  x^{\perp}
	\cdot (iw_\tau ,\, \nabla w_\tau)
	+\displaystyle\frac{\Omega ^2}{4}|x|^2 |w_\tau|^2\Big)dx -
	\displaystyle\frac{a}{2}\inte |w_\tau|^4dx\\[3mm]
	&=&\displaystyle\frac{\tau^2}{2\|w\|_2^2} \Big( 1-\frac{a}{a^*}\Big)
	\inte w^4(x)dx +\frac{\Omega ^2}{4\|w\|_2^2\tau ^2}\inte |x|^2w^2(x)dx+O(  \tau ^{-\infty})
\end{array}\]
as  $\tau\to\infty $. Thus, we conclude from above that
\begin{equation}\arraycolsep=1.5pt \begin {array}{lll}
	\displaystyle
	e_F(a) &\leq & \displaystyle\frac{\tau^2}{2 (a^*)^2} \left( a^*-a \right) \inte w^4(x)dx +
	\displaystyle\frac{C_0   }{\tau^{p}}
	\int_{\R^2} |x|^p w^2(x) dx\\[3mm]
	&&  +\displaystyle\frac{\Omega ^2}{4\|w\|_2^2\tau ^2}
	\inte |x|^2w^2(x)dx+O( \tau ^{-\infty})\quad \mbox{as}\,\ \tau\to\infty.
	\end{array}
	\label{3:energy-P}
\end{equation}
Optimizing the right hand side of (\ref{3:energy-P}) by taking
$$\tau =(a^*-a)^{-\frac{1}{\gamma +2}},\quad \text{where}\
	\gamma :=\min \big\{p,2\big\}>0,$$
we arrive at the desired upper bound of (\ref{4:con:1}), and the proof is therefore complete.
\qed
\medskip

For a given $0<\Omega <\Omega ^*$, if $u_a$ is a minimizer of $e_F(a)$ for $0<a<a^*$, then $u_a$ is a solution of the following elliptic equation
\begin{equation}
	-\Delta u_a+V(x)u_a+i\, \Omega \, (x^{\perp}\cdot \nabla u_a)=\mu_a u_a+a|u_a|^2u_a\quad \mbox{in}
	\ \, \R^2,  \label{3:eqn}
\end{equation}
where $\mu _a\in \R$ is a Lagrange multiplier. One can check that $\mu _a$ satisfies
\begin{equation}
	\mu _a=e_F(a)-\frac{a}{2}\inte |u_a|^4dx,\,\ 0<a<a^*.
	\label{3:mu}
\end{equation}
The following lemma gives the refined estimates of $u_a$ and $\mu_a$ as $a\nearrow a^*$.

\begin{lem}\label{lem3.3}
Set
\begin{equation}\label{eea}
		\eps_{a}:=\Big(\inte|\nabla u_{a}|^{2}dx\Big)^{-\frac{1}{2}}>0,
\end{equation}
where $u_{a}$ is a minimizer of $e_F(a)$. Suppose that the assumption of Theorem \ref{thm1.2} holds. Then we have
\begin{enumerate}
\item  $\eps_{a}>0$ satisfies $\eps_{a}\rightarrow 0$ and $\mu_a\eps ^2_{a}\to -1$ as $a\nearrow a^{\ast}$.
\item  Define the normalized function as follows:
\begin{equation}\label{3.9}
			w_a(x):=\eps_a u_a\big(\eps_a x+x_a\big)e^{-i (\theta_a+\frac{\eps _a \Omega}{2} x\cdot x_{a}^{\bot})},
			\end{equation}
where $x_a$ is a global maximal point of $|u_a|$ and $\theta_a\in [0,2\pi)$ is a proper constant. Then there exists a constant $\eta>0$, independent of $0<a<a^*$, such that
\begin{equation}\label{NEW6}
			\int_{B_2(0)} |w_a(x)|^2dx\geq \eta>0\quad \mbox{as}\ \ a\nearrow a^*.
\end{equation}
			
\item  Any global maximal point $x_a$ of $|u_a|$ satisfies  $\lim_{a\nearrow a^*}V_\Omega (x_a)=0$, and $w_{a}(x)$ satisfies
			\begin{equation} \label{3NEW7}
			\lim_{a\nearrow a^*} w_{a}(x)=\frac{w(x) }{\sqrt{a^*}}\quad \hbox{strongly in $H^1(\R^2,\C)$. }
			\end{equation}
\end{enumerate}
\end{lem}

\noindent{\bf Proof.} 1. On the contrary, suppose that $\eps _{a}\nrightarrow0$ as $a\nearrow a^{\ast}$. By the definition of $e_F(a)$, then we derive from (\ref{GNineq}) and the diamagnetic inequality (\ref{2:2:1A}) that
\begin{equation}\arraycolsep=1.5pt
\begin {array}{lll}
	e_F(a)&= &
	\displaystyle \inte |(\nabla -i\mathcal{A} )u_a|^2dx -
	\displaystyle\frac{a}{2} \inte |u_a|^4dx+ \inte V_{\Omega}(x)|u_a|^2dx
	\ge \displaystyle\inte V_{\Omega}(x)|u_a|^2dx.
	\end{array}\label{New1}
\end{equation}
Following the upper estimate of $e_F(a)$ in \eqref{4:con:1}, we deduce from the above inequality that the minimizer $u_a$ is bounded uniformly in $\h$. Furthermore, since the embedding $\h  \hookrightarrow L^{q}(\R^2, \mathbb{C})$ is compact for $2\le q<\infty$, we obtain that there exists a subsequence of $\{u_a\}$, still denoted by $\{u_a\}$,  such that
\[
	u_a\rightharpoonup u_0\ \ \mbox{weakly in }\   \h,\quad  u_a\to u_0\ \ \mbox{strongly in }\   L^q(\R^2, \mathbb{C})
\]
for some $u_0\in \h$, where $q\in[2,\infty)$. Since $u_a\to u_0$ strongly in $L^2(\R^2, \mathbb{C})$
as $a\nearrow a^*$, we now have $\|u_0\|^2_2=1$. However, it follows from \eqref{4:con:1} and \eqref{New1} that
\begin{equation}\arraycolsep=1.5pt \begin {array}{lll}
	0=\lim_{a\nearrow a^*} e_F(a)\geq  \liminf_{a\nearrow a^*}\displaystyle\inte V_{\Omega}(x)|u_a|^2dx\geq \displaystyle\inte V_{\Omega}(x)|u_0|^2dx\geq 0,
\end{array}\label{New2}
\end{equation}
which implies that $u_0\equiv 0$ in $\R^2$ by the assumption $(V)$, a contradiction. Therefore, we have $\eps _{a}\to 0$ as $a\nearrow a^{\ast}$.

We next prove that $\mu_a\eps_a^2\to -1$ as $a\nearrow a^*$. Actually, since $0<\Omega <\Omega ^*$ is fixed,  the definition (\ref{Omega}) of $\Omega ^*$ implies that
\begin{equation}
|x|^2 \le C(\Omega) V_{\Omega}(x) \quad \mbox{for sufficiently large}\, \ |x|>0.
\label{2:add:1A}
\end{equation}
By (\ref{2:add:1A}) we  then obtain from  \eqref{4:con:1} and \eqref{New1} that for any given large constant $M>0$,
\begin{equation}\label{New3}\
\begin{split}
	\inte \frac{\Omega^2|x|^2}{4}|u_a|^2dx&=\int_{|x|\leq M} \frac{\Omega^2|x|^2}{4}|u_a|^2dx+\int_{|x|> M} \frac{\Omega^2|x|^2}{4}|u_a|^2dx\\
	&\leq C+C(\Omega)\int_{|x|> M} V_{\Omega}(x)|u_a|^2dx\leq C.
\end{split}
\end{equation}
It then follows from above that for any given $\sigma>0$,
\begin{equation}\label{New4}
	\begin{split}
	\Omega \Big|\inte x^\perp\cdot(iu_a,\nabla u_a)dx\Big| &\leq  \frac{\sigma}{2}\inte |\nabla u_a|^2dx+\frac{\Omega^2}{2\sigma}\inte |x|^2|u_a|^2dx\\
	&\leq  \frac{\sigma}{2}\inte |\nabla u_a|^2dx+C(\sigma).
\end{split}\end{equation}
By the Gagliardo-Nirenberg inequality \eqref{GNineq} and  Lemma \ref{lem3.2}, we also have
\[
	\limsup_{a\nearrow a^*}\displaystyle\frac{\frac{a}{2}\inte |u_a|^4dx-e_F(a)}{\inte |\nabla u_a|^2dx}\leq 1.
\]
Suppose now that $\mu_a\eps ^2_{a}\nrightarrow -1$ as $a\nearrow a^{\ast}$. By taking a subsequence if necessary, we then obtain from above and (\ref{3:mu}) that there exists some constant $\gamma_0>0$ such that
\begin{equation}\label{New4A}
	-\liminf_{a\nearrow a^*}\mu_a\eps ^2_{a}=\liminf_{a\nearrow a^*}\displaystyle\frac{\frac{a}{2}\inte |u_a|^4dx-e_F(a)}{\inte |\nabla u_a|^2dx}\leq 1-\gamma_0.
\end{equation}
Applying \eqref{New2}--\eqref{New4A}, we thus deduce that
\begin{equation*}
	\begin{split}
	e_F(a)&=\inte |\nabla u_a|^2dx -
	\displaystyle\frac{a}{2} \inte |u_a|^4dx+ \inte V(x)|u_a|^2dx-\Omega \inte x^\perp\cdot(iu_a,\nabla u_a)\\
	&\geq \inte |\nabla u_a|^2dx-
	\displaystyle\frac{a}{2} \inte |u_a|^4dx-\frac{\sigma}{2}\inte |\nabla u_a|^2dx-C(\sigma)\\
	&\geq \big(\gamma_0-\frac{\sigma}{2}\big)\inte |\nabla u_a|^2dx-C(\sigma)= \frac{\gamma_0}{2}\inte |\nabla u_a|^2dx-C(\gamma_0),
	\end{split}
\end{equation*}
where we take $\sigma=\gamma_0>0$ in the last identity. Since we have $\eps _{a}\to 0$ as $a\nearrow a^{\ast}$, it follows from (\ref{eea}) that $\inte |\nabla u_a|^2dx\to \infty$  as $a\nearrow a^{\ast}$. The above inequality then yields that $e_F(a)\to\infty$ as $a\nearrow a^{\ast}$, which however contradicts to Lemma \ref{lem3.2}. This completes the proof of (1).

2. Denote $\bar w_{a}(x)=\eps_a u_a\big(\eps_a x+x_a\big)e^{-i (\frac{\eps _a \Omega}{2} x\cdot x_{a}^{\bot})}$ and $w_{a}(x):=\bar w_{a}(x)e^{i\theta_a}$, where the parameter $\theta_a\in [0,2\pi)$ is chosen properly such that
\begin{equation}\label{orth con}
	\Big\|  w_{a}-\frac{w}{\sqrt{a^*}}\Big\|_{L^2(\R^2)}=\min_{\theta\in [0,2\pi)}\Big\|e^{i\theta} \bar w_{a}-\frac{w}{\sqrt{a^*}}\Big\|_{L^2(\R^2)}.
	\end{equation}
Rewrite $w_{a}(x)=R_{a}(x)+iI_{a}(x)$, where $R_a(x)$ denotes the real part of $w_{a}(x)$, and $I_{a}(x)$ denotes the imaginary part of $w_{a}(x)$.  By \eqref{orth con}, we obtain the following orthogonality condition
\begin{equation}\label{new10}
	\inte  w(x) I_{a}(x)dx=0.
\end{equation}
Following (\ref{3:eqn}), we see that $w_a(x)$ satisfies the following Euler-Lagrange equation
\begin{equation}\label{NEW8}
	- \Delta w_a+i\eps^2_a\Omega\, \big(x^\bot\cdot\nabla w_a\big)+\Big[\frac{\eps_a^4\Omega^2}{4}|x|^2+\eps_a^2V_\Omega (\eps_a x+x_a)-\eps_a^2\mu _a-a |w_{a}|^2\Big] w_a=0
	\quad \mbox{in}\,\ \R^2,
\end{equation}
where $w_{a}$ is bounded uniformly in $H^1(\R^2,\C)$.
Define $W_a(x) =|w_a(x)|^2\ge 0$ in $\R^2$. We then derive from \eqref{NEW8} that
\begin{equation}\label{NEW9}\arraycolsep=1.5pt\begin{array}{lll}
	&&
	-\displaystyle\frac{1}{2} \Delta W_a+|\nabla w_a|^2-\eps^2_a\Omega\, x^\bot\cdot(iw_a,\nabla w_a)\\[2mm]
	&&+\displaystyle\Big[\frac{\eps_a^4\Omega^2}{4}|x|^2+\eps_a^2V_\Omega (\eps_a x+x_a)-\eps_a^2\mu _a-a W_a\Big] W_a=0
	\quad \mbox{in}\,\ \R^2.\end{array}
\end{equation}
Since
$$|\nabla w_a|^2-\eps^2_a\Omega\, x^\bot\cdot(iw_a,\nabla w_a)+\frac{\eps_a^4\Omega^2}{4}|x|^2 W_a \geq 0\ \ \mbox{in} \,\ \R^2,$$
which is due to the diamagnetic inequality (\ref{2:2:1A}), we have
\begin{equation}\label{NEW10}
	-\frac{1}{2} \Delta W_a-\eps_a^2\mu _aW_a-aW_a^2\leq 0 \quad \mbox{in}\,\ \R^2.
\end{equation}
Since $0$ is a global maximal point of $W_a(x)$ for all $a<a^*$, we have $-\Delta W_a(0)\geq 0$ for all $a<a^*$. Combining the fact that $\eps_a^2\mu _a\to -1$ as $a\nearrow a^*$, we then get that $W_a(0)\geq \beta>0$ holds uniformly in $a$ for some positive constant $\beta$. Following   De Giorgi-Nash-Moser theory \cite[Theorem 4.1]{HL}, we thus deduce from \eqref{NEW10} that
\begin{equation}\label{NEW11}
	\int_{B_2(0)}W_a^2dx\geq C\max_{x\in B_1(0)}W_a(x)\geq C_1(\beta)\quad \hbox{as}\,\   a\nearrow a^*.
\end{equation}
Moreover, since $w_a$ is bounded uniformly in $H^1(\R^2,\C)$, we obtain from (\ref{NEW11}) that
\begin{equation}\label{NEW12}
	W_a(0)=\max_{x\in\R^2}W_a(x)\leq  C\int_{B_2(0)}W_a^2dx\leq C.
	\end{equation}
Combining \eqref{NEW11} and \eqref{NEW12}, we get that
\begin{equation}\label{NEW13}
	\int_{B_2(0)}|w_a|^2dx\geq \int_{B_2(0)}\frac{W_a^2}{\max_{x\in\R^2}W_a(x)}dx\geq C_2(\beta)>0,
\end{equation}
which completes the proof of (\ref{NEW6}).

3. Since $|w_{a}|$ is bounded uniformly in $H^1(\R^2)$, we may assume that up to a subsequence if necessary, $|w_{a}|$ converges to $w_0$ weakly in $H^1(\R^2)$ as $a\nearrow a^*$ for some $0\le w_0\in H^1(\R^2)$. Note from  \eqref{NEW13} that $w_0\not\equiv  0$ in $\R^2$. By the weak convergence, we may assume that $|w_{a}|\to w_0$ almost everywhere in $\R^2$ as $a\nearrow a^*$. Applying Br\'{e}zis-Lieb  lemma gives that
\begin{equation}\label{3:C:1}
\|w_{a}\|_q^q=\|w_0\|_q^q+\big\||w_{a}|-w_0\big\|_q^q+o(1)\quad \mbox{as} \,\ a\nearrow a^*,\quad\hbox{where $2\leq q< \infty$,}
\end{equation}
	and
\begin{equation}\label{3:C:2}
\big\|\nabla |w_{a}|\big\|_2^2=\big\|\nabla  w_0 \big\|_2^2+\big\|\nabla (|w_{a}|- w_0)\big\|_2^2+o(1)\quad \mbox{as} \,\ a\nearrow a^*.
\end{equation}
This also implies that $\big\||w_{a}|-w_0\big\|_2^2\le 1$ holds uniformly as $a\nearrow a^*$. Recall from (\ref{4:con:1}) that $e_F(a)\to 0$ as $a\nearrow a^*$, and note also that
\[
	1=-\lim_{a\nearrow a^*}\mu_{a}\eps ^2_{a}=\lim_{a\nearrow a^*}\displaystyle\frac{\frac{a}{2}\inte |u_{a}|^4dx-e_F(a)}{\inte |\nabla u_{a}|^2dx}.
\]
By the definition of $w_a(x)$, we therefore get that
\begin{equation}\label{NEW22}
	\lim_{a\nearrow a^*}\inte |w_{a}|^4dx=\frac{2}{a^*}.
\end{equation}
Thus, we have
\begin{equation}\label{3:18} \begin{split}
	&\lim_{a\nearrow a^*}\Big\{\inte \big|\nabla |w_{a}|\big|^2dx-\frac{a^*}{2}\inte |w_{a}|^4dx\Big\}\\
	&\le  \lim_{a\nearrow a^*}\Big\{\inte \big|\nabla  w_{a} \big|^2dx-\frac{a^*}{2}\inte |w_{a}|^4dx\Big\} =0.
\end{split}\end{equation}
By the Gagliardo-Nirenberg inequality \eqref{GNineq}, we then get from (\ref{3:C:1}), (\ref{3:C:2}) and (\ref{3:18}) that
\begin{equation*}
	\begin{split}
	0&\ge\lim_{a\nearrow a^*}\Big\{\inte \big|\nabla |w_{a}|\big|^2dx-\frac{a^*}{2}\inte |w_{a}|^4dx\Big\} \\
	&=\inte |\nabla w_{0}|^2dx-\frac{a^*}{2}\inte |w_{0}|^4dx\\
	&~\quad +\lim_{a\nearrow a^*}\Big\{\inte \big|\nabla(|w_{a}|-w_{0})\big|^2dx-\frac{a^*}{2}\inte \big||w_{a}|-w_{0}\big|^4dx\Big\}\\
	& \geq \inte |\nabla w_{0}|^2dx-\frac{a^*}{2}\inte |w_{0}|^4dx\\
	&~\quad +\frac{a^*}{2}\lim_{a\nearrow a^*}\Big(\big\||w_{a}|-w_0\big\|^{-2}_2-1\Big)\inte \big||w_{a}|-w_{0}\big|^4dx \geq 0.
\end{split}
\end{equation*}
Using \eqref{GNineq} again, the above inequality implies that $\|w_0\|^{2}_2=1$ and
\begin{equation}\label{18}
	|w_{a}|\to w_0\quad \hbox{strongly in}\, \ L^2(\R^2)\quad \mbox{as} \,\ a\nearrow a^*.
\end{equation}
Due to the uniform boundedness of $|w_{a}|$ in $H^1(\R^2)$, we obtain that $|w_{a}|\to w_0$ strongly in $L^4(\R^2)$ as $a\nearrow a^*$.  By the weak lower semicontinuity and \eqref{GNineq}, we further get that $\nabla|w_{a}|\to\nabla w_0$ strongly in $L^2(\R^2)$ as $a\nearrow a^*$. Therefore, $\sqrt{a^*}w_0$ must be an optimizer of the Gagliardo-Nirenberg inequality \eqref{GNineq}, and there exists $y_0\in\R^2$ such that up to a subsequence if necessary,
\[ |w_{a}(x)|\to w_0(x)=\frac{w(x+y_0)}{\sqrt{a^*}}\quad\hbox{strongly in $H^1(\R^2)$}\quad \mbox{as} \,\ a\nearrow a^*.\]
Since the origin is a global maximal point of $|w_{a}|$ for all $a\in(0,a^*)$, it must be a global maximal point of $w(x+y_0)$, which implies that $y_0=0$, and up to a subsequence if necessary,
\begin{equation}\label{NEW15}
	 |w_{a}(x)|\to \frac{w(x)}{\sqrt{a^*}}\quad\hbox{strongly in}\,\ H^1(\R^2)\quad \mbox{as} \,\ a\nearrow a^*.
\end{equation}
Moreover, because the convergence \eqref{NEW15} is independent of the subsequence $\{|w_{a}|\}$, we conclude that \eqref{NEW15} holds for the whole sequence. Then it is not difficult to derive from (\ref{NEW15}) that up to a subsequence if necessary,
\begin{equation}
\label{ffin}
\lim_{a\nearrow a^*}w_a=\frac{w}{\sqrt{a^*}}e^{i\sigma}~\mbox{strongly in}~H^1(\mathbb{R}^2,\mathbb{C})
\end{equation}
for some $\sigma\in\mathbb{R}$. Moreover, we have $\sigma =0$ in view of (\ref{new10}). Since the convergence of \eqref{ffin} is independent of the choice of the subsequence, we deduce that \eqref{ffin} holds essentially true for the whole sequence, and hence (\ref{3NEW7}) holds true.

Following \eqref{4:con:1}, we also obtain that $\lim_{a\nearrow a^*}e_F(a)=0$, and hence $\lim_{a\nearrow a^*}V_{\Omega}(x_a)=0$ in view of (\ref{New1}) and \eqref{ffin}, where $x_a$ is a global maximum point of $|u_a|$.
This then completes the proof of Lemma \ref{lem3.3}.
\qed

\subsection{$L^\infty-$uniform convergence as $a\nearrow a^*$}	

In this subsection, we first derive $L^\infty-$uniform convergence  of $ w_{a}(x)$ as $a\nearrow a^*$, based on which we finally finish the proof of Theorem \ref{thm1.2}.  	
	
\begin{prop}\label{4:prop}  Under the assumptions of Theorem~\ref{thm1.2}, let $u_{a}$ be a minimizer of $e_F(a)$, and consider the sequences $\{w_{a}\}$ and $\{x_{a}\}$ defined in Lemma \ref{lem3.3}. Then we have
\begin{enumerate}
\item [(i).] There exists a large constant  $R>0$ such that as $a\nearrow a^*$,
\begin{equation}\label{new4}
			|w_{a}(x)|\leq Ce^{-\frac{2}{3}|x|}\quad \hbox{in}\ \ \R^2 / B_R(0).
\end{equation}

\item [(ii).]  The global maximal point $x_a$ of $|u_a|$ is unique as  $a\nearrow a^*$, and $w_{a}(x)$ satisfies
\begin{equation} \label{NEW7}
		 w_{a}(x)\to \frac{w(x) }{\sqrt{a^*}}
\ \ \hbox{uniformly in $L^\infty (\R^2, \mathbb{C})$ as}\,\ a\nearrow a^*.
\end{equation}

\item [(iii).] The following estimate holds
			\begin{equation}\label{NEW17}
			\displaystyle  \Omega\inte x^\perp\Big(iw_{a}(x),\nabla w_{a}(x)\Big)dx=o(\eps_{a}^{1+\frac{\gamma}{2}})\ \ \mbox{as}\,\ a\nearrow a^*,
			\end{equation}
where $\eps_{a}=\eps_{a}>0$ is as in  Lemma \ref{lem3.3} and $\gamma =\min \{p, 2\}>0$.
\end{enumerate}
\end{prop}


\noindent{\bf Proof.} 1. Set $W_{a}=|w_{a}|^2$, so that $W_{a} $ satisfies (\ref{NEW10}). Recall from \eqref{NEW15} that $W_{a}\to \frac{1 }{a^*}w^2(|x|)$ strongly in $L^2(\R^2)$ as $a\nearrow a^*$.
Applying De Giorgi-Nash-Moser theory \cite[Theorem 4.1]{HL} to (\ref{NEW10}), then there exists a sufficiently large $R>0$ such that
\begin{equation}\label{new6}
	W_{a}(x)\leq \frac{1}{18a^*}\quad\hbox{in $ \R^2/ B_R(0)$ }
\end{equation}
holds uniformly as $a\nearrow a^*$.
By Lemma \ref{lem3.3}, we then derive from (\ref{NEW10}) that for sufficiently large $R>0$,
\begin{equation}\label{new6A}
	- \Delta W_{a}(x)+ \frac{16}{9}W_{a}(x)\leq 0 \quad\hbox{in $ \R^2/ B_R(0)$}
\end{equation}
where we used that $-\eps_{a}^2\mu _a\to 1$ as $a\nearrow a^*$.
Applying the comparison principle to (\ref{new6A}), the exponential decay \eqref{new4} is then proved.

2. We first claim that $w_a(x)$ converges to $\frac{w(x)}{\sqrt{a^*}}$ uniformly in $L^\infty(\R^2, \mathbb{C})$ as $a\nearrow a^*$. Indeed, by the exponential decay of (\ref{4:exp}) and \eqref{new4}, we only need to show the $L^\infty-$uniform convergence of $w_a(x)$ on any compact domain of $\R^2$ as $a\nearrow a^*$. Since $w_a(x)$ satisfies \eqref{NEW8}, denote
$$G_a(x):=-i\eps^2_{a}\Omega (x^\bot\cdot\nabla w_{a})-\Big[\frac{\eps^4_{a}\Omega^2}{4}|x|^2+\eps_{a}^2V_\Omega (\eps_{a} x+x_a)-\eps_{a}^2\mu_a-a |w_{a}|^2\Big] w_{a},$$
so that
\begin{equation}\label{rev-32}
-\Delta w_{a}(x)=G_a(x)\quad\hbox{in $H^1(\R^2, \mathbb{C})$.}
\end{equation}
Because $w_{a}(x)$ is bounded uniformly in $H^1(\R^2,\mathbb{C})$ as $a\nearrow a^*$, $G_a(x)$ is also bounded uniformly in $L^2_{loc}(\R^2, \mathbb{C})$ as $a\nearrow a^*$.
For any large $R>0$, it then follows from \cite[Theorem 8.8]{GT} that
\begin{equation}\label{rev-31}
\|w_a(x)\|_{H^2(B_R)}\leq C\Big(\|w_a(x)\|_{H^1(B_{R+1})}+\|G_a(x)\|_{L^2(B_{R+1})}\Big),
\end{equation}
where $C>0$ is independent of $a>0$ and $R>0$.
Therefore, $w_a(x)$ is also bounded uniformly in $H^2_{loc}(\R^2, \mathbb{C})$. Since the embedding $H^2(B_R)\hookrightarrow L^\infty(B_R)$ is compact, cf. \cite[Theorem 7.26]{GT}, we obtain that there exists a subsequence $\{w_{a_k}(x)\}$ of $\{w_a(x)\}$ such that
\[
\lim_{a_k\nearrow a^*}w_{a_k}(x)=w_0(x)\quad\hbox{uniformly in $L^\infty(B_R,\mathbb{C})$. }
\]
Since $R>0$ is arbitrary, we get from \eqref{3NEW7} that
\begin{equation}\label{RRev-1}
\lim_{a_k\nearrow a^*}w_{a_k}(x)=\frac{w}{\sqrt{a^*}}\quad\hbox{uniformly in $L^\infty_{loc}(\R^2,\mathbb{C})$. }
\end{equation}
Because the above convergence is independent of the subsequence that we choose, \eqref{RRev-1} holds essentially for the whole sequence. This further implies that  (\ref{NEW7}) holds
uniformly in $L^\infty (\R^2, \mathbb{C})$ as $a\nearrow a^*$.

We next prove the uniqueness of $x_a$ as $a\nearrow a^*$, where $x_a$ is a global maximum point of $|u_a|$. By \eqref{rev-31}, one can get that $|\nabla w_a|$ is bounded uniformly in $L^{q}_{loc}(\R^2)$ as $a\nearrow a^*$ for any $q\geq 2$. The $L^p$ estimate \cite[Theorem 9.11]{GT}   applied to (\ref{rev-32}) then implies that $w_a$ is bounded uniformly in $W^{2,q}_{loc}(\R^2)$ as $a\nearrow a^*$. The standard Sobolev embedding thus gives that $w_a$ is bounded uniformly in $C^{1,\alpha}_{loc}(\R^2)$. Furthermore, since $\eps_a^2V_{\Omega}(\eps_ax+x_a)\in C^{\alpha}_{loc}(\R^2)$, it follows from the Schauder estimate \cite[Theorem 6.2]{GT} that $w_a$ is bounded uniformly in $C^{2,\alpha}_{loc}(\R^2)$ as $a\nearrow a^*$. Hence, there exists $w_0\in C_{loc}^2(\R^2)$ such that
\[
w_a\to w_0\quad\hbox{in $C^2_{loc}(\R^2)$\,\ as \,\ $a\nearrow a^*$.}
\]
Noting from (\ref{NEW7}) that $w_a\to \frac{w}{\sqrt{a^*}}$ in $L^\infty(\R^2)$ as $a\nearrow a^*$, we conclude  that $w_0\equiv \frac{w}{\sqrt{a^*}}$, and hence
\[
|w_a|\to \frac{w}{\sqrt{a^*}}\quad\hbox{in $C^2_{loc}(\R^2)$\,\ as\,\  $a\nearrow a^*$.}
\]
Because the origin is the unique global maximum point of $w$, the above convergence shows that all global maximum points of $|w_a|$ must stay in a small ball $B_{\delta}(0)$ as $a\nearrow a^*$ for some small $\delta >0$. Since $w''(0)<0$, we conclude that $w''(r)<0$ for $0\leq r<\delta$. It then follows from \cite[Lemma 4.2]{NT} that as $a\nearrow a^*$,  each $|w_a|$ has a unique global maximum point which is just the origin. This further proves the uniqueness of global maximum points for $|u_a|$ as $a\nearrow a^*$.

3. By the definition of $e_F(a)$, we obtain from Lemma \ref{lem3.3} that
\begin{equation}\label{new7}\arraycolsep=1.5pt\begin{array}{lll}
e_F(a)
&=&F_{a}(u_{a})=
F_{a}\Big(\frac{1}{\eps_{a}}e^{i (\frac{ \Omega}{2}  x\cdot x_a^{\bot}-\theta)}w_{a}
\big(\frac{x-x_a}{\eps_{a}}\big)\Big)\\[3mm]
	&=&\displaystyle \frac{1}{\eps_{a} ^2}\inte |\nabla  w_{a}(x) |^2dx-
	\displaystyle \frac{a}{2\eps_{a} ^2} \inte
	| w_{a}(x)|^4dx+\displaystyle\frac{\eps_{a}^2\Omega^2}{4}\inte |x|^2| w_{a}(x)|^2dx\\[3mm]
	&&+\displaystyle   \inte V_\Omega (\eps_{a} x+x_a)
	| w_{a}(x)|^2dx-\displaystyle  \Omega\inte x^\perp\cdot\Big(i w_{a}(x),\nabla w_{a}(x)\Big)dx.
\end{array}
\end{equation}
Applying the convergence \eqref{NEW7} and the exponential decay \eqref{new4}, we have
\begin{equation*}
	\begin{split}
	\quad\Omega\inte x^\perp\cdot\Big(i w_{a}(x),\nabla  w_{a}(x)\Big)dx
	&=\Omega\inte x^\perp\cdot\Big(R_{a}\nabla I_{a}-I_{a} \nabla R_{a} \Big)dx\\
	&=2\Omega\inte x^\perp\cdot\Big(R_{a}\nabla I_{a}\Big)dx\leq C\|\nabla I_{a}\|_{L^2},
	\end{split}
\end{equation*}
where $w_{a}(x):=R_{a}(x)+iI_{a}(x)$ is defined after (\ref{orth con}).
We then derive from above that
\begin{equation}\label{new8}
	\eps_{a}^2 e_F(a) \geq \inte (|\nabla R_{a}|^2+|\nabla I_{a}|^2)dx-\displaystyle \frac{a^*}{2}\inte (R_{a}^4+I_{a}^4+2R^2_{a}I^2_{a})dx-C\eps_{a} ^2\|\nabla I_{a}\|_{L^2}.
\end{equation}
Since it follows from (\ref{NEW7}) that $R_{a}\to \frac{w}{\sqrt{a^*}}$ and $I_{a}\to 0$ uniformly in $\R^2$ as $a\nearrow a^*$, we now have
\begin{equation*}
	\begin{split}
	&\quad\inte |\nabla R_{a}|^2dx-\displaystyle \frac{a^*}{2}\inte \big(R_{a}^4+I_{a}^4+2R^2_{a}I^2_{a}\big)dx\\
	&\geq \displaystyle \inte |\nabla R_{a}|^2dx\Big(1-\inte |R_{a}|^2dx\Big)-\displaystyle \frac{a^*}{2}\inte \big(I_{a}^2+2R^2_{a}\big)I^2_{a}dx\\
	&=\big(1+o(1)\big)\inte I^2_{a}dx-\big(1+o(1)\big)\inte w^2I^2_{a}dx \ \ \mbox{as} \,\ a\nearrow a^*,
\end{split}
\end{equation*}
where the Gagliardo-Nirenberg inequality \eqref{GNineq} is used. Then, we derive from \eqref{new8} that
\begin{equation}\label{NEW20}
\begin{split}
	\eps_{a}^2 e_F(a)  &\geq  \inte |\nabla I_{a}|^2dx+\inte I^2_{a}dx-\inte w^2I^2_{a}dx-o(1)\inte I^2_{a}dx-C\eps_{a} ^2\|\nabla I_{a}\|_{L^2} \\
	& =\big(\mathcal{L}I_{a},I_{a}\big)-o(1)\inte I^2_{a}dx-C\eps_{a} ^2\|\nabla I_{a}\|_{L^2}
\quad \mbox{as} \,\ a\nearrow a^*,
	\end{split}
\end{equation}
where the operator $\mathcal{L}$ is defined by
\begin{equation}\label{new12}
	\mathcal{L}:=-\Delta-w^2+1 \quad\mbox{in} \,\ L^2(\R^2).
\end{equation}

Note from \cite[Corollary 11.9 and Theorem 11.8]{Lieb} that $(0, w)$ is  the first eigenpair of $\mathcal{L}$ and $ker(\mathcal{L})=\{w\}$. Moreover, since the essential spectrum of the operator $\mathcal{L}$ satisfies $\sigma_{ess}(\mathcal{L})=[1,+\infty)$, we have $0\in \sigma_d(\mathcal{L})$, where $\sigma_d$ denotes its discrete spectrum. Because $I_{a}$ is orthogonal to $w(x)$ by \eqref{new10},  we deduce that
\begin{equation}\label{new13}
	\big(\mathcal{L}I_{a},I_{a}\big)\geq \|I_{a}\|^2_{L^2(\R^2)}\quad \mbox{as} \,\ a\nearrow a^*.
\end{equation}
On the other hand, we obtain from (\ref{new12}) that
\begin{equation}\label{new14}
	(\mathcal{L}I_{a},I_{a})\geq \|\nabla I_{a}\|^2_{L^2(\R^2)}-\|w\|^2_{L^{\infty}}\|I_{a}\|^2_{L^2(\R^2)}.
\end{equation}
We thus deduce from \eqref{new13} and \eqref{new14} that there exists a constant $\rho>0$ such that  as $a\nearrow a^*$,
\begin{equation}\label{new15CC}
	(\mathcal{L}I_{a},I_{a})\geq \rho\|I_{a}\|^2_{H^1(\R^2)}.
	\end{equation}
Inserting \eqref{new15CC} into \eqref{NEW20}, we have
\begin{equation}\label{new123}
	\eps_{a}^2 e_F(a)\geq \frac{\rho}{2}\|I_{a}\|^2_{H^1(\R^2)}-C\eps_{a} ^2\|\nabla I_{a}\|_{L^2(\R^2)}
	\quad \mbox{as} \,\ a\nearrow a^*.
\end{equation}
Similar to (\ref{New1}), applying Lemma \ref{lem3.2} yields that for $\gamma =\min \{p, 2\}>0$,
\[
	\arraycolsep=1.5pt \begin {array}{lll}
	&\quad&M(a^*-a)^\frac{\gamma}{\gamma +2} \ge e_F(a)
	\\
&= &
	\displaystyle \inte |(\nabla -i\mathcal{A} )u_{a}|^2dx -
	\displaystyle\frac{{a}}{2} \inte |u_{a}|^4dx+ \inte  V_\Omega (x)|u_{a}|^2dx\\[3mm]
	&\ge & \displaystyle \inte \big| \nabla |u_{a}|\big|^2dx -
	\displaystyle\frac{{a}}{2} \inte |u_{a}|^4dx\geq \displaystyle \Big(1-\frac{{a}}{a^*}\Big) \inte \big| \nabla |u_{a}|\big|^2dx\\[3mm]
	&\geq& \displaystyle \frac 1 2\Big(1-\frac{{a}}{a^*}\Big) \eps_{a}^{-2}\quad \mbox{as} \,\ a\nearrow a^*.
\end{array}
\]
This further implies that
\begin{equation}\label{new1234}
	(a^*-a)^\frac{1}{\gamma +2}\leq C\e,\ \ e_F(a)\leq M(a^*-a)^\frac{\gamma}{\gamma +2}\leq C\eps_{a}^\gamma \quad\mbox{as} \,\ a\nearrow a^*.
\end{equation}
Following (\ref{new123}) and (\ref{new1234}), we obtain that
\begin{equation}\label{new15}
	\|I_{a}\|_{H^1(\R^2)}\leq C\eps_{a}^{1+\frac{\gamma}{2}}\quad \mbox{as} \,\ a\nearrow a^*.
	\end{equation}
Applying \eqref{new4} and \eqref{new15} now yields that
\begin{equation}\label{new16}
	\begin{split}
	&\displaystyle   \inte x^\perp\cdot \Big(i w_{a}(x),\nabla  w_{a}(x)\Big)dx\\
	=&2\inte x^\perp\cdot (R_{a}\nabla I_{a})dx=2\inte x^\perp\cdot \Big(\frac{w(x)}{\sqrt{a^*}}\nabla I_{a}\Big)dx+o(\eps_{a}^{1+\frac{\gamma}{2}})\\
	=&-2\inte x^\perp\cdot \Big(I_{a}\nabla \frac{w(x)}{\sqrt{a^*}} \Big)dx+o(\eps_{a}^{1+\frac{\gamma}{2}})=o(\eps_{a}^{1+\frac{\gamma}{2}}) \ \ \mbox{as} \,\ a\nearrow a^*,
	\end{split}
\end{equation}
where the convergence \eqref{NEW7} is also used.
The proof of Proposition \ref{4:prop} is therefore complete. \qed
\medskip


\vskip 0.05truein

\noindent{\bf Completion of the proof for Theorem \ref{thm1.2}.} In view of (\ref{3.9}) and (\ref{NEW7}), to establish Theorem \ref{thm1.2}, the rest is to prove that
\begin{equation}\label{3A:DDE}
	\eps_a:=\Big(\int_{\R^2}|\nabla u_{a}|^{2}dx\Big)^{-\frac{1}{2}}=\displaystyle
\frac{(a^*-a)^{\frac{1}{2+\gamma}}}{\lambda}+o\big((a^*-a)^{\frac{1}{2+\gamma}}\big)>0\ \ \mbox{as}\ \,a\nearrow a^*,
\end{equation}
where $\gamma =\min \{2, p\}>0$ and $\lam >0$ is as in (\ref{def:li}).

{In order to prove} (\ref{3A:DDE}), we take
\[u_\beta(x)=\frac{\beta }{\|w\|_2(a^*-a)^\frac{1}{2+\gamma}}
	w\Big(\frac{ \beta}{(a^*-a)^\frac{1}{2+\gamma}}x-y_0\Big)e^{\frac{i\Omega y_0^\perp\cdot x}{2\tau}},  \]
where $\beta \in (0,\infty)$ is to be determined later, as a test function of the energy $e_F(a)$ as $a\nearrow a^*$, and minimize it over $\beta >0$. By the similar calculations of Lemma \ref{lem3.2}, one then gets that
\begin{equation}\label{3:DDD}
	e_F(a)\le \min_{\beta \in (0,\infty)}F_a(u_\beta)=\Big( 1 + \frac 2 \gamma\Big)\frac{ \lambda^2}{ a^*}\big(a^*- a\big)^{\gamma /(2+\gamma )}\quad \mbox{as}\ \  a\nearrow a^*.
\end{equation}
On the other hand, following \eqref{new7} and \eqref{NEW17}, we derive from Lemma \ref{lem3.3} that
\begin{equation}\label{new18}
\begin{split}
	\quad e_F(a)&=F_{a}(u_{a})\\
	&=\displaystyle \frac{1}{\eps_{a} ^2}\Big[\inte \big|\nabla w_{a}(x)\big|^2dx-
	\displaystyle \frac{a^*}{2 } \inte
	|w_{a}(x)|^4dx\Big]+
	\displaystyle \frac{a^*-a}{2\eps_{a} ^2} \inte
	|w_{a}(x)|^4dx\\
	&\quad +\displaystyle\frac{\eps_{a}^2\Omega^2}{4}\inte |x|^2| w_{a}(x)|^2dx +\displaystyle   \inte V_{\Omega}(\eps_{a} x+x_a)
	| w_{a}(x)|^2dx\\
	&\quad-\Omega\inte x^\perp\Big(iw_{a}(x),\nabla w_{a}(x)\Big)dx.
	\end{split}
\end{equation}
The term in the square bracket is non-negative and can be dropped for a lower bound of $e_F(a)$. Moreover, we infer from \eqref{NEW17} that	 as $a\nearrow a^*$,
\[\Omega\inte x^\perp\Big(iw_{a}(x),\nabla w_{a}(x)\Big)dx= \left\{
	\arraycolsep=1.5pt \begin{array}{ll}
	&o(\e^p),\,\  \hbox{if}\ \ 0<p\leq 2; \\[2mm]
	&o(\e^2),\, \ \hbox{if}\ \ p> 2.
	\end{array}
	\right.\]
Because $w_{a}\to \frac{w}{\sqrt{a^*}}$ uniformly in $\R^2$ as $a\nearrow a^*$, we get from the exponential decay \eqref{new4} that
\[
	 \frac{a^*-a}{2\eps_{a} ^2} \inte
|w_{a}(x)|^4dx=\big[1+o(1)\big]\frac{a^*-a}{2(a^*)^2\eps_{a} ^2} \inte
|w|^4dx\ \ \mbox{as}\ \,a\nearrow a^*,
\]	
\[
\frac{\eps_{a}^2\Omega^2}{4}\inte |x|^2| w_{a}(x)|^2dx=\big[1+o(1)\big]\frac{\eps_{a}^2\Omega^2}{4a^*}\inte |x|^2| w(x)|^2dx\ \ \mbox{as}\ \,a\nearrow a^*,
\]
and
\[
\begin{aligned}
&\quad\inte V_{\Omega}(\eps_{a} x+x_a)
| w_{a}(x)|^2dx\\
&=\int_{B_{\frac{1}{\sqrt{\e}}}}V_{\Omega}(\eps_{a} x+x_a)
| w_{a}(x)|^2dx+\int_{\R^2/B_{\frac{1}{\sqrt{\e}}}} V_{\Omega}(\eps_{a} x+x_a)
| w_{a}(x)|^2dx\\
&=\eps_{a}^p\int_{B_{\frac{1}{\sqrt{\e}}}}h( x+\e^{-1}x_a)
| w_{a}(x)|^2dx+o(\e^p)\\
&=\big[1+o(1)\big]\frac{\e^p}{a^*}\inte h(x+\e^{-1}x_a)
| w(x)|^2dx\ \ \mbox{as}\ \,a\nearrow a^*.
\end{aligned}
\]	
Note that	for $p=2$,
\begin{equation}\label{3T:DDE}
\begin{aligned}
&\frac{a^*-a}{2(a^*)^2\eps_{a} ^2} \inte
 |w|^4dx+\frac{\eps_{a}^2\Omega^2}{4a^*}\inte |x|^2| w(x)|^2dx+\frac{\e^p}{a^*}\inte h(x+\e^{-1}x_a)
 | w(x)|^2dx\\
 \geq &\,\Big( 1 + \frac 2 \gamma\Big)\frac{ \lambda^2}{ a^*}(a^*- a)^{\gamma /(2+\gamma )}\ \ \mbox{as}\ \,a\nearrow a^*,
\end{aligned}
\end{equation}
where the identity holds if and only if (\ref{3A:DDE}) holds true. Applying (\ref{3:DDD}), we then conclude from above that for $p=2$,
\[
e_F(a) \approx \Big( 1 + \frac 2 \gamma\Big)\frac{ \lambda^2}{ a^*}(a^*- a)^{\gamma /(2+\gamma )}\ \ \mbox{as}\ \,a\nearrow a^*,
\]
and (\ref{3A:DDE}) is hence true for $p=2$. Similarly, one can also obtain from above that (\ref{3A:DDE}) holds true for the cases where $p\not =2$. This completes the proof of Theorem \ref{thm1.2}.   \qed

The proof of Theorem \ref{thm1.2} implies, see those around (\ref{3T:DDE}), that for all $p>0$,
\begin{equation}\label{3:DDE}
	\eps_a:=\Big(\int_{\R^2}|\nabla u_{a}|^{2}dx\Big)^{-\frac{1}{2}}=\displaystyle
\frac{(a^*-a)^{\frac{1}{2+\gamma}}}{\lambda}+o\big((a^*-a)^{\frac{1}{2+\gamma}}\big)>0\ \ \mbox{as}\ \,a\nearrow a^*,
\end{equation}
\begin{equation}\label{4:DDE}
	0<\displaystyle \frac{1}{\eps_{a} ^2}\Big[\inte \big|\nabla w_{a}(x)\big|^2dx-
	\displaystyle \frac{a^*}{2 } \inte
	|w_{a}(x)|^4dx\Big]=o(\e^{\gamma})\ \ \mbox{as}\ \,a\nearrow a^*,
\end{equation}
and
\begin{equation}\label{5:DDE}
	\lim_{a\nearrow a^*}\frac{x_a}{\e}=y_0,\ \ \hbox{if}\,\ 0<p\leq 2,
\end{equation}
where $y_0$  is a global minimum point of $H(y)=\inte h(x+y)w^2(x)dx$ defined in the assumption $(V)$.

\section{Refined Limiting Profiles}
In this section we always assume $0<\Omega<\Omega^*$, where $\Omega ^*>0$ is
defined as in (\ref{Omega}), and suppose $V(x)\in L^\infty_{\rm loc}(\R^2)$ satisfying
(\ref{A:V}) and
\begin{enumerate}
	\item [\rm($V_1$).]
	$V_{\Omega}(x)=V(x)-\frac{\Omega^2}{4}|x|^2\geq0$,  $V_{\Omega}(x)\in C^2(\R^2)$ is a homogeneous function of degree $2$ and $H(y)=\inte V_{\Omega}(x+y)w^2(x)dx$ admits a
	unique critical point $y_0\in \R^2$.
\end{enumerate}
It is not difficult to get from the above assumption that Theorem \ref{thm1.2} is applicable as $a\nearrow a^*$.
As before, let $u_a$ be a minimizer of $e_F(a)$ as $a\nearrow a^*$, and set
\begin{equation}
\arraycolsep=1.5pt
 \begin{array}{lll}
 w_a(x)&: =& \displaystyle\eps_a u_a\Big(\eps_ax  +x_a\Big)e^{-i \big(
 \frac{\eps_a\Omega}{2 } x\cdot x_a^\perp-\theta_a\big)} \\[3mm]
 &:=&R_a(x)+iI_a(x)=\displaystyle\frac{w(x)}{\sqrt{a^*}}+ \hat{R}_a(x)+iI_a(x),\ \ \eps_a:=\Big(\int_{\R^2}|\nabla u_a|^{2}dx\Big)^{-\frac{1}{2}}>0,
 \end{array}
\label{4:1}
\end{equation}
where $\theta_a\in [0,2\pi)$ is chosen properly such that \eqref{orth con} holds,  and $x_a\in\R^2$ denotes the unique global maximal point of $|u_a|$ as $a\nearrow a^*$. Here $R_a ,\ \hat{R}_a$ and $  I_a $ are real-valued functions in $H^1(\R^2)$.
Following the proof of Theorem \ref{thm1.2}, we reduce from (\ref{new10}) and (\ref{NEW7}) that
\begin{equation}
\inte wI_adx\equiv 0,\ \ \hat{R}_a\to 0 \ \ \mbox{and}\ \ I_a\to 0  \ \mbox{uniformly in}\ \,  \R^2  \ \mbox{as}\ \ a\nearrow a^*.
\label{4:2}
\end{equation}
Based on (\ref{4:2}), {\em the main purpose of this section} is to establish the following refined estimates of $\hat{R}_a$ and $I_a$ as $a\nearrow a^*$.

 \begin{prop}\label{prop4.3}
Suppose $V(x)\in L^\infty_{\rm loc}(\R^2)$ satisfies (\ref{A:V}) and $(V_1)$, and assume $0<\Omega<\Omega^*$, where $\Omega ^*>0$ is defined as in (\ref{Omega}). Suppose $u_a$ is a minimizer of $e_F(a)$, and let $\hat{R}_a$ and $I_a$ be defined by (\ref{4:1}). Then there exists a constant $C>0$, independent of $0<a<a^*$, such that as $a\nearrow a^*$,	
\begin{itemize}
\item [(i).] $|\hat{R}_a(x)|, \, |\nabla \hat{R}_a(x)|\leq C\e^{4}e^{-\frac 12|x|}\,\ \hbox{in\,  $\R^2$,}$
\item [(ii).] $|I_a(x)|, \, |\nabla I_a(x)|\leq C\e^{6}e^{-\frac 18|x|}\,\ \hbox{in\,  $\R^2$,}$	
\end{itemize}
where $\e >0$ is as in (\ref{4:1}).
\end{prop}

In order to prove Proposition \ref{prop4.3}, we first give some observations. Because $\nabla |w_a(0)|\equiv 0$ holds for all $a\in(0,a^*)$, we calculate from (\ref{4:1}) that
\begin{equation}\label{4:3RT}
\nabla \hat{R}_a(0)=-\frac{I_a(0)\nabla I_a(0)}{\frac{w(0)}{\sqrt{a^*}}+\hat{R}_a(0)}\to 0   \ \ \mbox{uniformly in}\ \,  \R^2  \ \mbox{as}\ \ a\nearrow a^*.
\end{equation}
Recall from (\ref{NEW8}) that $w_a$ satisfies the following equation
\begin{equation}\label{4:weqn}
\Delta w_{a}-i\eps^2_a\Omega (x^\bot\cdot\nabla w_{a})-\Big[\frac{\eps_a^4\Omega^2}{4}|x|^2+\eps_a^2V_\Omega (\eps_a x+x_a)-\mu _a\eps_a^2-a |w_{a}|^2\Big] w_{a}=0
\ \ \mbox{in}\,\ \R^2,
\end{equation}
where $\eps_a>0$ is as in (\ref{4:1}), and $\mu _a\in \R$ is a Lagrange multiplier. Note from Lemma \ref{lem3.3} that $\mu _a\in \R$ satisfies $\mu _a\eps_a^2\to -1$ as $a\nearrow a^*$.
For simplicity, we denote the operator $\mathcal{L}_a$ by
\begin{equation}\label{D6}
\mathcal{L}_a:=-\Delta+\Big(\frac{\eps_a^4\Omega^2}{4}|x|^2+\eps_a^2V_\Omega (\eps_a x+x_a)-\e^2\mu _a-a|w_a|^2\Big).
\end{equation}
It then follows from (\ref{4:weqn}) that $(R_a, I_a)$ satisfies the following system
\begin{equation}\label{4:3}
\left\{
\begin{aligned}
\mathcal{L}_a R_a&=\eps^2_a\Omega (x^\bot\cdot\nabla I_a)
\ \ \mbox{in}\,\  \R^2,\\
\mathcal{L}_a I_a&=-\eps^2_a\Omega (x^\bot\cdot\nabla R_a)
\ \ \mbox{in}\,\  \R^2.\\
\end{aligned}
\right.
\end{equation}
By applying Lemma \ref{rev-7}, we start with the following  estimates:

\begin{lem}\label{lem4.1} Under the assumptions of Proposition \ref{prop4.3},  we have
\begin{enumerate}
\item There exists a constant $C>0$, independent of $0<a<a^*$, such that as $a\nearrow a^*$,
\begin{equation}\label{lem4.1:A}
  |\nabla  {R}_a (x)|, \ |\nabla  I_a (x)|\leq Ce^{-\frac{2}{3}|x|}\quad \hbox{in}\, \ \R^2,
\end{equation}
where ${R}_a$ and $I_a$ are defined by (\ref{4:1}).
\item The Lagrange multiplier $\mu _a$ of (\ref{4:weqn}) satisfies
\begin{equation}\label{lem4.1:B}
\big|1+\mu _a\eps_a^2\big|\leq C\eps^{4}_a
\ \ \mbox{as}\,\ a\nearrow a^*,
\end{equation}
where $\eps_a>0$ is as in (\ref{4:1}) and $C>0$ is independent of $0<a<a^*$.
\end{enumerate}
\end{lem}

\noindent{\bf Proof.} 1. We first claim that
\begin{equation}\label{rev-24}
\|\nabla w_a\|_{L^\infty(\R^2)}\leq C\,\ \hbox{and}\,\ \lim_{|x|\to\infty}|\nabla w_a(x)|=0 \,\ \hbox{as}\,\  a\nearrow a^*,
\end{equation}
where $C>0$ is independent of $0<a<a^*$. To prove the above claim, denote
\begin{equation}\label{rev-19}
\hat{V}_a(x)=\frac{\eps_a^4\Omega^2}{4}|x|^2+\eps_a^2V_\Omega (\eps_a x+x_a)-\mu _a\eps_a^2-a |w_{a}|^2,
\end{equation}
where $\mu _a\eps_a^2\to -1$ as $a\nearrow a^*$.
We then derive from \eqref{4:weqn} that $w_a$ satisfies
\begin{equation}\label{rev-20}
\big[-\Delta+\hat{V}_a(x)\big]w_a+i\e^2\Omega(x^\perp\cdot\nabla w_a)=0\,\ \hbox{in}\,\ \R^2.
\end{equation}
Following \eqref{rev-20}, we obtain that for $i\neq j$,
\[
\big[-\Delta+\hat{V}_a(x)\big]\partial_i w_a+\partial_i\hat{V}_a(x)w_a+i\e^2\Omega(x^\perp\cdot\nabla \partial_i w_a)+(-1)^{i+1}i\e^2\Omega\partial_jw_a=0\,\ \hbox{in}\,\ \R^2,
\]
which further implies that for $i\neq j$,
\begin{equation}\label{rev-21}
\begin{split}
&-\frac 12 \Delta|\partial_iw_a|^2+\hat{V}_a(x)|\partial_iw_a|^2+|\nabla \partial_iw_a|^2 +\partial_i\hat{V}_a(x)(w_a,\partial_iw_a)\\
&-\e^2\Omega x^\perp\cdot (i\partial_i w_a,\nabla \partial_i w_a)+(-1)^{i+1}\e^2\Omega(i\partial_jw_a,\partial_i w_a)=0\,\ \hbox{in}\,\ \R^2.
\end{split}
\end{equation}
Note from  (\ref{new4}) and (\ref{rev-19}) that for $i\neq j$,
\[
\e^2\Omega x^\perp\cdot (i\partial_i w_a,\nabla \partial_i w_a)\leq \frac{\eps_a^4\Omega^2}{4}|x|^2|\partial_i w_a|^2+|\nabla \partial_iw_a|^2,
\]
\[
\partial_i\hat{V}_a(x)(w_a,\partial_iw_a)\leq |\partial_i\hat{V}_a(x)|^2|w_a|^2+\frac 14|\partial_iw_a|^2\leq C\e^4e^{-|x|}+8a^2|w_a|^4|\partial_iw_a|^2+\frac 14|\partial_iw_a|^2,
\]
\[
\big|\e^2\Omega(i\partial_jw_a,\partial_i w_a)\big|\leq \frac{\e^2\Omega}{2}\big[(\partial_i w_a)^2+(\partial_jw_a)^2\big].
\]
It then follows from \eqref{rev-21} that as $a\nearrow a^*$,
\begin{equation}\label{rev-22}
\Big(-\frac 12\Delta-a|w_a|^2-8a^2|w_a|^4\Big) |\nabla w_a|^2\leq C\e^4 e^{-|x|}\ \ \hbox{in}\ \ \R^2,
\end{equation}
where the constant $C>0$ is independent of $0<a<a^*.$
Applying De Giorgi-Nash-Moser theory \cite[Theorem 4.1]{HL}, we obtain from \eqref{rev-22} that for any $y\in \R^2$,
\begin{equation}\label{rev-23}
\sup_{x\in B_{\frac 12}(y)}|\nabla w_a(x)|^2\leq C\Big(\|\nabla w_a\|^2_{L^2(B_{1}(y))}+\e^4 \big\|e^{-|x|}\big\|_{L^2(B_{1}(y))}\Big)
\end{equation}
holds  as $a\nearrow a^*$, where the constant $C>0$ is independent of $0<a<a^*.$
Since Theorem \ref{thm1.2} gives that $w_a\to \frac{w}{\sqrt{a^*}}$ strongly in $H^1(\R^2,\mathbb{C})$ as $a\nearrow a^*$, the claim \eqref{rev-24} is now proved in view of \eqref{rev-23}.

We next prove the estimate \eqref{lem4.1:A}. Since $(R_a,I_a)$ solves the system (\ref{4:3}),  $(R_a,I_a)$ satisfies the
system \eqref{rev-2} with
\[
V_1(x)=V_2(x)=\frac{\eps_a^4\Omega^2}{4}|x|^2+\eps_a^2V_\Omega (\eps_a x+x_a)-\e^2\mu _a,\, \,\delta=\frac {1}{18},
\]
\[
b_1=-b_2=\eps^2_a\Omega,\ \ f_1(x)=a|w_a|^2R_a,\ \ f_2(x)=a|w_a|^2I_a.
\]
We derive from (\ref{new4}) that as $a\nearrow a^*$,
 \[
 |f_1(x)|,\ \ |f_2(x)|\le C_1e^{-\frac{2|x|}{3}} \ \ \mbox{in}\ \ \R^2,
 \]
where $C_1>0$ is independent of $0<a<a^*$.
Together with (\ref{rev-24}), the above estimate of $f_1(x)$ and $ f_2(x)$ then yields that as $a\nearrow a^*$, the assumption (\ref{rev-9}) of Lemma \ref{rev-7} is satisfied for some $P>0$ and $Q=\frac 23$.
Since $ V_{\Omega}(x)$ is a homogeneous function of degree $2$, we have $2V_{\Omega}(x)=x\cdot \nabla V_{\Omega}(x)$, which implies that
$$|\nabla V_{\Omega}(x)|\geq \frac{2V_{\Omega}(x)}{|x|}\geq 2|x|\min_{x\in\R^2}V_{\Omega}\Big(\frac{x}{|x|}\Big):=C_1|x|.$$
On the other hand, we also have
$$\nabla V_{\Omega}(x)=\nabla \Big[|x|^2V_{\Omega}\Big(\frac{x}{|x|}\Big)\Big]=2V_{\Omega}\Big(\frac{x}{|x|}\Big)x+\frac{\nabla V_{\Omega}(y)\cdot x^\perp}{|x|}x^\perp\Big|_{y=\frac{x}{|x|}}, $$
which implies that
\[
|\nabla V_{\Omega}(x)|\leq C_2|x|, \ \ C_2:=3\|V_{\Omega}(x)\|_{C^1(\partial B_1)}>0.
\]
We thus obtain from above that
there exist positive constants $C_1>0$ and $C_2>0$ such that
\begin{equation}\label{D21}
C_1|x|\leq |\nabla V_{\Omega}(x)|\leq C_2|x|\,\ \hbox{in}\,\ \R^2.
\end{equation}
Applying (\ref{D21}), one can check that as $a\nearrow a^*$, $V_1(x)$ and $V_2(x)$ satisfy the assumption \eqref{rev-10} of Lemma \ref{rev-7} for $\delta=\frac {1}{18}$. Therefore,  the estimate
\eqref{lem4.1:A} is now proved by applying Lemma \ref{rev-7}.

2. Denoting $\hat{u}_a(x):=\e u_a(\e x+x_a) $,
we have $\|\nabla \hat{u}_a\|_2^2=1$ and
\begin{equation}\label{revise 3}
\begin{split}
\inte |\nabla w_a(x)|^2dx&=\inte \Big|\nabla \Big(\hat{u}_a(x)e^{-i\big(\frac{\Omega\e x\cdot x_a^{\perp}}{2}-\theta_a\big)}\Big)\Big|^2dx\\
&=\inte \Big|\nabla \hat{u}_a(x)e^{-i\big(\frac{\Omega\e x\cdot x_a^{\perp}}{2}-\theta_a\big)}-\frac{i\Omega\e x_a^{\perp}}{2}\hat{u}_a(x)e^{-i\big(\frac{\Omega\e x\cdot x_a^{\perp}}{2}-\theta_a\big)} \Big|^2dx\\
&=\inte \Big|\nabla \hat{u}_a(x)-\frac{i\Omega\e x_a^{\perp}}{2}\hat{u}_a(x) \Big|^2dx\\
&=\inte \Big\{|\nabla \hat{u}_a|^2+\frac{\e^2\Omega^2|x_a|^2}{4}|\hat{u}_a|^2-\Omega\e x_a^{\perp}\cdot (i\hat{u}_a, \nabla \hat{u}_a)\Big\}dx \\
&=\inte |\nabla \hat{u}_a|^2dx-III=1-III.
\end{split}
\end{equation}
We obtain from \eqref{new15} and \eqref{5:DDE} that as $a\nearrow a^*$,
\begin{equation}\label{revise 5}
\begin{split}
|III|:&=\Big|\inte \Big\{\frac{\e^2\Omega^2|x_a|^2}{4}|\hat{u}_a|^2-\Omega\e x_a^{\perp}\cdot (i\hat{u}_a, \nabla \hat{u}_a)\Big\}dx\Big|\\
&\leq C\e^4+\e\Big|\inte \Omega x_a^{\perp}\cdot (i\hat{u}_a, \nabla \hat{u}_a)dx\Big|\\
&\leq C\e^4+\e\Big|\inte\Omega x_a^{\perp}\cdot (iw_a, \nabla w_a)dx\Big|+ \frac{\Omega^2|x_a|^2\e^2}{2}\inte |w_a|^2dx \leq C\e^4.
\end{split}
\end{equation}
On the other hand, we infer from \eqref{3:mu} and \eqref{3.9} that
\begin{equation}\label{revise 4}
\begin{split}
\e^2\mu_a=\e^2e_F(a)-\frac{a}{2}\e^2\inte |u_a|^4dx=\e^2e_F(a)-\frac{a}{2}\inte |w_a|^4dx.\\
\end{split}
\end{equation}
We thus derive from \eqref{revise 3}-\eqref{revise 4} and \eqref{3:DDD} that as $a\nearrow a^*$,
\begin{equation}\label{B1}
\begin{split}
\quad |1+\varepsilon^2_a\mu _a|
&=\Big|\inte |\nabla w_a|^2+III+\e^2e_F(a)-\frac{a}{2}\inte |w_a|^4dx\Big|\\
&\leq \Big|\inte |\nabla w_a|^2-\frac{a^*}{2}\inte |w_a|^4dx+\frac{a^*-a}{2}\inte |w_a|^4dx\Big|+C\e^4 \leq C\e^4,
\end{split}
\end{equation}
where the estimate (\ref{4:DDE}) is also used in the last inequality. Therefore, the estimate \eqref{lem4.1:B} is proved. This completes the proof of Lemma \ref{lem4.1}.\qed
\medskip

The proof of Proposition \ref{prop4.3} also needs the following iteration results.

\begin{lem}\label{4:lem4.2}
Under the assumptions of Proposition \ref{prop4.3}, for any $\alpha\geq0$ suppose there exists a constant $C(\alpha)>0$, independent of $0<a<a^*$, such that as $a\nearrow a^*$,
\begin{equation}\label{prop4.1:A}
|x^{\perp}\cdot\nabla  {R}_a (x)|\leq C(\alpha)\e^{\alpha}e^{-\frac{1}{4}|x|}\ \ \mbox{in}\ \,\R^2.
\end{equation}
Then as $a\nearrow a^*$,
\begin{equation}\label{prop4.1:B}
|I_a(x)|,\ \, |\nabla I_a(x)|\leq CC(\alpha)\e^{2+\alpha}e^{-\frac{1}{8}|x|}\ \ \mbox{in}\, \ \R^2,
\end{equation}
where the constant $C>0$ is independent of $0<a<a^*$ and $\alpha \ge 0$.
\end{lem}

{\noindent \bf Proof.} Recall from (\ref{4:2}) and (\ref{4:3}) that $I_a$ satisfies
\begin{equation}\label{D1*}
\mathcal{L}_aI_a(x)=-\e^2\Omega\big(x^{\perp}\cdot \nabla R_a\big)\quad\hbox{in}\,\ \R^2,\quad \inte I_a(x)w(x)dx= 0.
\end{equation}
Multiplying \eqref{D1*} by $I_a$ and integrating over $\R^2$, we obtain from (\ref{prop4.1:A}) that
\begin{equation}\label{D10}
\inte (\mathcal{L}_aI_a)I_adx=-\eps^2_a\Omega \inte(x^\bot\cdot\nabla {R}_a)I_adx\leq CC(\alpha)\e^{2+\alpha}\|I_a\|_{L^2(\R^2)}.
\end{equation}
Following \eqref{new15CC}, we also get that
\begin{equation}\label{D11}
\begin{split}
\inte (\mathcal{L}_aI_a)I_adx&\ge  \inte \Big[(\mathcal{L}I_a)I_a- \big(1+\varepsilon^2_a\mu _a\big)I_a^2- \big(a|w_a|^2-w^2\big)I_a^2\Big]dx\\&
= \inte (\mathcal{L}I_a)I_adx+o(1)\inte I_a^2dx\geq \frac{\rho}{2}\|I_a\|^2_{H^1(\R^2)},
\end{split}\end{equation}
where the constant $\rho>0$ given in \eqref{new15CC} is independent of $0<a<a^*$.
Thus, we infer from \eqref{D10} and \eqref{D11} that
\begin{equation}\label{D17}
\|I_a\|_{L^2(\R^2)}\leq \|I_a\|_{H^1(\R^2)}\leq CC(\alpha)\e^{2+\alpha}.
\end{equation}

On the other hand, we derive from \eqref{D1*} that $|I_a|^2$ satisfies the following equation
\[
\begin{aligned}
&\Big[-\frac 12\Delta +\Big(\frac{\eps_a^4\Omega^2}{4}|x|^2+\eps_a^2V_\Omega (\eps_a x+x_a)-\mu _a\eps_a^2-a |w_{a}|^2\Big) \Big]|I_{a}|^2+|\nabla I_a|^2\\
&=-\eps^2_a\Omega (x^\bot\cdot\nabla  {R}_{a})I_a\,\ \hbox{in}\,\ \R^2,
\end{aligned}
\]
which implies that
\begin{equation}\label{D20}
-\frac 12\Delta|I_{a}|^2-\mu _a\eps_a^2|I_{a}|^2-a |w_{a}|^2|I_{a}|^2\leq -\eps^2_a\Omega (x^\bot\cdot\nabla  {R}_{a})I_a\,\ \hbox{in}\,\ \R^2.
\end{equation}
Since $\mu _a\eps_a^2\to -1$ as $a\nearrow a^*$, by De Giorgi-Nash-Moser theory \cite[Theorem 4.1]{HL}, we obtain from \eqref{D20} that for any $y\in \R^2$,
\begin{equation}\label{D-21}
\sup_{x\in B_{\frac 12}(y)}|I_a(x)|^2\leq C\Big(\|I_a\|^2_{L^2(B_1(y))}+\|\eps^2_a\Omega (x^\bot\cdot\nabla {R}_{a})I_a\|_{L^2(B_1(y))}\Big).
\end{equation}
We hence deduce from (\ref{new4}), \eqref{D17} and \eqref{D-21} that
\begin{equation}\label{rev-28}
\|I_a\|_{L^\infty(\R^2)}\leq  CC(\alpha)\e^{2+\alpha},
\end{equation}
and
\begin{equation}\label{rev-29}
\Big|a|w_a|^2I_a-\eps^2_a\Omega (x^\bot\cdot\nabla {R}_{a})\Big|\leq CC(\alpha)\e^{2+\alpha}e^{-\frac{|x|}{4}}\,\ \hbox{in\,\  $\R^2$},
\end{equation}
where the assumption (\ref{prop4.1:A}) is used.
Applying the comparison principle to \eqref{4:3}, we derive from \eqref{rev-28} and \eqref{rev-29} that
\[
|I_{a}(x)|\leq CC(\alpha)\e^{2+\alpha}e^{-\frac{|x|}{4}}\,\ \hbox{in}\,\ \R^2.
\]
Furthermore, applying gradient estimates (see (3.15) in \cite{GT}) to the equation \eqref{D1*}, we conclude from above that
\begin{equation}\label{D2}
|\nabla I_a(x)|\leq CC(\alpha)\e^{2+\alpha}e^{-\frac 18|x|}\ \ \hbox{in\,  $\R^2$,}
\end{equation}
which therefore completes the proof of Lemma \ref{4:lem4.2}. \qed
\medskip



\vskip 0.05truein

\noindent{\bf Proof of Proposition \ref{prop4.3}. }Note from (\ref{4:1}) and \eqref{4:3} that $(\hat {R}_a, I_a)$ satisfies the following system
\begin{equation}\label{D1}
\left\{
\begin{aligned}
\big(\mathcal{L}_a-w^2-\sqrt{a^*}R_aw\big)\hat{R}_a&=F_a(x)+\e^2\Omega\big(x^{\perp}\cdot \nabla I_a\big)\ \ \mbox{in}\,\  \R^2,\\
\mathcal{L}_a I_a&=-\eps^2_a\Omega (x^\bot\cdot\nabla \hat{R}_a)
\ \ \mbox{in}\,\  \R^2,\\
\end{aligned}
\right.
\end{equation}
where $F_a(x)$ is defined by
\begin{equation}\label{D9}
F_a(x):=-\Big[\frac{\eps_a^4\Omega^2}{4}|x|^2+\eps_a^2V_\Omega (\eps_a x+x_a)-\e^2\mu _a-1-a|I_a|^2+(a^*-a)R^2_a\Big]\frac{w}{\sqrt{a^*}}.
\end{equation}
Next, we shall divide the proof into three steps.

\vskip 0.05truein

{\em Step 1.} We first prove the following $C^1$ estimates  of $\hat{R}_a$ and $I_a$:
\begin{equation}\label{D22}
\| \hat{R}_a\|_{C^{1}(\R^2)}, \,\|I_a\|_{C^{1}(\R^2)}\leq C\e^{2}\ \ \mbox{as}\ \ a\nearrow a^*,
\end{equation}
where $C>0$ is independent of $0<a<a^*$.

Recall from Lemmas \ref{lem4.1} and  \ref{4:lem4.2} that as $a\nearrow a^*$,
\begin{equation}\label{D19}
| I_a(x)|,\ |\nabla I_a(x)|\leq C\e^2e^{-\frac 18|x|}\quad\hbox{in}\,\ \R^2.
\end{equation}
We claim that there exists $C>0$, independent of $0<a<a^*$, such that as $a\nearrow a^*$,
\begin{equation}\label{D7}
\|\hat{R}_a\|_{L^{\infty}(\R^2)}\leq C\e^2.
\end{equation}
On the contrary, suppose that the above claim is false, i.e.,
$\lim_{a\nearrow a^*}\frac{\|\hat{R}_a\|_{L^{\infty}(\R^2)}}{\e^2}=\infty.$
Set $\bar {R}_a=\frac{\hat {R}_a}{\|\hat {R}_a\|_{L^{\infty}(\R^2)}}$, then we get from
\eqref{D1} that $\bar {R}_a$ satisfies
\begin{equation}\label{D14}
\begin{split}
&\quad\Big(\mathcal {L}_a-w^2-\sqrt{a^*}R_aw\Big)\bar{R}_a
=\frac{F_a(x)+\e^2\Omega\big(x^{\perp}\cdot \nabla I_a\big)}{\|\hat {R}_a\|_{L^{\infty}(\R^2)}}\quad\hbox{in}\,\ \R^2.
\end{split}
\end{equation}
The estimate \eqref{D19} and  Lemma \ref{lem4.1} imply that the right hand side of \eqref{D14} satisfies
\[
\Big|\frac{F_a(x)+\e^2\Omega\big(x^{\perp}\cdot \nabla I_a\big)}{\|\hat {R}_a\|_{L^{\infty}(\R^2)}}\Big|\leq e^{-\frac {1}{16}|x|}\quad\hbox{in}\,\ \R^2.
\]
By the comparison principle, we then obtain from \eqref{D14} that
\[
|\bar {R}_a(x)|\leq e^{\frac {R}{16}-\frac {1}{16}|x|}\quad\hbox{in}\,\ \R^2/B_R.
\]
Let $z_a\in\R^2$ be a global maximal point of $\bar{R}_a$ in $\R^2$.
We then have $|z_a|\leq 2R$ uniformly in $0<a<a^*$ for sufficiently large $R>0$. Applying the elliptic regularity theory, there exists a subsequence, still denoted by $\{\bar{R}_a\}$, such that
$\bar{R}_a(x)\to {R}_0(x)$ in $C^1_{loc}(\R^2)$ as $a\nearrow a^*$. Hence, ${R}_0$ satisfies
\[
-\Delta {R}_0+(1-3w^2){R}_0=0\quad\hbox{in}\,\ \R^2,
\]
which implies that ${R}_0=\sum_{i=1}^2c_i\frac{\partial w}{\partial x_i}$. On the other hand, by \eqref{4:3RT} and \eqref{D19}, we have
\[
\big|\nabla {R}_0(0)\big|=\lim_{a\nearrow a^*}\Bigg|\frac{-I_a(0)\nabla I_a(0)}{\big[\frac{w(0)}{\sqrt{a^*}}+\hat{R}_a(0)\big]\|\hat{R}_a\|_{L^{\infty}(\R^2)}}\Bigg|=0.
\]
Because $\nabla {R}_0(0)=0$, by the non-degeneracy of $w(x)$ at the origin, we get that $c_1=c_2=0$ and hence ${R}_0\equiv 0$, which however contradicts to the fact that $1=\bar{R}_a(z_a)\to {R}_0(z_0)=1$ as $a\nearrow a^*$ for some $z_0\in \R^2$. Thus, the claim \eqref{D7} holds true.

Applying (\ref{lem4.1:B}), \eqref{D19} and \eqref{D7}, by  the comparison principle we derive from \eqref{D1} that $\hat{R}_a$ satisfies
\begin{equation}\label{D12}
|\hat{R}_a(x)|\leq C\e^{2}e^{-\frac {1}{16}|x|}\quad\hbox{in}\,\ \R^2.
\end{equation}
Furthermore, using the gradient estimates (see (3.15) in \cite{GT}) again, it follows from \eqref{D1} and (\ref{D12}) that
\begin{equation}\label{D8}
|\nabla \hat{R}_a|\leq C\e^{2}e^{-\frac {1}{32}|x|}\quad\hbox{in}\,\ \R^2.
\end{equation}
Combining \eqref{D19}, \eqref{D7} and \eqref{D8}, we therefore conclude that \eqref{D22} holds true.

\vskip 0.05truein
{\em Step 2.} We next prove that  as $a\nearrow a^*$,
\begin{equation}\label{D16}
|\hat{R}_a(x)|,\  |\nabla \hat{R}_a(x)|,\,|I_a(x)|,\  |\nabla I_a(x)|\leq C\e^{2}e^{-\frac 12 |x|}\  \, \hbox{in}\,\ \R^2,
\end{equation}
where the constant $C >0$ is independent of $0<a<a^*$.

Indeed, note from the system \eqref{D1} that $(\hat{R}_a,I_a)$ satisfies the system \eqref{rev-2} with
\[
V_1(x)=V_2(x)=\frac{\eps_a^4\Omega^2}{4}|x|^2+\eps_a^2V_\Omega (\eps_a x+x_a)-\e^2\mu _a,
\]
\[
f_1(x)=\big(w^2+\sqrt{a^*}R_aw\big)\hat{R}_a+F_a(x)+a|w_a|^2\hat{R}_a,\,f_2(x)=a|w_a|^2I_a,
\]
\[
b_1=-b_2=\eps^2_a\Omega,\,P=C\e^2,\,Q=\frac 12,\,\delta=\frac 14,
\]
where $F_a(x)$ is defined by (\ref{D9}).
By Step 1, one can check from \eqref{D21}, \eqref{D19}, \eqref{D12} and \eqref{D8} that the assumptions \eqref{rev-9} and \eqref{rev-10} of Lemma \ref{rev-7} are satisfied as $a\nearrow a^*$. We therefore conclude from Lemma \ref{rev-7} that the estimate \eqref{D16} of $(\hat{R}_a,I_a)$ holds true  as $a\nearrow a^*$.
\vskip 0.05truein

{\em Step 3.} Following Step 2, we deduce from Lemma \ref{4:lem4.2} that as $a\nearrow a^*$,
\begin{equation}\label{rev-25}
|I_a(x)|,\ |\nabla I_a(x)|\leq C\e^{4}e^{-\frac{1}{8}|x|}\ \ \hbox{in}\,\ \R^2,
\end{equation}
where $C>0$ is independent of $0<a<a^*$.
By the definition of $F_a(x)$ in \eqref{D9}, this further implies that as $a\nearrow a^*$,
\[
|F_a(x)|\leq C\e^4e^{-\frac{1}{16}|x|}\ \ \hbox{in}\,\ \R^2,
\]
where (\ref{new4}) is also used.
Similar to \eqref{D12} and \eqref{D8}, we then obtain that as $a\nearrow a^*$,
\begin{equation}\label{rev-26}
| \hat{R}_a(x)|, \ |\nabla \hat{R}_a(x)|\leq C\e^{4}e^{-\frac{1}{32}|x|}\,\ \hbox{in}\,\ \R^2.
\end{equation}
Similar to Step 2, we further conclude from \eqref{rev-25}, \eqref{rev-26}  and  Lemma \ref{rev-7} that as $a\nearrow a^*$,
\begin{equation}\label{D15}
|\hat{R}_a(x)|,\  |\nabla \hat{R}_a(x)|,\ |I_a(x)|,\  |\nabla I_a(x)|\leq C\e^{4}e^{-\frac 12 |x|}\ \ \hbox{in\ \  $\R^2$},
\end{equation}
where $C >0$ is independent of $0<a<a^*$.
Therefore, we infer from \eqref{D15} and Lemma \ref{4:lem4.2} that as $a\nearrow a^*$,
\begin{equation}\label{D27}
|I_a(x)|, \  |\nabla I_a(x)|\leq C\e^{6}e^{-\frac 18 |x|}\ \ \hbox{in\ \ $\R^2$},
\end{equation}
where $C >0$ is also independent of $0<a<a^*$. This completes the proof of Proposition \ref{prop4.3}.
\qed

\section{Nonexistence of Vortices}
Following the previous section's notations and results, in this section we shall prove Theorem \ref{thm1.2*} on the nonexistence of vortices, for which we shall employ the so-called method of inductive symmetry explained in Subsection 1.1. To prove Theorem \ref{thm1.2*}, we focus on the typical case where $V(x)=|x|^2$ and $0<\Omega<\Omega ^*$ is fixed. In this case, one can check that $\Omega ^*=2$ and $y_0=(0,0)$ is a unique and non-degenerate critical point of $H(y)$, see the assumption $(V_1)$ in Section 4, and therefore Proposition \ref{prop4.3} is applicable.

Under the above assumptions of Theorem \ref{thm1.2*}, suppose $u_a$ is a complex-valued minimizer of $e_F(a)$ as $a\nearrow a^*$, and let $w_a$ be defined by \eqref{4:1} for the rescaling form of $u_a$ with a constant phase. Denote
\begin{equation}
 w_a(x):=R_a(x)+iI_a(x)=\displaystyle\frac{w(x)}{\sqrt{a^*}}+ \hat{R}_a(x)+iI_a(x),\ \ \eps_a:=\Big(\int_{\R^2}|\nabla u_a|^{2}dx\Big)^{-\frac{1}{2}}>0,
\label{5:A:1}
\end{equation}
where $R_a(x),\ \hat{R}_a(x)$ and $  I_a (x)$ are real-valued functions in $H^1(\R^2)$. Proposition \ref{prop4.3} gives us the first estimates of $\hat R_a(x)$ and $I_a(x)$ as $a\nearrow a^*$. In this section, we first apply Proposition \ref{prop4.3} to deriving Proposition \ref{prop4.4} on the first estimate of the difference $R_a-v_a$ as $a\nearrow a^*$, instead of the difference $R_a-\frac{w}{\sqrt{a^*}}$, where $v_a>0$ defined in (\ref{5:eea}) denotes the  rescaled real-valued minimizer of $e_F(a)$ at zero rotation $\Omega =0$. Note from \cite{GLW} that $v_a>0$ is radially symmetric and unique as $a\nearrow a^*$. We then analyze the linearized operator $\mathcal{N}_a$ in Subsection 5.1, see \eqref{G10} for the definition of $\mathcal{N}_a$. To finish the proof of Theorem \ref{thm1.2*}, in Subsection 5.2 we shall apply the induction process to showing that the non-radial part of $R_a-v_a$ is arbitrarily small and
\[ \mbox{both}\ \, I_a(x)\equiv 0 \ \, \mbox{and}\ \, x_a\equiv 0\ \, \mbox{as}\ \,   a\nearrow a^*,\]
where $x_a\in\R^2$ denotes the unique global maximal point of $|u_a|$ as $a\nearrow a^*$. This is the most crucial part in proving Theorem \ref{thm1.2*}.

In order to establish Theorem \ref{thm1.2*}, we now denote $\hat {v}_a>0$ to be a  radially symmetric {\em real-valued} minimizer of the following functional
\[
e_{a}=\inf_{\{u\in \mathbb{H}, \, \|u\|^2_2=1 \} }E_{a}(u),
\]
where  the space $\mathbb{H}$ is defined by
\begin{equation}
    \mathbb{H}  :=  \Big \{u\in  H^1(\R ^2):\ \int _{\R ^2}
 |x|^2u^2 dx<\infty\Big \}, \label{H}
\end{equation}
and $E_{a}(u)$ satisfies
\[
E_{a}(u)=\inte \big(|\nabla u|^2+|x|^2u^2\big)-\frac{a}{2}\inte u^4dx.
\]
Since $\hat {v}_a>0$ is a minimizer of $e_{a}$, it satisfies the following Euler-Lagrange equation
\begin{equation}\label{G26}
-\Delta \hat {v}_a(x)+|x|^2\hat {v}_a(x)=\beta_a\hat {v}_a(x)+a\hat {v}_a^3(x) \ \ \mbox{in}\ \ \R^2,
\end{equation}
where $\beta_a\in\R$ is a Lagrange multiplier. Note that $\beta_a$ satisfies
\begin{equation}\label{Lagrange}
\beta_a=e_{a}-\frac{a}{2}\inte \hat {v}_a^4dx<0\quad\hbox{as }\ \,a\nearrow a^*,
\end{equation}
where $\hat {v}_a>0$ is a unique positive minimizer of $e_{a}$  as $a\nearrow a^*$ in view of \cite[Theorem 1.3]{GLW}.

Denote
\begin{equation}\label{5:eea}
v_a(x)=\e \hat {v}_a(\e x),\quad \e:=\Big(\inte|\nabla u_a|^{2}dx\Big)^{-\frac{1}{2}}>0.
 \end{equation}
It then follows from (\ref{G26}) that $v_a>0$ satisfies
\begin{equation}\label{G9}
-\Delta v_a(x)+\e^4|x|^2v_a(x)=\e^2\beta_a v_a(x)+av_a^3(x)\ \ \mbox{in}\, \ \R^2.
\end{equation}
Moreover, since $\e>0$ satisfies \eqref{5:eea}, we obtain from (\ref{3:DDE}) that
\begin{equation}\label{5:G10}
\e=\a+o(\a)\quad \hbox{as }\, \,a\nearrow a^*,\ \ \hbox{where }\ \, \a:=\frac{(a^*-a)^{\frac{1}{4}}}{\lambda}>0.
\end{equation}
Following Proposition \ref{prop4.3} and Lemma \ref{AP5}, we begin with the following estimates  as $a\nearrow a^*$:

\begin{prop}\label{prop4.4} Under the assumptions of Theorem \ref{thm1.2*}, let $\e>0$ and $\a>0$ be defined by (\ref{5:eea}) and (\ref{5:G10}), respectively. Then as $a\nearrow a^*$,
\begin{enumerate}
\item [(i).] $|{x_a}|\leq C\e^{5}$;
\item[(ii).]	
$|\e-\a|\leq C\a^{3}$;
\item[(iii).]
$\e^2|\beta_a-\mu _a|\leq C\a^{2}$;
\item [(iv).]
$|R_a-v_a|,\, |\nabla (R_a-v_a)|\leq C\a^2e^{-\frac 12 |x|}$ in $\R^2$,
\end{enumerate}
where $x_a\in\R^2$ denotes a global maximal point of $|u_a|$,   $R_a=R_a(x)$ is defined by \eqref{4:1}, and $\mu_a\in\R$ satisfies \eqref{3:mu}.
\end{prop}
{\noindent \bf Proof.} 1. Multiplying the first equation of \eqref{D1} by $\partial_iw:=\frac{\partial w}{\partial x_i}$ and then integrating over $\R^2$, we obtain that for $i=1,\, 2$,
\begin{equation}\label{F5}
\inte \partial_iw \Big(\mathcal{L}_a-w^2-\sqrt{a^*}R_aw\Big)\hat{R}_a=\inte (g_{1a}+g_{2a})\partial_iw,
\end{equation}
where
\[
g_{1a}(x)=-\Big[\frac{\eps_a^4\Omega^2}{4}|x|^2+\eps_a^4\Big(1-\frac{\Omega^2}{4}\Big)\big|x+\eps_a^{-1}x_a\big|^2-\e^2\mu _a-1+(a^*-a)R^2_a\Big]\frac{w}{\sqrt{a^*}},
\]
and
\[
g_{2a}(x)=\frac{aw|I_a|^2}{\sqrt{a^*}}+\e^2\Omega\big(x^{\perp}\cdot \nabla I_a\big).
\]
Applying Proposition \ref{prop4.3} gives that as $a\nearrow a^*$,
\begin{equation}\label{F8}
\Big| \inte g_{2a}(x)\partial_iw(x)dx\Big|\leq C\e^8,
\end{equation}
where $C>0$ is independent of $0<a<a^*$.

As for $\inte g_{1a}\partial_iw$,  we notice that
\begin{equation*}
\begin{split}
\inte g_{1a}\partial_iw&=-\inte \partial_iw\Big[\frac{\eps_a^4\Omega^2}{4}|x|^2+\eps_a^4\Big(1-\frac{\Omega^2}{4}\Big)
\big|x+\eps_a^{-1}x_a\big|^2-\e^2\mu _a-1+(a^*-a)R^2_a\Big]\frac{w}{\sqrt{a^*}}\\
&=-\Big(1-\frac{\Omega^2}{4}\Big)\eps_a^{4}\inte \partial_iw\Big[|x+\eps_a^{-1}x_a|^2-|x|^2\Big]\frac{w}{\sqrt{a^*}}\\
&\quad -\inte \partial_i w\Big[(a^*-a)\big(R_a+\frac{w}{\sqrt{a^*}}\big)\hat{R}_a\Big]
\frac{w}{\sqrt{a^*}},
\end{split}
\end{equation*}
due to the radial symmetry of $w$. We then obtain from above that as $a\nearrow a^*$,
\begin{equation}\label{F6}
\begin{split}
\inte g_{1a}\partial_iw&=-2\eps_a^{3}\Big(1-\frac{\Omega^2}{4}\Big)\inte \big[1+o(1)\big]\frac{\partial_i w^2}{2\sqrt{a^*}}x_a\cdot x\\
&\quad -\inte \partial_i w\Big[(a^*-a)(R_a+\frac{w}{\sqrt{a^*}})\hat{R}_a\Big]
\frac{w}{\sqrt{a^*}}\\
&:=-2\eps_a^{3}\Big(1-\frac{\Omega^2}{4}\Big)\inte \big[1+o(1)\big]\frac{\partial_i w^2}{2\sqrt{a^*}}x_a\cdot x-I.
\end{split}
\end{equation}
Note from Proposition \ref{prop4.3} that as $a\nearrow a^*$,
\[
|I|:=\Big|\inte \partial_i w\Big[(a^*-a)\big(R_a+\frac{w}{\sqrt{a^*}}\big)\hat{R}_a\Big]
\frac{w}{\sqrt{a^*}}\Big|\le C\eps_a^{8}.
\]
On the other hand, denoting the operator
\[\hat{\mathcal{L}}:=-\Delta+1-3w^2\quad\hbox{in $\R^2$,}\]
we have
\[
\begin{split}
II:&=\inte \partial_i w\Big(\mathcal{L}_a-w^2-\sqrt{a^*}R_aw\Big)\hat{R}_a\\
&=\inte \partial_i w\hat{\mathcal{L}}\hat{R}_a+\inte \partial_i w\Big(\mathcal{L}_a-w^2-\sqrt{a^*}R_aw-\hat{\mathcal{L}}\Big)\hat{R}_a\\
&=\inte \partial_i w\Big[\frac{\eps_a^4\Omega^2}{4}|x|^2+\eps_a^4\Big(1-\frac{\Omega^2}{4}\Big)|x+\eps_a^{-1}x_a|^2-\e^2\mu _a-1+2w^2\\
&\quad\qquad\quad\quad -a|w_a|^2-\sqrt{a^*}R_aw\Big]\hat{R}_a.
\end{split}
\]
Hence, we obtain from above that as $a\nearrow a^*$,
\begin{equation}\label{F7}
|II|:=\Big|\inte \partial_i w\Big(\mathcal{L}_a-w^2-\sqrt{a^*}R_aw\Big)\hat{R}_a\Big|\le C\eps_a^{8},
\end{equation}
where \eqref{lem4.1:B} and Proposition \ref{prop4.3} are used.
Therefore, we derive from \eqref{F5}--\eqref{F7} that as $a\nearrow a^*$,
\[
\begin{aligned}
&\quad\Big|2\eps_a^{3}\inte \big[1+o(1)\big]\frac{\partial_i w^2}{2\sqrt{a^*}}x_a\cdot x\Big|\leq C\e^{8}, \ \ i=1,\, 2,
\end{aligned}
\]
which further implies that Proposition \ref{prop4.4}(i) holds true.

2. On the contrary, suppose Proposition \ref{prop4.4}(ii) is false. Then there exists a constant $C(a)>0$, where $C(a)\to \infty$ as $a\nearrow a^*$, such that
\begin{equation}\label{G25}
|\e-\a|\geq C(a)\a^{3}\quad\hbox{as}\,\ a\nearrow a^*.
\end{equation}
Select the test function $w_{\tau}(x):= \frac{\tau w(\tau x)}{\|w\|_2}$, where $\tau=\lambda (a^*-a)^{-\frac{1}{4}}$ and $\lambda=\big(\inte |x|^2w^2(x)dx\big)^{\frac 14}$ is defined by \eqref{def:li}.
We then get that as $a\nearrow a^*$,
\begin{equation}\label{G23}
e_F(a) \leq F_a(w_\tau)\leq\int _{\R ^2} \big(|\nabla w_{\tau}|^2+|x|^2|w_{\tau}|^2\big)dx-\frac{a}{2}\int _{\R ^2}|w_{\tau}|^4dx=\frac{2\lambda^2(a^*-a)^{\frac{1}{2}}}{a^*}.
\end{equation}

On the other hand, we infer from Proposition \ref{prop4.3} that
\begin{equation*}
\begin{split}
\quad\Omega\inte x^\perp\cdot\Big(i w_{a}(x),\nabla  w_{a}(x)\Big)dx
&=\Omega\inte x^\perp\cdot\Big(R_{a}\nabla I_{a}-I_{a} \nabla R_{a} \Big)dx=O(\e^{10})\quad\hbox{as}\,\ a\nearrow a^*.
\end{split}
\end{equation*}
{Then from Proposition \ref{prop4.4}(i), (\ref{new7}) and Proposition \ref{prop4.3} we derive that} as $a\nearrow a^*$,
\begin{equation}\label{OK:prop:new1}
\begin{split}
e_F(a)
&\geq\displaystyle \frac{1}{\eps_{a} ^2}\inte \big|\nabla w_{a}(x)\big|^2dx-
\displaystyle \frac{a}{2\eps_{a} ^2} \inte
|w_{a}(x)|^4dx\\
&\quad+\displaystyle \eps_{a}^2\inte \Big(\frac{\Omega^2}{4}|x|^2+\frac{4-\Omega^2}{4}\Big|x+\frac{x_a}{\e}\Big|^2\Big) |w_{a}(x)|^2dx-C\e^{10}\\
&\geq \frac{a^*-a}{2\e^2}\inte |w_{a}(x)|^4dx+\displaystyle \eps_{a}^2\inte |x|^2|w_{a}(x)|^2dx-C\e^{6}\\
&\geq \frac{a^*-a}{2\e^2}\inte \Big|\frac{w}{\sqrt{a^*}}\Big|^4dx+\displaystyle \eps_{a}^2\inte|x|^2 \big|\frac{w}{\sqrt{a^*}}\big|^2dx-C\e^6:=\frac{1}{a^*}f(\e)-C\e^6,
\end{split}
\end{equation}
where $f(t):=(a^*-a)t^{-2}+t^2\inte |x|^2w^2(x)dx$ for $t>0$. Since  $f'(\a)=0$,
applying Taylor's expansion, we obtain from (\ref{OK:prop:new1}) that as $a\nearrow a^*$,
\begin{equation}
\begin{split}
e_F(a)
&\geq\frac{1}{a^*}f(\e)-C\e^6=\frac{1}{a^*}\Big[f(\a)+8\lambda^4(\e-\a)^2+O(|\e-\a|^3)\Big]-C\e^6\\
&\geq \frac{2\lambda^2(a^*-a)^{\frac{1}{2}}}{a^*}+\Big(\frac{4\lambda ^4C^2(a)}{a^*}-C\Big)\e^6,
\end{split}
\end{equation}
where both (\ref{5:G10}) and \eqref{G25} are used.
Since $C(a)\to \infty$ as $a\nearrow a^*$, the above estimate however  contradicts to \eqref{G23}. Therefore,  Proposition \ref{prop4.4}(ii) is proved.

3. Following Proposition \ref{prop4.4}(ii) and Lemma \ref{AP5}(i), one can deduce from \eqref{lem4.1:B} that as $a\nearrow a^*$,
\[
\begin{aligned}
 \quad|\e^2(\beta_a-\mu _a)|\leq |1+\e^2\beta_a|+|1+\e^2\mu _a| \leq |1+\alpha_a^2\beta_a|+|(\alpha_a^2-\e^2)\beta_a|+C\a^4\leq C\a^2,
\end{aligned}
\]
which then implies that Proposition \ref{prop4.4}(iii) holds true.

4. By Lemma \ref{AP5}, we deduce from Proposition \ref{prop4.3} that as $a\nearrow a^*$,
\[
\begin{aligned}
&\quad\Big|R_a(x)-v_a(x)\Big|\\
&\leq\Big|R_a(x)-\frac{w(x)}{\sqrt{a^*}}\Big|+\Big|\frac{w(x)}{\sqrt{a^*}}-\a\hat{v}_a(\a x)\Big|+\big|\a\hat{v}_a(\a x)-v_a(x)\big|\\
&\leq C\e^4e^{-\frac{1}{2}|x|}+C\a^{4}e^{-\frac{2}{3}|x|}+C\Big|\frac{(\a-\e)}{\a}\a\hat{v}_a(\a x)\Big|\\
&\quad +C\Big|\frac{(\a-\e)x}{\a}\Big|\big|\nabla\big(\a\hat {v}_a(\a x+\theta(\alpha_a-\e)x)\big)\big|\\
&\leq C\e^2e^{-\frac{1}{2}|x|}\quad\mbox{in $\R^2$},
\end{aligned}
\]
{where $\theta=(\theta_1,\theta_2)$ is a vector with $|\theta|<1.$} Similarly, one can obtain the similar estimate of $\nabla (R_a-v_a)$ as $a\nearrow a^*$. This completes the proof of Proposition \ref{prop4.4}(iv), and we are done.   \qed

\subsection{Analysis of the linearized operator $\mathcal{N}_a$}

In view of (\ref{G9}),  this subsection is devoted to the analysis of the linearized operator $\mathcal{N}_a:\, D(\mathcal{N}_a)\subset L^2(\R^2) \longmapsto L^2(\R^2)$ defined as
\begin{equation}\label{G10}
\mathcal{N}_a:=-\Delta+\big(\e^4|x|^2-\e^2\beta_a-3av_a^2\big),
\end{equation}
where $v_a>0$ and $\e >0$ are defined in (\ref{5:eea}), and the domain $D(\mathcal{N}_a)$ satisfies
\begin{equation}\label{5A:G10}
D(\mathcal{N}_a)=H^2(\R^2)\cap \Big\{u\in L^2(\R^2):\inte |x|^4u^2dx< \infty \Big\}.
\end{equation}
It then follows from \cite[Corollary 1.5]{BO} that $\mathcal{N}_a$ is non-degenerate and $ker(\mathcal{N}_a)=\{0\}$. The main aim of this subsection is to establish the following estimate:

\begin{lem}\label{lem4.6}
Under the assumptions of Theorem \ref{thm1.2*},
for any $\alpha\geq 0$ suppose $T_{\alpha,a}\in D(\mathcal{N}_a)$ is the unique solution of the following problem
\begin{equation}\label{G31}
\mathcal{N}_aT_{\alpha,a}=f_{\alpha,a}(x)\in L^2(\R^2),
\end{equation}
\begin{equation}\label{G39}
|\nabla T_{\alpha,a}(0)|\leq C(\alpha)\e^{2\alpha},\ \,
|f_{\alpha,a}(x)|\leq C(\alpha)\e^{2\alpha}e^{-\frac {1}{16}|x|}\ \ \hbox{in} \, \ \R^2,
\end{equation}
where $C(\alpha)>0$ depends only on $\alpha \ge0$.
Then there exists a positive constant $C>0$, independent of $\alpha$ and $a$, such that as $a\nearrow a^*$,
\[
|T_{\alpha,a}(x)|,\ |\nabla T_{\alpha,a}(x)|\leq CC(\alpha)\e^{2\alpha}e^{-\frac {1}{32}|x|}\ \ \hbox{in} \, \ \R^2.
\]
\end{lem}

To prove Lemma \ref{lem4.6}, motivated by \cite[Lemma 4.2]{NT2} we first employ Proposition \ref{prop4.4}  to discuss the following analytical properties of $\mathcal{N}_a$.

\begin{thm}\label{AP7}
Under the assumptions of Theorem \ref{thm1.2*}, we have
\begin{enumerate}
\item [(i).] $\mathcal{N}^{-1}_a$: $L^2(\R^2)\longmapsto L^2(\R^2)$ exists and is a continuous linear operator. Moreover, the following estimate holds true:
\begin{equation}\label{AP10}
\|\psi(x)\|_{H^2(\R^2)}\leq C\|\mathcal{N}_a\psi(x)\|_{L^2(\R^2)},\quad\hbox{if $\psi\in \big\{\frac{\partial w}{\partial x_1}, \frac{\partial w}{\partial x_2}\big\}^{\perp}$,}
\end{equation}
where $C>0$ is a constant independent of $0<a<a^*$.
\item [(ii).]
If $\phi(x)\in L^2(\R^2)$ is radially symmetric, then $\psi(x)=\mathcal{N}_a^{-1}\phi(x)\in D(\mathcal{N}_a)$ is also radially symmetric, where the space $D(\mathcal{N}_a)$ is defined by (\ref{5A:G10}).
\end{enumerate}
\end{thm}

\noindent{\bf Proof.} (i). For any given $f(x)\in L^2(\R^2)$, consider the following equation
\begin{equation}\label{AP8}
-\Delta u+\big(\e^4|x|^2-\e^2\beta_a-3av_a^2\big)u=f \ \ \hbox{in}\ \, \R^2.
\end{equation}
By \cite[Theorem 4.1]{OK} or \cite[Theorem XIII.67]{RS}, the operator $\big(-\Delta +(\e^4|x|^2-\e^2\beta_a )\big)^{-1}$ is compact from $L^2(\R^2)$ onto $L^2(\R^2)$. Hence, the equation \eqref{AP8} is solvable, if and only if  the following holds
\begin{equation}\label{AP9}
u-3a\Big(-\Delta +\big(\e^4|x|^2-\e^2\beta_a\big)\Big)^{-1}(v_a^2u)=\Big(-\Delta +\big(\e^4|x|^2-\e^2\beta_a\big)\Big)^{-1}f.
\end{equation}
One can obtain from Proposition \ref{prop4.4}(iv) that $v_a\in L^{\infty}(\R^2)$, which further implies that $\big(-\Delta +(\e^4|x|^2-\e^2\beta_a)\big)^{-1}\circ v_a^2$ is a compact operator. By Riesz-Fredholm theory of compact linear operators, we thus deduce that the equation \eqref{AP9} is solvable, if and only if $ker \mathcal{N}_a=\{0\}$, but the latter is true in view of \cite[Corally 1.5]{BO}. This gives that $\mathcal{N}_a^{-1}$ exists. Therefore, $\mathcal{N}_a$ is a one-to-one and  continuous linear  mapping from $D(\mathcal{N}_a)$ onto $ L^2(\R^2)$, which implies that $\mathcal{N}_a^{-1}$ is also continuous by the open mapping theorem.

We next prove the estimate \eqref{AP10} as follows. By the above argument, for any given $\phi(x)\in L^2(\R^2)$, assume $\psi\in D(\mathcal{N}_a)$ is the unique solution of the following equation
\begin{equation}\label{AP11}
\mathcal{N}_a\psi(x)=\phi(x)\ \ \hbox{in}\ \, \R^2.
\end{equation}
For simplicity, we denote
\[
d_a(x)=-\e^2\beta_a-1+3w^2-3av^2_a, \ \ \hat{\mathcal{L}}=-\Delta+1-3w^2\ \ \hbox{in}\ \, \R^2,
\]
so that
\[
\mathcal{N}_a=\hat{\mathcal{L}}+\e^4|x|^2+d_a(x).
\]
By (1) of Lemma \ref{lem3.3} and (iii)-(iv) of Proposition \ref{prop4.4}, we obtain that
\[\|d_a\|_{C^1(\R^2)}\to 0\ \ \hbox{as}\ \, a\nearrow a^*.\]
Thus, we infer from \eqref{AP11} that
\begin{equation}\label{AP12}
\begin{split}
\inte |\phi(x)|^2dx
&=\inte \big(\hat{\mathcal{L}}\psi\big)^2dx+\inte \big(\e^4|x|^2\psi+d_a\psi\big)^2dx+2\inte \hat{\mathcal{L}}\psi\big(\e^4|x|^2\psi+d_a\psi\big)dx\\
&=I_1+I_2+I_3.
\end{split}
\end{equation}
Applying \cite[Lemma 4.2]{NT2}, there exists a constant $\delta >0$ such that $I_1$ satisfies
\[
\begin{aligned}
I_1:=\inte \big(\hat{\mathcal{L}}\psi\big)^2dx\geq \delta\|\psi\|^2_{H^2(\R^2)}, \ \ \hbox{if\ $\psi\in\big\{\frac{\partial w}{\partial x_1}, \frac{\partial w}{\partial x_2}\big\}^{\perp}$}.
\end{aligned}
\]
We also derive that $I_2:=\inte \big(\e^4|x|^2\psi+d_a\psi\big)^2dx\geq 0$, and however $I_3$ satisfies
\[
\begin{aligned}
I_3:&=\inte 2\hat{\mathcal{L}}\psi\big(\e^4|x|^2\psi+d_a\psi\big)dx\\
&=\inte \Big\{2\e^4\big(|x|^2|\nabla \psi|^2+2\psi x\cdot\nabla\psi\big)+2d_a|\nabla\psi|^2+2\psi\nabla\psi \nabla d_a\\
&\quad+2(1-3w^2)\big(\e^4|x|^2+d_a\big)\psi^2\Big\}dx\\
&\geq \inte \big[-2\e^4|\psi|^2-2|d_a||\nabla\psi|^2-|\nabla d_a|\big(|\psi|^2+|\nabla\psi|^2\big)\\
&\quad-2|d_a||1-3w^2|\psi^2-6\e^4w^2|x|^2\psi^2\big]dx\\
&\geq -\inte \big[C\max\big\{\e^4,\,|\nabla d_a|,\,|d_a|\big\}\psi^2+\big(2|d_a|+|\nabla d_a|\big)|\nabla\psi|^2\big]dx \\
&\geq -\frac{\delta}{2}\|\psi\|^2_{H^2(\R^2)} \ \ \, \hbox{as\quad $a\nearrow a^*$}.
\end{aligned}
\]
By applying above estimates, we conclude from (\ref{AP12}) that if  $\psi\in\big\{\frac{\partial w}{\partial x_1}, \frac{\partial w}{\partial x_2}\big\}^{\perp}$, then as $a\nearrow a^*$,
\[
\inte |\phi(x)|^2dx=\inte \big(\hat{\mathcal{L}}\psi+\e^4|x|^2\psi+d_a\psi\big)^2dx\geq \frac{\delta}{2}\|\psi\|^2_{H^2(\R^2)},
\]
which thus implies that \eqref{AP10} holds.

(ii). Set $\bar {\psi}(x)=\frac{1}{2\pi|x|}\int_{|y|=|x|}{\psi}(y)dy $ for $x\neq (0,0)$, and $\bar {\psi}(x)={\psi}(0)$ for $x=(0,0)$, where $\psi$ satisfies (\ref{AP11}).
We then   obtain from (\ref{AP11}) that $\bar {\psi} (x)\in D(\mathcal{N}_a)$ satisfies
$
\mathcal{N}_a\bar {\psi}(x)=\overline{\mathcal{N}_a {\psi}(x)}=\bar{\phi}(x)=\phi(x)=\mathcal{N}_a{\psi}(x)\ \ \hbox{in}\ \, \R^2.
$
By the non-degeneracy of $\mathcal{N}_a$ (cf. \cite[Corollary 1.5]{BO}), it further implies that ${\psi}(x)$ is also radially symmetric. This completes the proof of Theorem \ref{AP7}. \qed

 \medskip

{\noindent \bf Proof of Lemma \ref{lem4.6}.} Choose suitable constants $a_{\alpha,a}$ and $b_{\alpha,a}$ such that
\begin{equation}\label{G32}
T_{\alpha,a}=a_{\alpha,a}\partial_1w+b_{\alpha,a}\partial_2w+\hat{T}_{\alpha,a},\,\ \inte\hat{T}_{\alpha,a}\partial_i w=0\,\ \hbox{for}\, \ i=1, 2,
\end{equation}
where and below we denote $\partial_i w=\frac{\partial w}{\partial x_i}$ for $i=1$ and $2$. Without loss of generality, we may assume $|a_{\alpha,a}|\geq |b_{\alpha,a}|$.

We first claim that as $a\nearrow a^*$,
\begin{equation}\label{G36}
|b_{\alpha,a}|\leq |a_{\alpha,a}|\leq C C(\alpha)\e^{2\alpha}.
\end{equation}
Indeed, following \eqref{G32}, we have
\[
\begin{aligned}
\mathcal{N}_a\hat{T}_{\alpha,a}&=f_{\alpha,a}-\mathcal{N}_a(a_{\alpha,a}\partial_1w+b_{\alpha,a}\partial_2w)\\
&=f_{\alpha,a}+\big(-\e^4|x|^2+\e^2\beta_a+1-3w^2+3av_a^2\big)\big(a_{\alpha,a}\partial_1w+b_{\alpha,a}\partial_2w\big)\ \ \hbox{in} \, \ \R^2.
\end{aligned}
\]
By \eqref{AP10} and \eqref{G39}, we infer from the above equation that as $a\nearrow a^*$,
\begin{equation}\label{G33}
\begin{split}
&\quad\big\|\hat{T}_{\alpha,a}\big\|_{L^{\infty}(\R^2)}\leq C\big\|\hat{T}_{\alpha,a}\big\|_{H^2(\R^2)}\\
 &\leq C\big\|f_{\alpha,a}+(-\e^4|x|^2+\e^2\beta_a+1-3w^2+3av_a^2)(a_{\alpha,a}\partial_1w+b_{\alpha,a}\partial_2w)\big\|_{L^2(\R^2)}\\
&\leq C\|f_{\alpha,a}\|_{L^2(\R^2)}+C\e^2\big(|a_{\alpha,a}|+|b_{\alpha,a}|\big)\\
&\leq CC(\alpha)\e^{2\alpha}+C\e^2|a_{\alpha,a}|,
\end{split}
\end{equation}
where Proposition \ref{prop4.4} and Lemma \ref{AP5} are also used in the second inequality. By the gradient estimates, we then obtain from (\ref{G39}) and \eqref{G33} that as $a\nearrow a^*$,
\begin{equation}\label{G34}
\begin{split}
&\quad|\nabla\hat{T}_{\alpha,a}(0)|\\
&\leq C\big\|f_{\alpha,a}+(-\e^4|x|^2+\e^2\beta_a+1-3w^2+3av_a^2)(a_{\alpha,a}\partial_1w+
b_{\alpha,a}\partial_2w)\big\|_{L^{\infty}(B_1)}\\
&\quad+C\big\|\hat{T}_{\alpha,a}\big\|_{L^{\infty}(\partial B_1)}\\
&\leq CC(\alpha)\e^{2\alpha}+C\e^2|a_{\alpha,a}|.
\end{split}
\end{equation}
On the other hand, we have for $\lambda_{\alpha,a}:=\frac{b_{\alpha,a}}{a_{\alpha,a}}$,
\begin{equation}\label{G35}
\begin{split}
&\quad\big|\nabla (a_{\alpha,a}\partial_1w+b_{\alpha,a}\partial_2w)(0)\big|\\
&=|a_{\alpha,a}|\sqrt{\big[\partial_1\partial_1w(0)+\lambda_{\alpha,a}\partial_1\partial_2w(0)\big]^2+\big[\partial_1\partial_2w(0)+\lambda_{\alpha,a}\partial_2\partial_2w(0)\big]^2}\\
&\geq |a_{\alpha,a}|\frac{\big|\partial_1\partial_1w(0)\partial_2\partial_2w(0)-
[\partial_1\partial_2w(0)]^2\big|}{\sqrt{[\partial_1\partial_2w(0)]^2+[\partial_2\partial_2w(0)]^2}}:=\eta |a_{\alpha,a}|,
\end{split}
\end{equation}
where the constant $\eta >0$ is independent of $\alpha\ge 0$ and $0<a<a^*$.
Since $w(0)$ is non-degenerate, i.e.,
\[
\partial_1\partial_1w(0)\partial_2\partial_2w(0)-[\partial_1\partial_2w(0)]^2\neq 0,
\]
we have $\eta>0$. Therefore, we deduce from \eqref{G34} and \eqref{G35} that as $a\nearrow a^*$,
\[
\begin{aligned}
|\nabla{T}_{\alpha,a}(0)|&=\big|\nabla\hat{T}_{\alpha,a}(0)+\nabla (a_{\alpha,a}\partial_1w+b_{\alpha,a}\partial_2w)(0)\big|\\
&\geq \big|\nabla (a_{\alpha,a}\partial_1w+b_{\alpha,a}\partial_2w)(0)\big|-\big|\nabla\hat{T}_{\alpha,a}(0)\big|\\
&\geq \frac {\eta}{2}|a_{\alpha,a}|-CC(\alpha)\e^{2\alpha}.
\end{aligned}\]
Under the assumption \eqref{G39}, we therefore infer from the above inequality that \eqref{G36} holds true.

Under the assumption \eqref{G39}, we now obtain from \eqref{G32}--\eqref{G33} that
\begin{equation}\label{G38}
\|T_{\alpha,a}\|_{L^{\infty}(\R^2)}\leq C|a_{\alpha,a}|+
\|\hat T_{\alpha,a}\|_{L^{\infty}(\R^2)}\leq CC(\alpha)\e^{2\alpha}.
\end{equation}
Note from \eqref{G31} and \eqref{G39} that $T_{\alpha,a}$ satisfies
\[
|\mathcal{N}_aT_{\alpha,a}|\leq CC(\alpha)\e^{2\alpha}e^{-\frac{1}{16}|x|}\ \ \hbox{in} \ \, \R^2/B_R,
\]
where $R>0$ is large enough. Since $-\e^2\beta_a\to 1$ as $a\nearrow a^*$, by the comparison principle we deduce from \eqref{G38} and Proposition \ref{prop4.4}(iv) that as $a\nearrow a^*$,
\[
|T_{\alpha,a}(x)|\leq CC(\alpha)\e^{2\alpha}e^{-\frac{1}{16}|x|}\ \ \hbox{in} \ \, \R^2/B_R.
\]
By gradient estimates, together with \eqref{G38}, we further conclude from above that as $a\nearrow a^*$,
\[
|\nabla T_{\alpha,a}(x)|\leq CC(\alpha)\e^{2\alpha}e^{-\frac{1}{32}|x|}\ \ \hbox{in} \ \, \R^2,
\]
which completes the proof of Lemma \ref{lem4.6}. \qed

\subsection{Proof of Theorem \ref{thm1.2*}}

The main purpose of this subsection is to complete the proof of Theorem \ref{thm1.2*}. Before that, we need the following induction argument.

\begin{lem}\label{lem4.5}
Under the assumptions of Theorem \ref{thm1.2*}, for any $m\geq 2$ suppose ${R}_a(x)-v_a(x)$  can be decomposed as
\begin{equation}\label{D35}
Q_a(x):={R}_a(x)-v_a(x)=\psi_{m,a}(|x|)+T_{m,a}(x),
\end{equation}
\begin{equation}\label{D33}
|T_{m,a}(x)|,\ \, |\nabla T_{m,a}(x)|\leq A^m\e^{2m}e^{-\frac{1}{2}|x|}\ \ \hbox{in\ \,$\R^2$}\, \ \mbox{as}\ \, a\nearrow a^*,
\end{equation}
where $\psi_{m,a}(|x|)$ is radially symmetric and the constant $A>0$ is independent of $m $ and  $0<a<a^*$. Then the decomposition of (\ref{D35}) and (\ref{D33}) holds  for $m+1$.
\end{lem}

{\noindent \bf Proof.} Recall first that ${R}_a(x)$ and $v_a(x)$ are defined in (\ref{4:1}) and (\ref{5:eea}), respectively.  Under the assumptions of Theorem \ref{thm1.2*}, suppose $R_a=v_a+\psi_{m,a}+T_{m,a}$ holds for some $m\ge 2$, where $v_a$ and $\psi_{m,a}(|x|)$ are radially symmetric, and $T_{m,a}$ satisfies (\ref{D33}). For convenience, we always use the symbol $C$ to denote some positive constant independent of $A>0$, $m\ge 2$ and  $0<a<a^*$. We also consider sufficiently small $\e >0$, so that $0<A\e<1$ as $a>0$ is sufficiently close to $a^*$ from below.

For any $m\geq 2$, by Proposition \ref{prop4.4}(iv), we obtain from  \eqref{D35} and \eqref{D33} that
\begin{equation}\label{D40}
\begin{split}
|\psi_{m,a}(x)|,\,|\nabla \psi_{m,a}(x)|&\leq CA^m\e^{2m}e^{-\frac{1}{2}|x|}+C\e^{2}e^{-\frac{1}{2}|x|}\\
&\leq C\e^{m}e^{-\frac{1}{2}|x|}+C\e^{2}e^{-\frac{1}{2}|x|}\\
&\leq C\e^{2}e^{-\frac{1}{2}|x|} \quad\hbox{in}\,\ \R^2.
\end{split}
\end{equation}
Since $v_a(|x|)+\psi_{m,a}(|x|)$ is radially symmetric, we derive from \eqref{D35} and \eqref{D33} that for any $m\geq 2$,
\[
\big|x^{\perp}\cdot \nabla R_{a}\big|=\big|x^{\perp}\cdot \nabla T_{m,a}\big|\leq CA^{m}\e^{2m}e^{-\frac{1}{4}|x|}\ \ \mbox{in\,\ $\R^2$}.
\]
We then obtain from Lemma \ref{4:lem4.2} that
\begin{equation}\label{G14}
|I_a(x)|,\ |\nabla I_a(x)|\leq CA^{m}\e^{2(m+1)}e^{-\frac{1}{8}|x|}\ \ \mbox{in\,\ $\R^2$}.
\end{equation}
Since $R_a=v_a+\psi_{m,a}+T_{m,a}$, we deduce from (\ref{4:3}) and \eqref{G10} that
\begin{equation}\label{G1}
\begin{split}
\mathcal{N}_aQ_a:&=\frac{(4-\Omega^2)\e^{4}}{4}\big(|x|^2-|x+{\e^{-1}}x_a|^2\big)R_a-\e^2(\beta_a-\mu _a)R_a\\
&\quad+a(R_a+2v_a)Q^2_a+\e^2\Omega\big(x^{\perp}\cdot\nabla I_a\big)+aI_a^2R_a\\
&=\e^2\Omega\big(x^{\perp}\cdot\nabla I_a\big)+aI_a^2R_a+B_{1,a}(x)+B_{2,a}(x)\quad\hbox{in}\,\ \R^2,
\end{split}
\end{equation}
where the radially symmetric function $B_{1,a}(x)$ satisfies
\begin{equation}\label{G4}
\begin{split}
B_{1,a}(x)&=-\e^2(\beta_a-\mu _a)(v_a+\psi_{m,a})+a(\psi_{m,a}^3+3v_a\psi_{m,a}^2),
\end{split}
\end{equation}
and the non-radially symmetric function $B_{2,a}(x)$ satisfies
\begin{equation}\label{G5}
\begin{split}
B_{2,a}(x)&=\frac{(4-\Omega^2)\e^{4}}{4}\big(|x|^2-|x+{\e^{-1}}x_a|^2\big)R_a-\e^2(\beta_a-\mu _a)T_{m,a}\\
&\quad+a(3\psi_{m,a}^2+3\psi_{m,a}T_{m,a}+T^2_{m,a}+6v_a\psi_{m,a}+3v_aT_{m,a})T_{m,a}.
\end{split}
\end{equation}
Following Theorem \ref{AP7} and (\ref{G4}), there exists a unique function $\psi_{m+1,a}\in C^2(\R^2)\cap L^\infty(\R^2)$ of
\begin{equation}\label{G6}
\begin{split}
\mathcal{N}_a\psi_{m+1,a}=B_{1,a}(x)\quad\hbox{in\  $\R^2$,}
\end{split}
\end{equation}
and moreover, $\psi_{m+1,a}$ is also radially symmetric. We now define
\begin{equation}\label{G7}
T_{m+1,a}(x):=Q_a(x)-\psi_{m+1,a}(x),
\end{equation}
so that $T_{m+1,a}(x)$ satisfies
\begin{equation}\label{G16}
\mathcal{N}_aT_{m+1,a}=\e^2\Omega\big(x^{\perp}\cdot\nabla I_a\big)+aI_a^2R_a+B_{2,a}
\ \ \mbox{in}\,\  \R^2,
\end{equation}
due to (\ref{G1}).
Therefore, under the decomposition of (\ref{G7}), the rest is to prove that (\ref{D33}) holds  for $m+1$. We shall prove it via the following three steps.

\vskip 0.05truein

{\em Step 1.} Since $T_{m,a}$ satisfies (\ref{D33}), we claim that
\begin{equation}\label{G11}
|{x_a}|\leq CA^{m}\e^{2m-1}\quad \mbox{as}\ \ a\nearrow a^*,
\end{equation}
where $x_a$ is the unique global maximal point of $|u_a|$ as $a\nearrow a^*$.

Indeed, since $B_{1,a}(x)$ defined in (\ref{G4}) is radially symmetric, we deduce from \eqref{G1} that
\begin{equation}\label{G12}
\begin{split}
&\inte \partial_iw\mathcal{N}_aQ_a=\Omega\e^2\inte\partial_iw\big( x^{\perp}\cdot\nabla I_a\big)+a\inte \partial_i wI^2_aR_a+\inte \partial_i wB_{2,a}.
\end{split}
\end{equation}
Following \eqref{G5} and \eqref{G12} yields that
\begin{equation}\label{G27}
\begin{split}
&\quad \inte \frac{(4-\Omega^2)\e^{4}}{4}\partial_i w\big(|x|^2-|x+{\e^{-1}}x_a|^2\big)R_a\\
&=\inte \partial_iw\mathcal{N}_aQ_a-\Omega\e^2\inte\partial_iw\big( x^{\perp}\cdot\nabla I_a\big)-a\inte \partial_i wI^2_aR_a\\
&\quad-\inte \partial_i w\Big[\e^2(\mu_a-\beta_a)T_{m,a}+a\big(3\psi_{m,a}^2+3\psi_{m,a}T_{m,a}\\
&\quad\qquad \qquad \quad+T_{m,a}^2+6v_a\psi_{m,a}+3v_aT_{m,a}\big)T_{m,a}\Big].
\end{split}
\end{equation}
We obtain from (\ref{G7}) that for $\hat{\mathcal{L}}=-\Delta +1-3w^2$,
\begin{equation}\label{G13}
\begin{split}
\inte \partial_iw\mathcal{N}_aQ_a&=\inte \partial_iw\hat{\mathcal{L}}Q_a+\inte \partial_iw(\mathcal{N}_a-\hat{\mathcal{L}})Q_a\\
&=\inte \partial_iw\Big[{\e^4|x|^2}-\e^2\beta_a-1+3w^2-3av_a^2\Big]Q_a\\
&=\inte \partial_iw\Big[{\e^4|x|^2}-\e^2\beta_a-1+3w^2-3av_a^2\Big]T_{m,a}.
\end{split}
\end{equation}
By Proposition \ref{prop4.4} and Lemma \ref{AP5}, we get from \eqref{D33} and (\ref{G13}) that  as $a\nearrow a^*$,
\[
\Big|\inte \partial_iw\mathcal{N}_aQ_a\Big|=\Big|\inte \partial_iw\big({\e^4|x|^2}-\e^2\beta_a-1+3w^2-3av_a^2\big)T_{m,a}\Big|\leq CA^m\e^{2(m+1)}.
\]
Note  from  Proposition \ref{prop4.4} and \eqref{D33} that  as $a\nearrow a^*$,
\[
\Big|\e^2(\mu_a-\beta_a)\inte \partial_i wT_{m,a}\Big|\leq CA^m\e^{2(m+1)}.
\]
Because $\e>0$ is small enough such that $0<A\e<1$, we also have as $a\nearrow a^*$,
\[
a\Big|\inte \partial_i wI^2_aR_a\Big|\leq CA^m\e^{2(m+1)},
\]
and
\[
a\Big|\inte \partial_i w\big(3\psi_{m,a}^2+3\psi_{m,a}T_{m,a}+T_{m,a}^2+6v_a\psi_{m,a}+3v_aT_{m,a}\big)T_{m,a}\Big|\leq CA^m\e^{2(m+1)},
\]
due to (\ref{D40}) and (\ref{G14}).
Using above estimates, we then conclude from (\ref{G14}) and (\ref{G27}) that  as $a\nearrow a^*$,
\begin{equation}\label{G28}
\Big|\inte \frac{(4-\Omega^2)\e^{4}}{4}\partial_i w\big(|x|^2-|x+{\e^{-1}}x_a|^2\big)R_a\Big|\leq CA^m\e^{2(m+1)}.
\end{equation}
Since
\[
\begin{split}
&\inte  \frac{(4-\Omega^2)\e^{4}}{4}\partial_i w\big(|x|^2-|x+{\e^{-1}}x_a|^2\big)R_a\\
&=-\inte \frac{(4-\Omega^2)\e^{3}}{2}\partial_i w(x_a\cdot x)R_a
-\inte  \frac{(4-\Omega^2)\e^{4}}{4}\partial_i w\e^{-2}x_a^2T_{m,a}\\
&=\inte \frac{(4-\Omega^2)\e^{3}}{4} w^2x^{(i)}_a+O(A^m\varepsilon_a^{2(m+1)}),\ \ i=1,\, 2,
\end{split}\]
where $x_a=(x^{(1)}_a, x^{(2)}_a)$, we further deduce from \eqref{G28} that the claim \eqref{G11} holds true.

\vskip 0.05truein

{\em Step 2.} We claim that as $a\nearrow a^*$,
\begin{equation}\label{G17*}
|T_{m+1,a}|,\ |\nabla T_{m+1,a}|\leq CA^{m}\e^{2(m+1)}e^{-\frac{1}{32}|x|}\ \ \hbox{in\,\ $\R^2$.}
\end{equation}

Actually, following \eqref{G11} we obtain from (\ref{D33}), \eqref{new4} and \eqref{G5}  that   as $a\nearrow a^*$,
\begin{equation}\label{5:G14}
|B_{2,a}(x)|,\ |\nabla B_{2,a}(x)|\leq CA^m\e^{2(m+1)}e^{-\frac{1}{2}|x|}\ \ \hbox{in}\, \ \R^2.
\end{equation}
Since
\[T_{m+1,a}(x)=Q_a(x)-\psi_{m+1,a}(x)=R_a(x)-[v_a(x)+\psi_{m+1,a}(x)],
\]
where $v_a(x)+\psi_{m+1,a}(x)$ is radially symmetric, we derive from \eqref{G14} that as $a\nearrow a^*$,
\begin{equation}\label{G37}
\begin{split}
|\nabla T_{m+1,a}(0)|=|\nabla R_a(0)|=\Big|-\frac{I_a(0)\nabla I_a(0)}{{R}_a(0)}\Big|\leq CA^{m}\e^{2(m+1)}.
\end{split}
\end{equation}
In view of (\ref{G14}), (\ref{5:G14}) and \eqref{G37}, applying Lemma \ref{lem4.6} we therefore conclude from (\ref{G16}) that the claim (\ref{G17*}) holds true.


{\em Step 3.} In this step, we prove that (\ref{D33}) holds for $m+1$.
Recall from (\ref{4:3}), (\ref{G7}) and (\ref{G16}) that $(T_{m+1,a},I_a)$ satisfies the following system
\begin{equation}\label{G8}
\left\{
\begin{aligned}
\mathcal{N}_aT_{m+1,a}&=\e^2\Omega\big(x^{\perp}\cdot\nabla I_a\big)+aI_a^2R_a+B_{2,a}(x)
\ \ \mbox{in}\,\ \R^2,\\
\mathcal{L}_a I_a&=-\eps^2_a\Omega (x^\bot\cdot\nabla T_{m+1,a})
\ \ \mbox{in}\,\ \R^2.\\
\end{aligned}
\right.
\end{equation}
Hence, we obtain from \eqref{G8} that $(T_{m+1,a},I_a)$ satisfies the system
\eqref{rev-2} with
\[
b_1=-b_2=\e^2\Omega,\,P=CA^m\e^{2(m+1)},\,Q=2\delta=\frac 12,
\]
\[
f_1(x)=3av_a^2T_{m+1,a}+aI_a^2R_a+B_{2,a}(x),\,f_2(x)=a|w_a|^2I_a,
\]
\[
V_1(x)=\e^4|x|^2-\e^2\beta_a,\,V_2(x)=\frac{\e^4\Omega^2|x|^2}{4}+\e^4\Big(1-\frac{\Omega^2}{4}\Big)|x+x_a\e^{-1}|^2-\e^2\mu_a.
\]
By Steps 1 and 2, we obtain from \eqref{G14} that $(T_{m+1,a},I_a)$ satisfies \eqref{rev-8}, and $V_1(x)$ and $V_2(x)$ satisfy the assumption \eqref{rev-10} of  Lemma \ref{rev-7}. Moreover, we infer from \eqref{G14}, \eqref{5:G14} and \eqref{G17*} that $f_1(x)$ and $f_2(x)$ satisfy the assumption \eqref{rev-9} of  Lemma \ref{rev-7}. Therefore, we deduce from Lemma \ref{rev-7} that
\begin{equation}\label{rev-18}
|I_a|,\ |\nabla I_a|,\,|T_{m+1,a}|,\ |\nabla T_{m+1,a}|\leq CA^{m}\e^{2(m+1)}e^{-\frac{1}{2}|x|}\ \ \hbox{in\,\ $\R^2$.}
\end{equation}
Because the positive constant $C$ in (\ref{rev-18}) is independent of $A>0$, $m\ge 2$ and  $0<a<a^*$, one can choose a sufficiently large constant $A$ such that $A>C$. Therefore, we conclude from (\ref{rev-18}) that (\ref{D33}) holds for $m+1$, and the proof is complete. \qed
\medskip


\vskip 0.05truein
{\noindent \bf Proof of Theorem \ref{thm1.2*}.} Set $T_{1,a}(x):=R_a(x)-v_a(x)$, where ${R}_a(x)$ and $v_a(x)$ are defined in (\ref{4:1}) and (\ref{5:eea}), respectively. We then obtain from Proposition \ref{prop4.4}(iv) that
\begin{equation}\label{5A:G29}
|T_{1,a}(x)|,\ |\nabla T_{1,a}(x)|\leq C_1\e^2e^{-\frac 12 |x|}\ \ \mbox{in}\ \, \R^2,
\end{equation}
where $C_1>0$ is independent of  $0<a<a^*$.
Stimulated by (\ref{G1}), let $\psi_{2,a}\in C^2(\R^2)\cap L^\infty(\R^2)$  be the unique solution of the following equation
\[
\mathcal{N}
_a\psi_{2,a}=-\e^2(\beta_a-\mu _a)v_a\ \ \mbox{in}\ \, \R^2,
\]
and set $T_{2,a}:=R_a-v_a-\psi_{2,a}$. Theorem \ref{AP7}(ii) then gives that $\psi_{2,a}(|x|)$ is radially symmetric. Moreover, based on (\ref{5A:G29}), the same argument of proving Lemma \ref{lem4.5} gives that there exists a constant $C_2>0$, independent of $0<a<a^*$, so that
\[
|T_{2,a}(x)|,\ |\nabla T_{2,a}(x)|\leq C_2\e^4e^{-\frac 12 |x|}\ \ \mbox{in}\ \, \R^2\quad \mbox{as}\ \ a\nearrow a^*.
\]
Take $A>0$ large enough that $A^2>C_2$, from which we then have
\begin{equation}\label{5B:G29}
|T_{2,a}(x)|,\ |\nabla T_{2,a}(x)|\leq A^2\e^4e^{-\frac 12 |x|}\ \ \mbox{in}\ \, \R^2\quad \mbox{as}\ \ a\nearrow a^*.
\end{equation}
By Lemma \ref{lem4.5}, we thus deduce from (\ref{5B:G29}) that for any $m\geq 2$,
\[
| T_{m,a}(x)|,\ |\nabla T_{m,a}(x)|\leq A^{m}\e^{2m}e^{-\frac 12|x|}\ \ \mbox{in}\ \, \R^2\quad \mbox{as}\ \ a\nearrow a^*.
\]
Recall from (\ref{D35}) that ${R}_a(x)=[v_a(x)+\psi_{m,a}(|x|)]+T_{m,a}(x)$, where $v_a(x)+\psi_{m,a}(|x|)$ is radially symmetric.
Applying Lemma \ref{4:lem4.2}, we then derive from above that for any $m\geq 2$,
\begin{equation}\label{5C:G29}
|I_a(x)|,\ |\nabla I_a(x)|\leq CA^m\e^{2m+2}e^{-\frac 18|x|}\ \ \mbox{in}\ \, \R^2\quad \mbox{as}\ \ a\nearrow a^*.
\end{equation}
Moreover, the proof of  Lemma \ref{lem4.5}, see (\ref{G11}), gives that  for any $m\geq 2$,
\begin{equation}\label{M:G11}
|{x_a}|\leq CA^{m}\e^{2m-1}\quad \mbox{as}\ \ a\nearrow a^*,
\end{equation}
where $x_a$ is the unique global maximal point of $|u_a|$ as $a\nearrow a^*$.
Finally, let $a<a^*$  be sufficiently close to $ a^*$ such that $0<A\e<1$. We then conclude from (\ref{5C:G29}) and (\ref{M:G11}) that for any $m\geq 2$,
\[
\|I_a\|_{C^1(\R^2)}\leq C\e^{m+2}\ \ \mbox{in}\ \, \R^2 \ \ \mbox{and}\ \ |{x_a}|\leq C \e^{m-1}\quad \mbox{as}\ \ a\nearrow a^*,
\]
which therefore implies that both $I_a(x)\equiv 0$ and $x_a\equiv 0$ hold as $a\nearrow a^*$.

Since $I_a(x)\equiv 0$ and $x_a\equiv 0$ as $a\nearrow a^*$, we obtain from (\ref{4:1}) and (\ref{4:3}) that $R_a:=\displaystyle\eps_a u_a ( \eps_ax   )e^{ i\Omega\theta_a } $ is a real-valued function and satisfies
\begin{equation}\label{5D:G29}
-\Delta R_a+\big(\eps_a^4|x|^2-\e^2\mu _a-aR_a^2\big)R_a=0\ \ \mbox{in}\,\ \R^2.
\end{equation}
This also gives the absence of vortices for $u_a$ as $a \nearrow a^*$.	
Following (\ref{5D:G29}), the same argument of \cite[Theorem 1.3]{GLW} gives the uniqueness of $R_a$  as $a\nearrow a^*$, which then implies the uniqueness of $u_a$, up to a constant phase. This completes the proof of Theorem \ref{thm1.2*}.\qed

\appendix
\section{Appendix}

In this appendix, we shall establish Lemmas \ref{rev-7} and \ref{AP5} which are used in the proof of Theorem \ref{thm1.2*}.

\begin{lem}\label{rev-7}
For $i=1,2$, suppose that $f_i(x)\in C^1(\R^2)$ satisfies
\begin{equation}\label{rev-9}
	|f_i(x)|,\,|\nabla f_i(x)|\leq Pe^{-Q|x|}\,\ \hbox{in} \,\ \R^2
\end{equation}
for some constants $P>0$ and $Q>0$, and assume that there exists $b_i\in\R$ such that $V_i(x)\in C^1(\R^2)$ satisfies
\begin{equation}\label{rev-10}
	V_i(x)-|\nabla V_i|^2- \frac{b_i^2|x|^2}{4}-\big(|b_1|+|b_2|\big)\geq 2Q^2+\delta\,\ \hbox{in}\,\ \R^2
\end{equation}
for some  constant $\delta>0$. Let $(w_1,w_2)\in C^3(\R^2)\times C^3(\R^2)$ be a real solution of
\begin{equation}\label{rev-2}
\left\{
\begin{aligned}
	-\Delta w_1 +V_1(x)w_1&=b_1\big(x^{\perp}\cdot \nabla w_2\big)+f_1(x)\,\ \hbox{in} \,\ \R^2,\\
	-\Delta w_2+V_2(x)w_2&=b_2\big(x^{\perp}\cdot \nabla w_1\big)+f_2(x)\,\ \hbox{in} \,\ \R^2
\end{aligned}
\right.
\end{equation}
satisfying
\begin{equation}\label{rev-8}
	|w_i(x)|,\ |\nabla w_i(x)|\to 0\,\ \hbox{as}\,\ |x|\to\infty,\ \ i=1,\, 2.
\end{equation}
Then for $i=1,\, 2$, $w_i(x)$ and $\nabla w_i(x)$ satisfy
\begin{equation}\label{rev-11}
	|w_i(x)|,\ |\nabla w_i(x)|\leq C(\delta)Pe^{-Q|x|}\,\ \hbox{in} \,\ \R^2,
	\end{equation}
where the constant $C(\delta)>0$ depends only on $\delta>0$.
\end{lem}

\noindent{\bf Proof.} We first prove that (\ref{rev-11}) holds for $w_i(x)$ with $i=1 $ and $ 2$. Actually, we obtain from (\ref{rev-2}) that
\begin{equation}\label{rev-3}
\left\{
\begin{aligned}
\Big[-\frac 12\Delta+V_1(x)\Big]w^2_1+|\nabla w_1|^2&=b_1\big(x^{\perp}\cdot \nabla w_2\big)w_1+f_1w_1\,\ \hbox{in} \,\ \R^2,\\
\Big[-\frac 12\Delta+V_2(x)\Big]w^2_2+|\nabla w_2|^2&=b_2\big(x^{\perp}\cdot \nabla w_1\big)w_2+f_2w_2\,\ \hbox{in} \,\ \R^2.
\end{aligned}
\right.
\end{equation}
For $i,j=1,2$, we get that
\[
b_i\big(x^{\perp}\cdot \nabla w_j\big)w_i(x)\leq \frac{b_i^2|x|^2|w_i|^2}{4}+|\nabla w_j|^2,\,\ f_i(x)w_i(x)\leq \frac\delta 2 w_i^2+C_1(\delta)f_i^2,
\]
where $\delta >0$ is as in (\ref{rev-10}).
Applying \eqref{rev-9} and \eqref{rev-10}, it then follows from (\ref{rev-3}) that
\begin{equation}\label{rev-4}
(-\Delta+4Q^2+\delta)(w^2_1+w^2_2)\leq C_1(\delta)P^2e^{-2Q|x|}\,\ \hbox{in} \,\ \R^2.
\end{equation}
Applying the comparison principle to \eqref{rev-4}, we thus deduce from  \eqref{rev-8} that
\[
w^2_1+w^2_2\leq C_2(\delta)P^2e^{-2Q|x|}\,\ \hbox{in} \,\ \R^2,
\]
which implies that (\ref{rev-11}) holds for $w_i(x)$ with $i=1 $ and $ 2$.

We next prove that (\ref{rev-11}) holds for $|\nabla w_i|$ with $i=1 $ and $ 2$. For $ i,\,j=1,\, 2,$ denote $\partial_iw_j(x)=\frac{\partial w_j(x)}{\partial x_i}$, and let $\delta_{ij}$ be the Kronecker function satisfying $\delta_{ij}=1$ for $i=j$ and $\delta_{ij}=0$ for $i\neq j$. We then derive from \eqref{rev-2} that $(\partial_iw_1,\partial_iw_2)$ satisfies
\begin{equation}\label{rev-5}
\left\{
\begin{aligned}
&\big[-\Delta +V_1(x)\big]\partial_iw_1+\partial_i V_1(x)w_1=b_1\big[\big(x^{\perp}\cdot \nabla \partial_iw_2\big)-\delta_{i,2}\partial_1w_2+\delta_{i1}\partial_2w_2\big]+\partial_if_1(x),\\
&\big[-\Delta +V_2(x)\big]\partial_iw_2+\partial_i V_2(x)w_2=b_2\big[\big(x^{\perp}\cdot \nabla \partial_iw_1\big)-\delta_{i,2}\partial_1w_1+\delta_{i1}\partial_2w_1\big]+\partial_if_2(x).\\
\end{aligned}
\right.
\end{equation}
The above system further yields that
\begin{equation}\label{rev-6}
\left\{
\begin{aligned}
&\Big[-\frac 12\Delta +V_1(x)\Big]|\nabla w_1|^2+\sum_{i=1}^{2}|\nabla \partial_iw_1|^2+\sum_{i=1}^{2}\partial_i V_1(x)w_1\partial_iw_1\\
&=b_1\Big[\sum_{i=1}^{2}\big(x^{\perp}\cdot \nabla \partial_iw_2\big)\partial_iw_1+\partial_2w_2\partial_1w_1-\partial_1w_2\partial_2w_1\Big]+\sum_{i=1}^{2}
\partial_if_1(x)\partial_iw_1\,\ \hbox{in} \,\ \R^2,\\
&\Big[-\frac 12\Delta +V_2(x)\Big]|\nabla w_2|^2+\sum_{i=1}^{2}|\nabla \partial_iw_2|^2+\sum_{i=1}^{2}\partial_i V_2(x)w_2\partial_iw_2 \\
&=b_2\Big[\sum_{i=1}^{2}\big(x^{\perp}\cdot \nabla \partial_iw_1\big)\partial_iw_2+\partial_2w_1\partial_1w_2-\partial_1w_1\partial_2w_2\Big]+\sum_{i=1}^{2}
\partial_if_2(x)\partial_iw_2\,\ \hbox{in} \,\ \R^2.\\
\end{aligned}
\right.
\end{equation}
Note that for $i,j,l=1,2$ with $j\neq l$,
\[
\partial_i V_j(x)w_j\partial_iw_j\leq  |\partial_i V_j(x)|^2|\partial_iw_j|^2+\frac 14w^2_j,\quad\partial_if_j(x)\partial_iw_j\leq \frac\delta 2 (\partial_iw_j)^2+C_1(\delta)|\partial_if_j(x)|^2,
\]
\[
b_l(x^{\perp}\cdot \nabla \partial_iw_j\big)\partial_iw_l\leq \frac{b_l^2|x|^2|\partial_i w_l|^2}{4}+|\nabla \partial_iw_j|^2,\quad
b_l\partial_1w_j\partial_2w_l\leq \frac 12|b_l|(|\nabla w_1|^2+|\nabla w_2|^2),
\]
where $\delta >0$ is again as in (\ref{rev-10}).
Following above estimates, we infer from \eqref{rev-10} that
\[
\begin{aligned}
\Big(-\Delta +4Q^2+\delta \Big)(|\nabla w_1|^2+|\nabla w_2|^2)\leq(|w_1|^2+|w_2|^2)+2C_1(\delta)(|\nabla f_1|^2+|\nabla f_2|^2)\ \ \hbox{in} \,\ \R^2.
\end{aligned}
\]
Applying \eqref{rev-9} and \eqref{rev-10}, since (\ref{rev-11}) holds for $w_i(x)$ with $i=1 $ and $ 2$, we deduce from the above equation that
\begin{equation}\label{rev-17}
\Big[-\Delta +4Q^2+\delta \Big](|\nabla w_1|^2+|\nabla w_2|^2)\leq C_3(\delta)P^2e^{-2Q|x|}\ \ \hbox{in} \ \ \R^2.
\end{equation}
Applying the comparison principle to \eqref{rev-17}, we derive from \eqref{rev-8}   that
\[
|\nabla w_1|^2+|\nabla w_2|^2\leq C_4(\delta)P^2e^{-2Q|x|}\ \ \hbox{in} \ \ \R^2.
\]
We therefore conclude from above that (\ref{rev-11}) holds for $|\nabla w_i|$ with $i=1 $ and $ 2$, and we are done. \qed

\vskip 0.05truein

\medskip

The rest part of this appendix is to derive some estimates used in the proof of  Proposition \ref{prop4.4}, for which we consider  the following minimization problem
\begin{equation}\label{AP1}
e_a=\inf_{\{u\in \mathbb{H}, \, \|u\|^2_2=1 \} }E_a(u),
\end{equation}
where $\mathbb{H}  :=  \big \{u\in  H^1(\R ^2, \R):\ \int _{\R ^2}
 |x|^2u^2 dx<\infty\big \}$, and $E_a(u)$ is defined by
\begin{equation}\label{AP2}
E_a(u)=\inte \big(|\nabla u|^2+|x|^2u^2\big)dx-\frac{a}{2}\inte  u^4dx,\ \ a>0.
\end{equation}
As stated in Theorem \ref{thmA}, $e_a$ admits minimizers if and only if $0<a <a^*=\|w\|^2_{2}$.
%
Let $\hat {v}_a>0$ be a real minimizer of $e_{a}$ as $a\nearrow a^*$. Then $\hat {v}_a>0$ satisfies the following Euler-Lagrange equation
\begin{equation}\label{6:AP2}
-\Delta \hat {v}_a+|x|^2\hat {v}_a=\beta_a\hat {v}_a+a\hat {v}_a^3\ \ \mbox{in}\,\ \R^2,
\end{equation}
where $\beta_a\in \R$ is the Lagrange multiplier and satisfies
$$\beta_a=e_{a}-\frac{a}{2}\inte \hat {v}_a^4 <0\ \ \mbox{as}\,\ a\nearrow a^*.$$
We also denote
$$ \a=\frac{(a^*-a)^{\frac{1}{4}}}{\lambda}>0.$$
The following lemma gives the estimates of $\beta_a$ and $\hat {v}_a$ as $a\nearrow a^*$:

\begin{lem}\label{AP5}
Let $\hat {v}_a>0$ be a real minimizer of $e_{a}$. Then as $a\nearrow a^*$,
\begin{enumerate}
\item [(i).] $1+\a^2\beta_a =O(\a ^4)$;
\item [(ii).]  $\big|\a\hat {v}_a(\a x)\big|\leq Ce^{-\frac 34|x|},\,\Big|\nabla\big(\a\hat {v}_a(\a x)\big)\Big|\leq Ce^{-\frac 23|x|}$ in $\R^2$;
\item [(iii).] $\max_{i,j=1,2}\Big|\partial_i\partial_j\Big(\a\hat {v}_a(\a x)\Big)\Big|\leq Ce^{-\frac{7}{12}|x|}$\, in  $\R^2$;
\item [(iv).]
$\Big|\a\hat {v}_a(\a x)-\frac{w(x)}{\sqrt{a^*}}\Big|\leq C\a^{4}e^{-\frac{2}{3}|x|}$,\,\, $\Big|\nabla \Big(\a\hat {v}_a(\a x)-\frac{w(x)}{\sqrt{a^*}}\Big)
\Big|\leq C\a^{4}e^{-\frac{1}{2}|x|}$\, in  $\R^2$;
\end{enumerate}
All above constants $C>0$ are independent of $0<a<a^*$.
\end{lem}


{\noindent \bf Proof.} 1. The estimate (i) follows directly from (3.1) and (3.37) in \cite{GLW}.

2. Denote $\bar{v}_a(x)=\a \hat v(\a x)$. It then follows from (\ref{6:AP2}) that $\bar{v}_a$ satisfies
\begin{equation}\label{6:A:2}
-\Delta \bar {v}_a+\a^4|x|^2\bar {v}_a=\a^2\beta_a\bar {v}_a+a\bar {v}_a^3\ \ \mbox{in}\,\ \R^2.
\end{equation}
Similar to  (3.13) in \cite{GWZZ}, one can obtain from (\ref{6:A:2}) that $\|\bar {v}_a\|_{L^{\infty}(\R^2)}\leq C$ and $\bar {v}_a(x)\to 0$ as $|x|\to\infty$ for all $0<a<a^*$. Since the estimate (i) gives that $\a^2\beta_a\to -1$ as $a\nearrow a^*$, we derive from (\ref{6:A:2}) that  as $a\nearrow a^*$,
\begin{equation}\label{6:A:1}
-\Delta \bar {v}_a+\frac{9}{16}\bar {v}_a\leq 0 \ \ \mbox{in}\, \ \R^2/B_R,
\end{equation}
where $R>0$ is large enough.
By the comparison principle, we then deduce from (\ref{6:A:1}) that  as $a\nearrow a^*$,
\[
|\bar {v}_a(x)|\leq Ce^{-\frac 34|x|} \ \ \mbox{in}\, \ \R^2/B_R,
\]
which implies that there exists a sufficiently large constant $C>0$ such that  as $a\nearrow a^*$,
\begin{equation}\label{6:AA:1}
|\bar {v}_a(x)|\leq Ce^{-\frac 34|x|} \ \ \mbox{in}\, \ \R^2.
\end{equation}

Next, we give the estimate of $\nabla \bar {v}_a(x)$ as follows. Denoting $\bar{v}_{a,j}:=\frac{\partial \bar{v}_a}{\partial x_j},$ $j=1, 2$, we then infer from \eqref{6:A:2} that $\bar{v}_{a,j}$ satisfies
\begin{equation}\label{AP13}
-\Delta \bar{v}_{a,j}+\big(\a^4|x|^2-\a^2\beta_a-3a\bar{v}^2_{a}\big)\bar{v}_{a,j}=-2\a^4x_j\bar{v}_a\quad\hbox{in \ $\R^2$.}
\end{equation}
Applying (i) and (\ref{6:AA:1}), we have for $a\nearrow a^*$,
\[
-\a^2\beta_a\to 1,\ |x_j\bar{v}_a|\leq Ce^{-\frac{2}{3}|x|}\ \ \hbox{in}\,\ \R^2/B_R,
\]
where $R>0$ is large enough.
Therefore, by the comparison principle, we deduce from (\ref{AP13}) that  as $a\nearrow a^*$,
\[
|\bar{v}_{a,j}(x)|\leq Ce^{-\frac{2}{3}|x|}\ \ \hbox{in}\ \ \R^2/B_R.
\]
On the other hand, similar to (3.13) in \cite{GWZZ} again,  one can get from (\ref{AP13}) that $\bar{v}_{a,j}(x)$ is bounded uniformly in $\R^2$ for $k$, which further implies that  as $a\nearrow a^*$,
\begin{equation}\label{A:AP13}
|\bar{v}_{a,j}(x)|\leq Ce^{-\frac{2}{3}|x|}\ \ \hbox{in}\,\ \R^2,\ \ j=1,\,2,
\end{equation}
from which we then obtain the desired estimate of $\nabla \bar{v}_{a}$  as $a\nearrow a^*$. This gives the estimate of (ii).

3. By the exponential estimate (\ref{A:AP13}), applying gradient estimates \cite[(3.15)]{GT} to \eqref{AP13} gives that  as $a\nearrow a^*$,
\[
|\nabla \bar {v}_{a,j}(x)|\leq Ce^{-\frac{7}{12}|x|} \ \ \mbox{in\, $\R^2$},
\]
where $\bar{v}_{a,j}:=\frac{\partial \bar{v}_a}{\partial x_j}$ for $j=1, 2$.
Taking $j=1,\,2$, we then obtain that for any $i,\,j=1,\,2$,
\[
|\partial_i\partial_j\bar {v}_{a}(x)|\leq Ce^{-\frac{7}{12}|x|} \ \ \mbox{in\, $\R^2$}
\]
as $a\nearrow a^*$, which therefore gives the estimate of (iii).

4. Denote $\tilde{v}_a(x):=\bar v_a(x)-\frac{w}{\sqrt{a^*}}=\a\hat {v}_a(\a x)-\frac{w}{\sqrt{a^*}}$. We then infer from \cite[Theorem 1.4]{GLW} that as $a\nearrow a^*,$
\begin{equation}\label{AP6}
\|\tilde{v}_a(x)\|_{L^{\infty}(\R^2)}\leq C\a^{4}.
\end{equation}
Moreover, because $\tilde{v}_a(x)$ satisfies the following equation
\[
\begin{aligned}
&\quad-\Delta \tilde{v}_a-\a^2\beta_a\tilde{v}_a-a\big(\bar{v}^2_a+\frac{w\bar{v}_a}{\sqrt{a^*}}+\frac{w^2}{a^*}\big)\tilde{v}_a\\
&=-\a^4|x|^2\bar{v}_a+\big(1+\a^2\beta_a\big)\frac{w}{\sqrt{a^*}}+\frac{(a-a^*)w^3}{(a^*)^{\frac 32}}\ \ \hbox{in}\,\  \R^2,
\end{aligned}
\]
by the comparison principle, we obtain from \eqref{AP6} and (i) that  as $a\nearrow a^*$,
\[
\big|\tilde{v}_a(x)\big|\leq C\a^{4}e^{-\frac{2}{3}|x|}\ \ \hbox{in}\,\  \R^2/B_R,
\]
where $R>0$ is large enough. This further implies that  as $a\nearrow a^*$,
\[
\big|\tilde{v}_a(x)\big|\leq C\a^{4}e^{-\frac{2}{3}|x|}\ \ \hbox{in}\,\  \R^2.
\]
Applying the gradient estimate (3.15) in \cite{GT} and the above exponential decay of $\hat{v}_a(x)$, we finally obtain that
\[
|\nabla\tilde{v}_a(x)|\leq C\a^{4}e^{-\frac{|x|}{2}}\ \ \hbox{in}\,\  \R^2
\]
 as $a\nearrow a^*$, which therefore completes the proof of Lemma \ref{AP5}. \qed

 \vskip 0.16truein
\noindent {\bf Acknowledgements:}
The authors are  very grateful to the referees for many valuable suggestions which lead to the great improvements of the present paper. The authors also thank Prof. Robert Seiringer very much for his helpful discussions and interests on the subject of the present paper.


\begin{thebibliography}{GNN}

\bibitem{Abo} J. R. Abo-Shaeer, C. Raman, J. M. Vogels and W. Ketterle, {\em Observation of vortex lattices in Bose-Einstein condensate}, Science {\bf 292} (2001), 476.








\bibitem{A} A. Aftalion,  Vortices in Bose-Einstein condensates, Progress in Nonlinear Differential Equations and their Applications, 67. Birkh$\ddot{a}$user Boston, Inc., Boston, MA,  2006.





\bibitem{AA} A. Aftalion, S. Alama and L. Bronsard, {\em Giant vortex and the breakdown of strong pinning in a rotating Bose-Einstein condensate}, Arch. Rational Mech. Anal. {\bf 178} (2005), 247--286.





\bibitem{AJ} A. Aftalion, R. L. Jerrard and J. Royo-Letelier, {\em Non-existence of vortices in the small density region of a condensate}, J. Funct. Anal. {\bf 260} (2011), 2387--2406.







\bibitem{Anderson} M. H. Anderson, J. R. Ensher, M. R. Matthews, C. E. Wieman and E. A. Cornell, {\em Observation of Bose-Einstein condensation in a dilute atomic vapor}, Science {\bf 269} (1995), 198--201.



\bibitem{ANS} J. Arbunich, I. Nenciu and C. Sparber, {\em Stability and instability properties of rotating Bose-Einstein condensates}, Lett. Math. Phys. {\bf 109} (2019), 1415--1432.









\bibitem{AS} G. Arioli and A. Szulkin, {\em A semilinear Schr$\ddot{o}$dinger equation in the presence of a magnetic field}, Arch. Ration. Mech. Anal. {\bf 170} (2003), 277--295.



\bibitem{BC} W. Z. Bao and Y. Y. Cai, {\em Ground states of two-component Bose-Einstein condensates with an internal atomic Josephson junction}, East Asia J. Appl. Math. {\bf 1} (2011), 49--81.





\bibitem{B} I. Bloch, J. Dalibard and W. Zwerger, {\em Many-body physics with ultracold gases}, Reviews of Modern Phys. {\bf 80} (2008), 885--964.







\bibitem{Hulet2} C. C. Bradley, C. A. Sackett and R. G. Hulet, {\it Bose-Einstein condensation of lithium: observation of limited condensate number}, Phys. Rev. Lett. {\bf 78} (1997), 985.



\bibitem{Hulet1} C. C. Bradley, C. A. Sackett, J. J. Tollett and R. G. Hulet, {\it Evidence of Bose-Einstein condensation in an atomic gas with attractive interactions}, Phys. Rev. Lett. {\bf 75} (1995), 1687. {\it Erratum} Phys. Rev. Lett. {\bf 79} (1997), 1170.





\bibitem{BO} J. Byeon and Y. Oshita, {\em Uniqueness of standing waves for nonlinear Schr$\ddot{o}$dinger equations}, Proc. Roy. Soc. Edinburgh Sect. A {\bf  138} (2008), 975--987.



\bibitem{BWang} J. Byeon and Z. Q. Wang, {\em Standing waves with a critical frequency for nonlinear Schr$\ddot{o}$dinger equations}, Arch. Rational Mech. Anal. {\bf 165} (2002), 295--316.



\bibitem{Cao} D. M. Cao and  Z. W. Tang, {\em Existence and uniqueness of multi-bump bound states of nonlinear  Schrodinger  equations with electromagnetic fields},  J. Differential Equations  {\bf 222}  (2006), no. 2, 381--424.


\bibitem{CC} L. D. Carr and C. W. Clark, {\em Vortices in attractive Bose-Einstein condensates in two dimensions}, Phys.
Rev. Lett. {\bf 97} (2006), 010403.


\bibitem{C} T. Cazenave, Semilinear Schr$\ddot{o}$dinger equations, Courant Lecture Notes in Mathematics  Vol. 10, Courant Institute of Mathematical Science/AMS, New York, 2003.



\bibitem{CD}  M. Correggi and D. Dimonte, {\em On the third critical speed for rotating Bose-Einstein condensates}, J. Math. Phys. {\bf 57} (2016), 071901.



\bibitem{CP} M. Correggi, F. Pinsker, N. Rougerie and J. Yngvason, {\em Critical rotational speeds for superfluids in homogeneous traps}, J. Math. Phys. {\bf 53} (2012), 095203.



\bibitem{CR} M. Correggi and N. Rougerie, {\em Boundary behavior of the Ginzburg-Landau order parameter in the surface superconductivity regime}, Arch. Rational Mech. Anal. {\bf 219} (2016), 553--606.



\bibitem{CRY} M. Correggi, N. Rougerie and J. Yngvason, {\em The transition to a giant vortex phase in a fast rotating Bose-Einstein condensate}, Comm. Math. Phys. {\bf 303} (2011), 451--508.


\bibitem{D} F. Dalfovo, S. Giorgini, L. P. Pitaevskii and S. Stringari, {\em Theory of Bose-Einstein condensation in trapped gases}, Reviews of Modern Phys. {\bf 71} (1999), 463--512.

\bibitem{EL} M. J. Esteban and P. L. Lions, {\em Stationary solutions of nonlinear Schr$\ddot{o}$dinger equations with an external magnetic field}, Partial differential equations and the calculus of variations, Vol. I,  401--449, Progr. Nonlinear Differential Equations Appl. 1, Birkhuser Boston, Boston, MA, 1989.

\bibitem{F} A. L. Fetter, {\em Rotating trapped Bose-Einstein condensates}, Reviews of Modern Phys. {\bf 81} (2009), 647--691.



\bibitem{GNN} B. Gidas, W. M. Ni and L. Nirenberg, {\em Symmetry of positive solutions of nonlinear elliptic equations in $\R^n$}, Mathematical analysis and applications  Part A, Adv. in Math. Suppl. Stud. Vol. {\bf 7}, Academic Press, New York  (1981), 369--402.


\bibitem{GT} D. Gilbarg and N. S. Trudinger, Elliptic Partial Differential Equations of Second Order, Springer, 1997.


\bibitem{G60} E. P. Gross, {\em Structure of a quantized vortex in boson systems}, Nuovo Cimento {\bf 20} (1961), 454--466.


\bibitem{G63} E. P. Gross, {\em Hydrodynamics of a superfluid condensate}, J. Math. Phys. {\bf 4} (1963), 195--207.

\bibitem{Grossi} M. Grossi, {\em On the number of single-peak solutions of the nonlinear Schr\"odinger equations}, Ann. Inst. H. Poincar\'e Anal. Non Lin\'eaire  {\bf 19} (2002), 261--280.


\bibitem{GLW} Y. J. Guo, C. S. Lin and J. C. Wei, {\em Local uniqueness and refined spike profiles of ground states for two-dimensional attractive Bose-Einstein condensates}, SIAM J. Math. Anal. {\bf 49} (2017), 3671--3715.


\bibitem{GS}  Y. J. Guo and R. Seiringer, {\em On the mass concentration for Bose-Einstein condensates with attractive interactions}, Lett. Math. Phys. {\bf 104} (2014), 141--156.



\bibitem{GWZZ} Y. J. Guo, Z. Q. Wang, X. Y. Zeng and H. S. Zhou,  {\em  Properties of ground states of attractive Gross-Pitaevskii equations with multi-well potentials}, Nonlinearity {\bf 31}  (2018),  957--979.



\bibitem{GZZ} Y. J. Guo, X. Y. Zeng and H. S. Zhou, {\em Energy estimates and symmetry breaking  in attractive Bose-Einstein condensates with ring-shaped potentials}, Ann. Inst. H. Poincar\'e Anal. Non Lin\'eaire {\bf 33} (2016), 809--828.



\bibitem{HL} Q. Han and F. H. Lin,  Elliptic Partial Differential Equations, Courant Lecture Note in Math. 1, Courant Institute of Mathematical Science/AMS, New York, 2011.



\bibitem{HM}  C. Huepe, S. Metens, G. Dewel, P. Borckmans and  M.E. Brachet, {\em Decay rates in attractive Bose-Einstein condensates}, Phys. Rev. Lett. {\bf 82}  (1999), 1616--1619.



\bibitem{IM-1} R. Ignat and V. Millot, {\em The critical velocity for vortex existence in a two-dimensional rotating Bose-Einstein condensate}, J. Funct. Anal. {\bf 233} (2006), 260--306.



\bibitem{IM-2}   R. Ignat and V. Millot, {\em Energy expansion and vortex location for a two-dimensional rotating Bose-Einstein condensate}, Rev. Math. Phys. {\bf 18} (2006), 119--162.





\bibitem{KM} Y. Kagan, A.E. Muryshev and G.V. Shlyapnikov, {\em Collapse and Bose-Einstein condensation in a trapped Bose gas with nagative scattering length}, Phys. Rev. Lett. {\bf 81} (1998), 933--937.




%




\bibitem{KTU} K. Kasamatsu, M. Tsubota and M. Ueda, {\em Giant hole and circular superflow in a fast rotating Bose-Einstein condensate}, Phys. Rev. B {\bf 66} (2002), 053606.







\bibitem{Ku} K. Kurata, {\em Existence and semi-classical limit of the least energy solution to a nonlinear Schr$\ddot{o}$dinger equation with electromagnetic fields}, Nonlinear Anal. {\bf 41} (2000), 763--778.





\bibitem{K} M. K. Kwong, {\em Uniqueness of positive solutions of $\Delta u-u+u^p=0$  in $\R^N$}, Arch. Rational Mech. Anal. {\bf 105} (1989), 243--266.



\bibitem{Lewin} M. Lewin, P. T. Nam and N. Rougerie, {\em Blow-up profile of rotating 2D focusing Bose gases}, Macroscopic Limits of Quantum Systems, Springer Verlag, 2018.



\bibitem{Lieb} E. H. Lieb and M. Loss, Analysis, Graduate Studies in Mathematics Vol. 14. Amer. Math. Soc., Providence, RI, second edition, 2001.



\bibitem{Lieb06} E. H. Lieb and R. Seiringer, {\em Derivation of the Gross-Pitaevskii equation for rotating Bose gases}, Comm. Math. Phys. {\bf 264} (2006), 505--537.





\bibitem {LSS} E. H. Lieb, R. Seiringer, J. P. Solovej and  J. Yngvason, {\em The mathematics of the Bose gas and its condensation}, Oberwolfach Seminars, {\bf 34} Birkh$\ddot{a}$user Verlag, Basel, 2005.



\bibitem {LSY} E. H. Lieb, R. Seiringer and J. Yngvason, {\em Bosons in a trap: A rigorous derivation of the Gross-Pitaevskii energy functional}, Phys. Rev. A {\bf 61} (2000), 043602.



\bibitem {LT} M. Loss and B. Thaller, {\em Optimal heat kernel estimates for Schr\"{o}dinger operators with magnetic field in two dimensions}, Comm. Math. Phys. {\bf 186} (1997), 95--107.


\bibitem {LC} E. Lundh, A. Collin, and K.-A. Suominen, {\em Rotational states of Bose gases with attractive interactions
in anharmonic traps}, Phys. Rev. Lett. {\bf 92} (2004), 070401.




\bibitem{NT} W.-M. Ni and I. Takagi, {\em On the shape of least-energy solutions to a semilinear Neumann problem}, Comm. Pure Appl. Math. {\bf 44} (1991), 819--851.



\bibitem{NT2} W.-M. Ni and I. Takagi, {\em Locating the peaks of least-energy solutions to a semilinear Neumann problem},  Duke Math. J. {\bf 70} (1993), 247--281.





\bibitem{OK} N. Okazawa, {\em An $L^p$ theory for Schr\"{o}dinger operator with nonnegative potentials}, J. Math. Soc. Japan. {\bf 36} (1984), 675--688.





\bibitem {P} L. P. Pitaevskii, {\em Vortex lines in an imperfect Bose gas}, Sov. Phys. JETP. {\bf 13} (1961), 451--454.


\bibitem {RS} M. Reed and B. Simon, Methods of modern mathematical physics. IV. Analysis of operators, Academic Press, New York-London, 1978.



\bibitem{Ro} N. Rougerie, {\em The giant vortex state for a Bose-Einstein condensate in a rotating anharmonic trap: extreme rotation regimes}, J. Math. Pures Appl. {\bf 9} (2011), 296--347.


\bibitem{Hulet3} C. A. Sackett, H. T. C. Stoof and R. G. Hulet, {\em Growth and collapse of a Bose-Einstein condensate with attractive interactions}, Phys. Rev. Lett. {\bf 80} (1998), 2031.



\bibitem{SS} E. Sandier and S. Serfaty, {\em Global minimizers for the Ginzburg-Landau functional below the first critical magnetic field}, Ann. Inst. H. Poincar\' e  Anal. Non Lin\' eaire {\bf 17} (2000), 119--145.


\bibitem{SSbook} E. Sandier and S. Serfaty,  Vortices in the Magnetic Ginzburg-Landau Model, Progress in Nonlinear Differential Equations and their Applications {\bf 70}, Basel: Birkh\'auser, 2007.




\bibitem {SS} S. Secchi and M. Squassina, {\em On the location of spikes for the Schr\"{o}dinger equation with electromagnetic field}, Comm. Contemp. Math. {\bf 7} (2005), 251--268.



\bibitem {S02} R. Seiringer, {\em Gross-Pitaevskii theory of the rotating Bose gas}, Comm. Math. Phys. {\bf 229} (2002), 491--509.






\bibitem {WT} M. Wadati and T. Tsurumi, {\em Critical number of atoms for the magnetically trapped Bose-Einstein condensate with negative s-wave scattering length}, Phys. Lett. A {\bf 247} (1998), 287--293.




\bibitem {Wang} X. F. Wang, {\em On concentration of positive bound states of nonlinear Schr$\ddot{o}$dinger equations}, Comm. Math. Phys. {\bf 153} (1993),  229--244.





\bibitem {W} M. I. Weinstein, {\em Nonlinear Schr$\ddot{o}$dinger equations and sharp interpolations estimates}, Comm. Math. Phys. {\bf 87} (1983), 567--576.

\bibitem {WG} N. K. Wilkin, J. M. F. Gunn and R. A. Smith, {\em Do attractive Bosons condense?}, Phys. Rev. Lett.
{\bf 80} (1998), 2265.

\bibitem {Zhang} J. Zhang, {\em Stability of standing waves for nonlinear Schr\"{o}dinger equations with unbounded potentials}, Z. Angew. Math. Phys. {\bf 51} (2000), 498--503.



\bibitem {Z} J. Zhang, {\em Stability of attractive Bose-Einstein condensates}, J. Stat. Phys. {\bf 101} (2000), 731--746.



\bibitem {ZW} M. W. Zwierlein, J. R. Abo-Shaeer, A. Schirotzek, C. H. Schunck and W. Ketterle, {\em Vortices and superfluidity in a strongly interacting fermi gas}, Nature {\bf 435} (2005), 1047--1051.


\end{thebibliography}
\end{document}